\documentclass[final,onefignum,onetabnum]{siamonline171218}

\usepackage{amsfonts}
\usepackage{graphicx}
\usepackage{epstopdf}
\usepackage{dsfont}
\usepackage{subfig}
\usepackage{amssymb}
\usepackage{algorithmic}
\ifpdf
  \DeclareGraphicsExtensions{.eps,.pdf,.png,.jpg}
\else
  \DeclareGraphicsExtensions{.eps}
\fi

\usepackage{amsmath}
\makeatletter
\renewcommand{\pod}[1]{\allowbreak\mathchoice
  {\if@display \mkern 18mu\else \mkern 8mu\fi (#1)}
  {\if@display \mkern 18mu\else \mkern 8mu\fi (#1)}
  {\mkern4mu(#1)}
  {\mkern4mu(#1)}
}

\usepackage{enumitem}
\setlist[enumerate]{leftmargin=.5in}
\setlist[itemize]{leftmargin=.5in}

\newsiamremark{remark}{Remark}
\newsiamremark{hypothesis}{Hypothesis}
\crefname{hypothesis}{Hypothesis}{Hypotheses}
\newsiamthm{claim}{Claim}
\newsiamthm{conjecture}{Conjecture}

\headers{Growth factors of random butterfly matrices}{J. Peca-Medlin and T. Trogdon}

\title{Growth factors of random butterfly matrices and the stability of avoiding pivoting\thanks{Submitted to the editors \today.
\funding{This work was partially funded by NSF DMS-1916492 (AGEP-GRS Supplement), NSF DMS-1945652 (TT)}}}

\author{John Peca-Medlin\thanks{Department of Mathematics, University of Arizona, Tucson, AZ 
  (\email{johnpeca@math.arizona.edu}).}
\and Thomas Trogdon\thanks{Department of Applied Mathematics, University of Washington, Seattle, WA 
  (\email{trogdon@uw.edu}).}
}

\usepackage{amsopn}

\makeatletter
\newcommand*{\addFileDependency}[1]{
  \typeout{(#1)}
  \@addtofilelist{#1}
  \IfFileExists{#1}{}{\typeout{No file #1.}}
}
\makeatother

\newcommand*{\myexternaldocument}[1]{%
    \externaldocument{#1}%
    \addFileDependency{#1.tex}%
    \addFileDependency{#1.aux}%
}

\def\V#1{{\mathbf #1}}

\def\O{\operatorname{O}}
\def\P{\mathbb P}
\def\E{\mathbb E}
\def\U{\operatorname{U}}
\def\SO{\operatorname{SO}}
\def\SU{\operatorname{SU}}
\def\B{\operatorname{B}}

\def\trace{\operatorname{Tr}}
\def\Haar{\operatorname{Haar}}

\def\Uniform{\operatorname{Uniform}}
\def\Cauchy{\operatorname{Cauchy}}
\def\t{\theta}

\ifpdf
\hypersetup{
  pdftitle={Growth factors of random butterfly matrices},
  pdfauthor={J. Peca-Medlin, T. Trogdon}
}
\fi

\myexternaldocument{ex_supplement}

\begin{document}

\maketitle

\begin{abstract}
Random butterfly matrices were introduced by Parker in 1995 to  remove the need for pivoting when using Gaussian elimination. The growing applications of butterfly matrices have often eclipsed the mathematical understanding of how or why butterfly matrices are able to accomplish these given tasks. To help begin to close this gap using theoretical and numerical approaches, we  explore the impact on the growth factor of preconditioning a linear system by butterfly matrices. These results are compared to other common methods found in randomized numerical linear algebra. In these experiments, we show preconditioning using butterfly matrices has a more significant dampening impact on large growth factors than other common preconditioners and a smaller increase to minimal growth factor systems. Moreover, we are able to determine the full distribution of the growth factors for a subclass of random butterfly matrices. Previous results by Trefethen and Schreiber relating to the distribution of random growth factors were limited to empirical estimates of the first moment for Ginibre matrices.
\end{abstract}




\section{Introduction and background}

Butterfly matrices are a recursively defined subclass of orthogonal matrices that were introduced by D. Stott Parker as a tool to accelerate common methods in computational linear algebra \cite{Pa95}. This recursive structure enables quick matrix-vector multiplication (see \Cref{sec:rbm} for the full definition of butterfly matrices). In particular,  for a dimension $N=2^n$ butterfly matrix $\Omega$ and a vector $\V x \in \mathbb R^N$,  $\Omega \V x$ can be computed in $3Nn$ floating-point operations (FLOPs) rather than the standard $2N^2(1+o(1))$ FLOPs needed to compute $A\V x$ for general $A \in \mathbb R^{N\times N}$. Parker exploited this property to introduce a method to efficiently remove the need for pivoting when using Gaussian elimination (GE) (see \Cref{subsec:ge}). In particular, Parker showed  butterfly matrices can be used to transform any nonsingular matrix to a block nondegenerate matrix\footnote{We say a square matrix is \emph{block degenerate} if a principal minor vanishes, in which case GE without pivoting would be forced to stop at an attempt to divide by 0; otherwise, a  matrix is called \emph{block nondegenerate}.}:

\begin{theorem}[\cite{Pa95}] \label{thm:parker}
If $A$ is a nonsingular matrix and $U,V$ are independent random butterfly matrices, then $UAV^*$ is  block nondegenerate with probability $1-\mathcal O(\epsilon_{machine})$.\footnote{If using exact arithmetic, then $\epsilon_{\operatorname{machine}} = 0$, in which case \Cref{thm:parker} yields an almost sure (a.s.) result.} 
\end{theorem}

Hence, when solving the linear system 
\begin{equation}\label{eq:standard linear}
    A\V x = \V b,
\end{equation} 
using random butterfly matrices $U,V$ one can instead, with very high probability, 
solve the equivalent linear system \begin{equation}\label{eq: precond linear}
    UAV^*\V y = U \V b \quad \text{with} \quad  \V x = V^*\V y
\end{equation}
using GE without pivoting. This requires only $\mathcal O(N^2n)$ FLOPs to transform \eqref{eq:standard linear} into \eqref{eq: precond linear}, not impacting the leading-order complexity of GE (i.e., $\frac23N^3+\mathcal O(N^2)$ FLOPs). 
The costs of moving large amounts of data using pivoting can be substantial on many high performance machines. 
In \cite{baboulin}, numerical experiments running GE with partial pivoting on an dimension 10,000 random matrix using a hybrid CPU/GPU setup resulted in pivoting accounting  for 20 percent of the total computation time.
Pivoting and the communication overhead to coordinate data movements is also  a hindrance to parallel architectures and block algorithms \cite{Pa95}. Removing the need for pivoting would clear  this potential bottleneck to enable faster computations.

The use of butterfly transformations in the machine learning and image processing communities has recently grown in popularity, with a lot of use in Convolutional Neural Network (CNN) architectures \cite{cnn,fast_alg,butterflynet}. Applications remain ahead of the theoretical understanding of the properties of random butterfly matrices. Building off of previous results in \cite{Pa95,Tr19}, this paper aims to fill in some of these missing pieces by outlining new theoretical and experimental results relating to the stability of GE using random butterfly transformations. In particular, we want to further explore the question as to whether removing pivoting is a good idea in practice in spite of the saved computational time. Additionally, in line with \cite{LiLuDo20}, we will explore combining butterfly transformations with iterative refinement. Our main contribution is to measure success of this combination against other random transformations taken from randomized numerical linear algebra, as well as comparing performance metrics against other pivoting schemes used with GE. 


\subsection{Overview of results}

The paper is structured to explore the impact of preconditioning linear systems with random butterfly matrices on the growth factors when using GE with different pivoting strategies.  After giving an overview and background of growth factors and butterfly matrices in the remainder of this section, we  focus on two main linear system models $A \V x = \V b$ of dimension $N = 2^n$. These models are extreme with respect to the growth factor $\rho(A)$ (see \Cref{subsec:gf}), which satisfies $1 \le \rho(A) \le 2^{N-1}$.  Main results will be given in \Cref{sec:gf}, which  includes both theoretical and numerical results comparing butterfly matrices to common preconditioners used in randomized numerical linear algebra. 

We  first consider the \emph{na\"ive model} in \Cref{sec:naive}, where  $A = \V I$, for which $\rho(A) = 1$ is minimized to get a cursory idea of how much the growth factor can be inflated. In this model, we show using numerical experiments that butterfly matrices lead to more modest overall growth and relative errors in minimal growth factor models using GENP, GEPP or GECP. Moreover, comparisons between different pivoting strategies show using GENP with butterfly preconditioners and one step of iterative refinement leads to similar accuracy as using GECP and butterfly preconditioners.

Of note, in \Cref{thm:gf} we give the full distribution of the growth factor for a subclass of random butterfly matrices. This presents, to our knowledge, the first full description of the growth factor for a non-trivial dense linear system. Prior results focus on first-moment empirical estimates or asymptotic bounds (cf. \cite{TrSc90, HiHi20}). The proof of \Cref{thm:gf} will be delayed until \Cref{sec:thm proofs}.

We  then consider a \emph{worst-case model} in \Cref{sec:wc}, where $\rho(A) = 2^{N-1}$ is maximized, as first introduced by Wilkinson. This model is used to test the dampening impact of preconditioning linear systems with large growth factors. Using numerical experiments, we show butterfly matrices as well as other common preconditioners except the Walsh transform have a strong dampening impact on the growth factor when using GENP. Using GEPP, the butterfly matrices lead to Gaussian logarithmic growth factors, which differs significantly from the dampening behavior of Haar orthogonal matrices. Additionally, we show using GENP with one step of iterative refinement leads to similar accuracy as using GECP with butterfly preconditioners.

\subsection{Notation and preliminaries}
\label{sec:prelim}
For convenience, we will use $N = 2^n$ throughout. For $A \in \mathbb R^{n\times m}$,  $A_{ij}$ denotes the entry in the $i$th row and $j$th column of $A$. Let $\V e_i$ denote the standard basis elements of $\mathbb R^n$ and $\V E_{ij} = \V e_i\V e_j^T$, the standard basis elements of $\mathbb R^{n\times m}$.   $\V I$ denotes the identity matrix and $\V 0$  the zero matrix or vector (with the dimensions implicit from context if not stated explicitly).
For $\sigma \in S_n$, the symmetric group on $n$ elements, let $P_\sigma$ denote the orthogonal permutation matrix such that $P_\sigma \V e_i = \V e_{\sigma(i)}$. 
Let $\|\cdot\|_{\max}$ denote the elementwise max norm of a matrix defined by $\|A\|_{\max} = \max_{i,j} |A_{ij}|$ and $\|\cdot\|_\infty$ denote the induced $\ell_\infty$ matrix norm. Note $\|\cdot\|_{\max}$ and $\| \cdot\|_\infty$ are invariant under row or column permutations or multiples of $\pm1$ but are not invariant under general orthogonal transformations.  

$\O(n),\U(n)$ denote the orthogonal and unitary groups of $n\times n$ matrices and  $\SO(n),\SU(n)$ denote the respective special orthogonal and special unitary subgroups; note $\O(n)$ will be used for the orthogonal matrices while $\mathcal O(n)$ is the classical ``big-oh'' notation. Using Gaussian Elimination (GE), we will further emphasize particular pivoting strategies, including GE without pivoting (GENP), GE with partial pivoting (GEPP), GE with rook pivoting (GERP) and GE with complete pivoting (GECP). See \Cref{sec: background} for additional background and \Cref{subsec:ge} for further discussions and results using GE.

We write $X \sim Y$ if $X$ and $Y$ are equal in distribution. Let $X \sim \Uniform([0,2\pi))$ denote $X$ is a uniform random variable with probability density $\frac1{2\pi}\mathds 1_{[0,2\pi)}$ and $Y \sim \Cauchy(1)$ denote $Y$ is a Cauchy random variable with probability density $\frac2\pi \frac1{1 + x^2}$. Let $\operatorname{Gin}(n,m)$ denote the $n\times m$ Ginibre ensemble, consisting of random matrices with independent and identically distributed (iid) standard Gaussian entries. Let $\mathbb S^{n-1} = \{\V x \in \mathbb R^n: \|\V x \|_2 = 1\}$. Recall that if $\V x \sim \operatorname{Gin}(n,1)$, then $\V x/\|\V x\|_2 \sim \Uniform(\mathbb S^{n-1})$. Let $\epsilon_{\operatorname{machine}}$ denote the machine-epsilon, which is the minimal positive number such that $\operatorname{fl}(1+\epsilon_{\operatorname{machine}} ) \ne  1$ using floating-point arithmetic.\footnote{We will use the IEEE standard model for floating-point arithmetic.} 
If using $t$-bit mantissa precision, then $\epsilon_{\operatorname{machine}} = 2^{-t}$. Our later experiments in \Cref{subsec: naive experiments,subsec: WC experiments} will use double precision in MATLAB, which uses a 52-bit mantissa.

{Standard models from randomized numerical linear algebra will be used for comparison in \Cref{subsec: naive experiments,subsec: WC experiments}. These will include the Walsh transformation and Discrete Cosine Transformations (DCT). Sampling for the following experiments will use native (deterministic) MATLAB functions (viz., the Fast Walsh-Hadamard transform \texttt{fwht} and the default Type II Discrete cosine transform \texttt{dct}) applied after an independent row sign transformation chosen uniformly from $\{\pm 1\}^N$. See \cite{St99,Tr11} for an overview of numerical properties of the Walsh and DCT transforms, and \cite{MaTr20} for a thorough survey that provides proper context for use of these transforms and other tools from randomized numerical linear algebra.}

{Additionally, we will utilize left and right invariance properties of the Haar measure on locally compact Hausdorff topological groups, first established by Weil \cite{We40}. For a compact group $G$, this measure can be normalized to yield a probability measure $\Haar(G)$, which inherits the invariance and regularity properties of the original measure. This then provides a uniform means to sample from the compact group, such as $\O(n)$. Stewart provided an outline to sample from $\Haar(\O(n))$ by using $\operatorname{Gin}(n,n)$: if $A \sim \operatorname{Gin}(n,n)$ and $A = QR$ is the $QR$ decomposition of $A$ where $R$ has positive diagonal entries, then $Q \sim \Haar(\O(n))$ \cite{stewart}. Our experiments will employ efficient sampling methods for $\Haar(\O(n))$ that use Gaussian Householder reflectors, in line with the $QR$ factorization of $\operatorname{Gin}(n,n)$ (see \cite{Mezz} for an outline of this method). }

\subsection{Growth factors}
\label{subsec:gf}

The growth factor of a matrix $A$ is a quantity determined by the $LU$ factorization using a fixed pivoting scheme. We will focus on three particular growth factor definitions in this article. See \cite{CoPe07} for an overview on other common definitions found in the literature, along with some explicit properties and relationships comparing different definitions. 

The first growth factor we will consider is related to the max-norm of the associated factors encountered during GE, given by
\begin{equation}
\label{eq: maxnorm gf}
    \rho(A) := \frac{\|L\|_{\max} \cdot {\displaystyle \max_k \|A^{(k)}}\|_{\max}}{\|A\|_{\max}},
\end{equation}
where $A^{(k)}$ denotes the matrix before the $k^{th}$ step in GENP with zeros below the first $k-1$ diagonal entries. This is the classical definition first used by Wilkinson in the 1960s in his error analysis on the backward stability of  GEPP\footnote{Wilkinson's original definition did not include the $\|L\|_{\max}$ factor above, so his results and subsequent studies apply to $\rho(A)/\|L\|_{\max}$; these definitions are consistent when using any pivoting strategy where $\|L\|_{\max} = 1$ always holds, which includes GEPP, GERP, and GECP but not GENP.} \cite{Wi61,Wi65}. This growth factor accounts for the maximal growth that can be encountered throughout the entire GE process. 
Another growth factor is derived from the $\ell_\infty$-induced matrix norm:
\begin{equation}
\label{eq: gf alt}
	\rho_{o}(A) := \frac{\| |L||U| \|_\infty}{\|A \|_\infty}.
\end{equation}
This growth factor only accounts for the final $LU$ factorization of $A$. Our experiments in \Cref{sec:wc} will focus on the following growth factor:
\begin{equation}
\label{eq: l_inf gf}
    \rho_\infty(A) := \frac{\|L\|_\infty \|U\|_\infty}{\|A\|_\infty}.
\end{equation}
Note the trivial bound $\rho_o \le \rho_\infty$, which follows from the submultiplicativity of $\|\cdot\|_\infty$. We will use $\rho_\infty$ for numerical experiments in \Cref{subsec: naive experiments,subsec: WC experiments} since it is computationally simpler to use than $\rho$ and has a more convenient form when applied to our model in \Cref{subsec: WC experiments}.


The growth factor is an important component in controlling the relative error in a computed solution to a linear system using floating-point arithmetic.  If $PAQ=LU$ is the computed $LU$ factorization used to compute the approximate solution $\hat {\V x}$ to the linear system $A\V x = \V b$ for nonsingular $A \in \mathbb R^{n \times n}$, then (cf. \cite{Wi61}) 
\begin{equation}\label{eq:ineq bound max}
    \frac{\|\V x - \hat{\V x}\|_\infty}{\|\V x\|_{\infty}} \le 4 n^2 \epsilon \kappa_\infty (A) \rho(A)
\end{equation}
where $\epsilon = \epsilon_{\operatorname{machine}}$ and 
    $\kappa_\infty(A) = \|A\|_\infty \|A^{-1}\|_\infty$
is the $\ell_\infty$-condition number and (cf. \cite[Section 9.7]{Hi02})
\begin{equation} \label{eq: gf bound inf}
    \frac{\|\V x - \hat{\V x}\|_\infty}{\|\V x\|_\infty} \le \gamma_{3n} \kappa_\infty(A) \rho_o(A)
\end{equation}
where $\displaystyle \gamma_m = \frac{m \epsilon}{1-m \epsilon}$. (Additionally, \cref{eq: gf bound inf} holds when using $\rho_\infty$ instead of $\rho_o$.) As such, error analysis of GE using different pivoting schemes often focuses on the study of the growth factors themselves. Note Wilkinson proved a weak form of backward stability of GEPP by establishing the bound $1 \le \rho \le 2^{n-1}$ along with \eqref{eq:ineq bound max} \cite{Wi61}. Despite this exponential growth factor, GEPP typically has much better performance than this worst-case behavior. Understanding why GEPP still behaves so well remains an outstanding problem in numerical analysis.

\begin{remark}
This paper will focus on the impact of orthogonal transformations on the the growth factor. \eqref{eq:ineq bound max} and \eqref{eq: gf bound inf} yield the corresponding growth factors along with $\kappa_\infty$  control the relative error. Even through $\kappa_\infty$ is not invariant under orthogonal transformations, the impact is relatively moderate. 
For example, we will see $\E \kappa_\infty(B) \approx N^{0.710719}$ for a particular subclass of  $N \times N$ butterfly matrices (cf. \Cref{cor: average cond}).
\end{remark}

\subsubsection{Previous results on growth factors of random matrices}
\label{subsec:prev_gf}

Surprisingly, the literature on random growth factors for dense matrices remains sparse. Interest in growth factors of non-random matrices dates back to Wilkinson's original work establishing the backward stability of GEPP, which established maximal exponential growth factors of dimension $2^{n-1}$ that could lead to a loss of $n-1$ bits of precision \cite{Wi61}. Early focus was on the worst-case behavior of growth factors, which led to very pessimistic views of precision using LU factorizations. In \cite{TrSc90}, Trefethen and Schreiber shifted the view away from the worst-case model to introduce average-case analysis of the stability of GE.  They were interested in why GEPP was successful in practice, with much higher precision than would be expected from the worst-case scenario. To accomplish this, they carried out experiments to compute the growth factors of a variety of random matrices.

Through statistical arguments and numerical experiments, they showed that the average growth factor $\rho$ using GEPP was no larger than $\mathcal O(n)$ for various random matrices with iid entries. Limiting their experiments to matrices of dimension at most $2^{10}$, they showed $\rho$ using GEPP is approximately $n^{2/3}$ while $\rho$ using GECP was approximately $n^{1/2}$ for ensembles with iid entries. They  conjecture further the growth factor should be asymptotically $\mathcal O(n^{1/2})$ for both GEPP and GECP. Additionally, they observe in the iid case that only a few intermediate steps of GE were needed until the remaining entries exhibited approximately normal behavior. Hence, the Ginibre ensemble was a good stand-in for an approximately universal growth factor model for iid matrices. 

In the same paper, Trefethen and Schreiber also experimented with Haar orthogonal matrices and observed the average growth factors were significantly larger than the iid models. This was not too surprising since orthogonal scaled Hadamard matrices have growth factors that are near the largest recorded using GECP \cite{Cryer,EdelmanGECP}. In \cite{HiHi20}, Higham, Higham and Pranesh establish 
\begin{equation}
\label{eq: haar lb}
\rho(A) \gtrsim \frac{n}{4\ln n}
\end{equation}
for $A \sim \Haar(\O(n))$ using any pivoting strategy; i.e., they show asymptotically a {lower bound} growing faster than the iid ensemble growth factors.

One common theme among these preceding works on random growth factors is that the analysis is limited to computing specific parameters (e.g., the first moment) or bounds relating to the growth factors. These only give a small glimpse at the distribution of these random growth factors. The recursive structure of butterfly matrices, as described in the next section, enables us to go beyond these prior limitations. We will provide the full distribution of the growth factors of Haar-butterfly matrices in \Cref{thm:gf}.

\subsection{Random butterfly matrices}
\label{sec:rbm}

Butterfly matrices are a family of recursive orthogonal transformations that were introduced by Parker in 1995 as a means of accelerating common computations in linear algebra \cite{Pa95}. This section will build on top of \cite{Tr19} in further analyzing additional numerical properties of random butterfly matrices.

A dimension $N=2^n$ butterfly matrix $B$ is formed using dimension $N/2 = 2^{n-1}$ symmetric matrices $C,S$ such that $[C,S] = \V 0$ and $C^2+S^2=\V I_{N/2}$
and dimension $N/2 = 2^{n-1}$ butterfly matrices $A_1,A_2$ by forming the product
\begin{equation}
\label{eq:bm_def}
    B = \begin{bmatrix} C & S\\-S & C\end{bmatrix} \begin{bmatrix} A_1 \\&A_2\end{bmatrix} = \begin{bmatrix} C A_1 & S A_2\\-S A_1 & C A_2\end{bmatrix}
\end{equation}
for $n \ge 1$ and starting with $\{1\}$ for $n=0$.
\begin{enumerate}
    \item If $A_1=A_2$ at each recursive step, then the resulting butterfly matrices are called \textit{simple} butterfly matrices.
    \item Let $\B(N)$ and $\B_s(N)$ denote the  $N\times N$ \textit{scalar} butterfly matrices and \textit{simple scalar} butterfly matrices, respectively, formed using scalar matrices  \begin{equation}
    (C_k,S_k)=(\cos\theta_k\V I_{2^{k-1}} ,\sin\theta_k\V I_{2^{k-1}} )    
    \end{equation}
    at each recursive step for $\theta_k \in [0,2\pi)$.
\end{enumerate} 
For example, the  $2\times 2$ butterfly matrices are comprised precisely of $\SO(2)$, the standard clockwise rotation matrices of the form
\begin{equation}
    \label{eq:so2_form}
    B(\t) := \begin{bmatrix}
    \cos\t & \sin\t\\-\sin\t & \cos\t
    \end{bmatrix}.
\end{equation}
Moreover, we write $B(\boldsymbol\t) \in \B_s(N)$ where $\boldsymbol\t \in [0,2\pi)^n$ with $\theta_i$ being the angle introduced in the $i^{th}$ recursive step in constructing $B(\boldsymbol\theta)$. Of particular note, using Kronecker products we can rewrite \eqref{eq:bm_def} in the case of the scalar butterfly matrices as
\begin{equation}
    B=\left(B(\t) \otimes \V I_{N/2} \right) (A_1 \oplus A_2)
\end{equation} for $B \in \B(N)$ with $A_1,A_2 \in \B(N/2)$, and
\begin{equation}\label{eq:sbm_def}
    B(\boldsymbol \t) = \bigotimes_{j=1}^n B(\t_{n-j+1})
\end{equation}
for $B(\boldsymbol \t) \in \B_s(N)$ and $\boldsymbol \t \in [0,2\pi)^n$, where $\t_j \in [0,2\pi)$ is the angle introduced in the $i^{th}$ recursive step to form $B$.


Note $(C,S)$ satisfying the above criteria necessarily satisfy $(C,S) = Q(\Lambda_1,\Lambda_2)Q^*$ for $Q \in \U(N)$ and  diagonal $\Lambda_1,\Lambda_2$ such that $((\Lambda_1)_{ii},(\Lambda_2)_{ii}) = (\cos\theta,\sin\theta) \in S^1$ for some $\theta$ for each $i$. By construction, it follows $\B_s(N) \subset \B(N) \subset \SO(N)$, with equality only for $N=2$. While $\B(N)$ is not closed under multiplication, $\B_s(N)$  comprises an abelian subgroup of $\SO(N)$ where 
\begin{equation}
\label{eq: Haar_b abelian}
B(\boldsymbol\t)B(\boldsymbol\psi) = B(\boldsymbol\t+\boldsymbol\psi) \quad \mbox{and} \quad B(\boldsymbol\t)^{-1} = B(-\boldsymbol\t).    
\end{equation} 
This last result is an immediate consequence of \eqref{eq:sbm_def} and the the mixed-product property of the Kronecker product. See \Cref{sec: background} for additional information and discussion relating to properties of Kronecker products. 

Let $\Sigma$ be a collection of dimension $2^k$ pairs $(C_k,S_k)$ of random symmetric matrices with $[C_k,S_k] = \V0$ and $C_k^2+S_k^2 = \V I_{2^{k}}$. We will write $\B(N,\Sigma)$ and $\B_s(N,\Sigma)$ to denote the ensembles of random butterfly matrices and random simple butterfly matrices formed by  independently sampling $(C,S)$ from $\Sigma$ at each recursive step. Let
\begin{align}
    \Sigma_S &= \{(\cos\theta^{(k)}\V I_{2^{k-1}} ,\sin\theta^{(k)}\V I_{2^{k-1}} ): \theta^{(k)} \mbox{ iid }  \operatorname{Uniform}([0,2\pi), k\ge 1\}
\quad \mbox{and}\\
    \Sigma_D &= \{ \bigoplus_{j=1}^{2^{k-1}}(\cos\t_j^{(k)},\sin\t_j^{(k)}): \theta^{(k)}_j \mbox{ iid }  \operatorname{Uniform}([0,2\pi),  k\ge 1\}.
\end{align}
A large focus for the remainder of this paper is on the \emph{Haar-butterfly matrices}, $\B_s(N,\Sigma_S)$, while numerical experiments in \Cref{sec:naive,sec:wc} will also use the other random scalar butterfly ensemble, $\B(N,\Sigma_S)$ along with the random diagonal butterfly ensembles, $\B(N,\Sigma_D)$ and $\B_s(N,\Sigma_D)$. 

The name of $\B_s(N,\Sigma_S)$ is purposefully suggestive due to the result:

\begin{proposition}[\cite{Tr19}]
\label{prop:haar_butterfly}
    $\B_s(N,\Sigma_S) \sim \Haar(\B_s(N))$.
\end{proposition}
Since $\B_s(N)$ is a closed subgroup of $\SO(N)$, then it has a Haar measure that can be used to uniformly sample a matrix from this ensemble; so this explicit construction of random butterfly matrices  provides a method to directly sample from this distribution. 



A particularly useful result that follows from \eqref{eq:sbm_def} and the mixed-product property is that matrix factorizations of $B(\t) \in \SO(2)$ translate directly to matrix factorizations of $B(\boldsymbol \t) \in \B_s(N)$. For example, expanding a result in \cite{Tr19} (see Lemma 2.4), one can easily compute the eigenvalue decomposition of $B(\t)$ as $B(\t) = U\Lambda_\t U^*$ for
\begin{equation}
    U = \frac1{\sqrt 2} \begin{bmatrix}
    1 & 1\\i & -i
    \end{bmatrix} \quad \mbox{and} \quad \Lambda_\t = \begin{bmatrix}
    e^{i\t}\\&e^{-i\t}
    \end{bmatrix}.
\end{equation}
It follows then an eigenvalue decomposition of $B(\boldsymbol\t)$ is given by $B(\boldsymbol\t) = U_n\Lambda_{\boldsymbol\t} U_n^*$ for
\begin{equation}\label{eq: rbm eigendecomp}
    U_n = \bigotimes^n U \quad \mbox{and} \quad \Lambda_{\boldsymbol\t} = \bigotimes_{j=1}^n \Lambda_{\t_{n-j+1}}.
\end{equation}
Note in particular $U_n$ is independent of $\boldsymbol \t$. This is natural since $\B_s(N)$ is an abelian group of orthogonal matrices
and hence can be simultaneously diagonalized, so all such matrices must share the same eigenvectors. This holds also for $\B(\boldsymbol \t) \sim \B_s(N,\Sigma_S)$. In this case, the stochasticity in the model is confined to the eigenvalues, which are necessarily of the form $e^{i \V v^T \boldsymbol \t}$ for each $\V v \in \{\pm 1\}^n$ by \eqref{eq: rbm eigendecomp}. Moreover, each eigenvalue is uniformly distributed on the unit circle (i.e., the one-point correlation of $\B(\boldsymbol \t)$), see \Cref{lemma:simple}. See \cite{Tr19} for other discussions and results relating to the spectral properties of butterfly matrices.

\section{The impact of random butterfly matrices on growth factors}
\label{sec:gf}

This section will focus on the impact on $\rho$ and $\rho_\infty$ by preconditioning linear systems using random transformations, as outlined in \eqref{eq: precond linear}. We will consider both 1-sided (i.e., $V = \V I$) and 2-sided preconditioning. The 1-sided model has been studied in \cite{baboulin,Tr19}, while the 2-sided model is consistent with the block nondegenerate transformation as originally outlined by Parker \cite{Pa95}. We will use two models, which are both extreme cases with respect to $1 \le \rho \le 2^{N-1}$. 
\begin{enumerate}
    \item \emph{Na\"ive model}: using $A = \V I$, with $\rho(A) = 1$, and
    \item \emph{Worst-case model}: using $A$ such that $\rho(A) = 2^{N-1}$.
\end{enumerate}

The na\"ive model is so named since it would be unnecessary to use any method to solve the trivial linear system $\V I\V x = \V b$. The main motivation to include this is as a test model to look at how much one can mess up a linear system by using preconditioning. For this model, we will only consider the 1-sided preconditioning case so that this  provides a scenario to directly study the growth factors of these random transformations. {Random growth factors for fixed matrix ensembles have been studied previously (e.g., \cite{HiHi20,TrSc90}). We are interested in continuing that line of study so that we can compare the results for our particular random ensembles to these previous studies.}


The worst-case model goes to the other extreme, using a construction first due to Wilkinson  \cite{Wi61}. Wilkinson provided an explicit matrix $A$ such that $\rho(A)=2^{N-1}$ satisfies the upper bound $\rho \le 2^{N-1}$. This model will then provide a sufficient scenario to study the dampening capacity of randomized preconditioning on the growth factors. For this, we will only consider the 2-sided random transformations through numerical experiments. 

\begin{remark}
{We also ran experiments for 2-sided preconditioning using the identity matrix, $U\V I V^* = UV^*$, (in line with Parker's motivation), as well as 1-sided preconditioning using the Wilkinson matrix, $\Omega A$. The results mostly mirrored those for the 1-sided na\"ive and 2-sided worst case scenarios and would not add significantly to the discussion;  so we will refrain from including these additional trials. Of note for the two Haar distributions (i.e., Haar-butterfly and Haar orthogonal) that will be used in the experiments, the invariance properties of the Haar measure then yield $UV^*\sim U$ for $U,V$ both iid Haar matrices, so that the 1-sided and 2-sided na\"ive experiments are equivalent for these specific ensembles. In particular, all of the results in \Cref{subsec:gf_haar_b} equally apply to the 1-sided and 2-sided na\"ive models using Haar-butterfly matrices.} 
\end{remark}

\begin{remark}
While both models represent extreme cases, most practical applications would be far from either model. Additionally, some applications will not be impacted by orthogonal transformations. For example, if $G \sim \operatorname{Gin}(N,N)$ and $U \in \O(N)$ then $UG \sim \operatorname{Gin}(N,N)$. Hence, $BG \sim G$ for \textit{any} random butterfly matrix $B$. Similarly, $AG \sim G$ for $A \sim \Haar(\O(N))$. Other iid models have empirically been shown to behave similarly to the Ginibre case (cf. \cite{TrSc90}). Hence, this suggests using random butterfly matrices as preconditioners in such a setting would have marginal benefits.
\end{remark}

\subsection{Na\"ive model}
\label{sec:naive}

For this model, we  study the growth factors of the random preconditioners directly. The Haar-butterfly models have the rare benefit that the distribution of the associated growth factors for these matrices can be computed explicitly. This will be established in \Cref{subsec:gf_haar_b}. 
\Cref{subsec: naive experiments}  uses numerical experiments to determine approximate distributions of growth factors for other random matrix models.

\subsubsection{Growth factors of 
Haar-butterfly matrices}
\label{subsec:gf_haar_b}

Computations using simple scalar butterfly or Haar-butterfly matrices are made significantly more tractable by combining \eqref{eq:sbm_def} along with \Cref{lemma: norm_mult,lemma: kron_factor}. \Cref{sec:thm proofs} will contain the proofs of the outstanding technical details. 

\begin{theorem}
    \label{thm:gf}
%

    Let $B \sim \B_s(N,\Sigma_S)$ and $X_j \sim \Cauchy(1)$ be iid for $j \ge 1$.  Let $Y_j = |X_j|$ if using GENP and $Y_j = |X_j| \mid |X_j| \le 1$ if using GEPP or GERP.  Then $B$ has an $LU$ factorization using GENP a.s., $B$ has unique factors using GEPP that also determine the GERP factorization, and 
    \begin{align}
        \rho(B) &\sim \prod_{j=1}^n (1 + Y_j^2) \label{eq: thm maxnorm gf}\\
        \rho_o(B) & \sim \prod_{j=1}^n \left( 1 + \frac{2Y_j^2}{1 + Y_j}\right)\\
        \rho_\infty(B) &\sim \prod_{j=1}^n (1 + \max(Y_j,Y_j^2)),
    \end{align}
    where $1 \le \rho(B) \le \rho_\infty(B)$ and $1 \le \rho_o(B) \le \rho_\infty(B)$ when using GENP and $1 \le \rho(B) \le \rho_o(B) \le \rho_\infty(B) \le N$ when using GEPP or GERP. 
    
    Moreover, the GEPP max-norm growth factor of $B \sim \B_s(N,\Sigma_S)$ is minimal among all GE pivoting strategies, i.e., if $P,Q$ are permutation matrices such that $PBQ$ has a GENP factorization a.s., then $\rho^{\operatorname{GENP}}(PBQ) \ge \rho^{\operatorname{GEPP}}(B)$.
%
\end{theorem}
\begin{remark}
When using pivoting, then  $\max(Y_j,Y_j^2) = Y_j$ since $0\le Y_j \le 1$ so that  $\rho_\infty(B) \sim \prod_{j=1}^n (1 + Y_j)$. 
\end{remark}

\begin{proof}[Proof sketch of \cref{thm:gf}] {The main technical steps needed to prove \Cref{thm:gf} involve determining the explicit GENP and GEPP factorizations of a Haar-butterfly matrix along with the intermediate GE factors (see \Cref{prop: ge factors,lemma:Bk}). This is accomplished making heavy use of the Kronecker product structure of simple scalar butterfly matrices along with the mixed-product property. This structure is further exploited to show GEPP and GERP lead to the same $LU$ factorization (see \Cref{cor: gerp}). (This argument does not extend to the later Worst Case model (see \Cref{sec:wc}), which loses the Kronecker product structure.) Most of the technical expense is spent on establishing the maximal growth encountered during GENP, GEPP or GERP  occurs in the final GE step, i.e.,
    \begin{equation}
    \label{eq: max U}
        \max_k \|B^{(k)}\|_{\max} = \|U\|_{\max}
    \end{equation}
(see \Cref{prop:max Uk}). The remainder follows directly from the mixed-product property and the multiplicative property for Kronecker products with respect to $\|\cdot\|_{\max}$ and $\|\cdot\|_\infty$ (by \Cref{lemma: kron_factor,lemma: norm_mult}). The final statement that the GEPP growth factor minimizes the Haar-butterfly growth factor using any pivoting strategy follows from \Cref{thm: lower bd gf} (first found in \cite{HiHi89}) along with the fact $\|B\|_{\max} = \|B^{-1}\|_{\max} = \|U\|_{\max}^{-1}$ for $PB=LU$ the GEPP factorization of $B$. These details are filled out in \Cref{sec:thm proofs}.}
\end{proof}

\begin{remark}
    Explicit distributions can similarly be determined as in \Cref{thm:gf} for other growth factors depending only on the final $LU$ factorization and using a norm that is multiplicative with respect to Kronecker products (e.g., max-norm, induced norms, or $p$-Schatten norms). Additionally, other simple random scalar butterfly models using different distributions on the input angles would similarly result in a multiplicative distribution, but with the Cauchy terms changed to reflect the different angle distributions.
\end{remark}



Using GENP,  $\rho(B)$, $\rho_o(B)$ and $\rho_\infty(B)$ have no finite moments of any order since the absolute Cauchy has no finite moments. However, since $\rho(B)$, $\rho_o(B)$  and $\rho_\infty(B)$ are bounded when using pivoting,  we can calculate the average growth factors \textit{exactly} rather than being restricted to empirical estimates (as in \cite{TrSc90}):
\begin{corollary}
\label{cor: exp gf}
    Let $B \sim \B_s(N,\Sigma_S)$. Using GEPP or GERP, then $\E \rho(B) = N^\alpha$, $\E \rho_o(B) = N^\beta$ and $\E\rho_\infty(B) = N^\gamma$
    for $\alpha = \log_2(\frac4\pi) \approx 0.34850387$, $\beta = \log_2(\frac{6\log 2}\pi) \approx 0.404699998$ and  $\gamma = \log_2(1 + \frac{\log 4}\pi) \approx 0.52734183$. Moreover,  $\E \rho(B)^k = N^{\alpha_k}$ for $\alpha_k = \log_2(\frac2\pi | \operatorname{Im} {B}_2(k,\frac12)|)$.\footnote{Here $B_x(a,b) = \int_0^x t^{a-1}(1-t)^{b-1} \operatorname{d}t$ denotes the incomplete beta function. This arises via $I_k = \frac\pi4\E \rho(B(\t))^k = \int_0^{\pi/4} \sec^{2k}x \operatorname{d}x = \frac12 \int_1^2 \frac{u^{k-1}}{\sqrt{u-1}} \operatorname{d}u = -\frac12 \operatorname{Im} B_2(k,\frac12)$ using a $u = \sec^2x$ substitution. Alternatively, the recurrence $I_k = \frac{2^{k-1}}{2k-1} + \frac{2k-2}{2k-1} I_{k-1}$ with $I_1=1$ enables easy integer moment computations.}
\end{corollary}

We can now relate  the growth factors of Haar-butterfly matrices directly to the growth factors of other random ensembles of matrices studied in \cite{HiHi20,TrSc90}. Using only GEPP, \Cref{thm:gf,cor: exp gf} show $\rho(B)$ is sublinear and $\E \rho(B) \approx N^{0.34850387}$. Of note, this is smaller than the average-case growth factors of iid ensembles using both GEPP {and even GECP}, where Trefethen and Schreiber indicated these were, respectively, about $N^{2/3}$ using GEPP and $N^{1/2}$ using GECP (and asymptotically $\mathcal O(N^{1/2})$). 

Another application using \Cref{cor: exp gf}  and Markov's inequality is if $B \sim \B_s(N,\Sigma_S)$ then using GEPP we have
\begin{equation}
\label{eq: markov}
    \P\left(\rho(B) \ge \frac{N}{4 \ln N}\right) \le 4n \ln 2 N^{\alpha - 1} = N^{\alpha - 1 + o(1)} = o(1),
\end{equation}
since $\alpha < 1$.\footnote{Note \cref{eq: markov} holds if we replace $\frac{N}{4n}$ with $N^{c + o(1)}$ for any $c>\alpha$.} This compares to  $\rho(A) \gtrsim \frac{N}{4n}$ for $A \sim \Haar(\O(N))$ using GEPP \cite{HiHi20}. So for sufficiently large $n$, Haar-butterfly matrices have GEPP growth factors that are strictly less than the growth factors of Haar orthogonal matrices with high probability.


Although $|X|$ has no finite moments of any order for $X \sim \Cauchy(1)$, $\ln(1+X^2)$ has finite moments of any order: for example, $\E \ln(1+X^2) = \ln 4$ and $\E\ln(1+X^2)^2=\frac{\pi^2}3+(\ln 4)^2$ (so that $\operatorname{Var}\ln(1+X^2)=\frac{\pi^2}3$).  Since $\ln \rho(B)$ (and $\ln\rho_o(B)$ and $\ln \rho_\infty(B)$) are then \textit{sums} of iid terms, each of which is dominated by $1+\ln(1+X^2)$,  the Central Limit Theorem provides a method to analyze the growth factors for sufficiently large Haar-butterfly matrices. 
\begin{corollary}[Haar-butterfly CLT]
\label{cor:clt}
Let $B \sim \B_s(N,\Sigma_S)$ and $Z \sim N(0,1)$. Using GENP, GEPP or GERP, then for any $t \in \mathbb R$
\begin{equation}
    \lim_{n\to \infty}\P\left(\frac{\ln \rho(B) - n\mu}{\sqrt{n}\sigma} \le t\right) = \P( Z \le t) 
\end{equation}
where $\mu = \ln 4$, $\sigma^2 = \frac{\pi^2}3$ 
when using GENP and $\mu = \ln4-\frac{4G}\pi$, $\sigma^2 = \frac7{12}\pi^2 +(\ln2)^2 + \frac{4G \mu}{\pi}-\frac{16}\pi \operatorname{Im}(\operatorname{Li}_3(1+i))$ 
when using GEPP or GERP.\footnote{We are using Catalan's constant 
\begin{equation}
    G := \sum_{n\ge 0} \frac{(-1)^n}{(2n+1)^2} \approx 0.91596559,
\end{equation} 
along with the trilogarithmic function $\operatorname{Li}_3$, where
\begin{equation}
    \operatorname{Li}_s(z) := \sum_{k\ge 1} \frac{z^n}{k^s}.
\end{equation}
}
\end{corollary}
\begin{remark}
{\Cref{cor:clt} shows that $(\rho(B) e^{-n\mu})^{1/(\sqrt n \sigma)}$ converges in distribution to a standard lognormal random variable, $e^Z$ where $Z \sim N(0,1)$. Analogous results hold  when using $\rho_\infty$ or $\rho_o$.}
\end{remark}

GENP growth factors often are not studied extensively, which is not too surprising given their instability. \Cref{cor:clt} can be used to better understand typical behavior of GENP for Haar-butterfly matrices. Since $\rho(B)$ does not have a finite mean, the \emph{median}, $M_n$, would then be a more desirable statistic to gauge behavior of $\rho(B)$. For $n$ sufficiently large,  \Cref{cor:clt} guarantees
\begin{equation}
    \P(Z \le 0) = \frac12=\P(\rho(B) \le M_n) 
    \approx \P\left(Z \le \frac{\ln M_n - n\mu}{\sqrt{n} \sigma}\right).
\end{equation}
Since the convergence of the distribution functions in \Cref{cor:clt} is uniform, then $\ln M_n - n \mu = o( n^{1/2})$. Using GENP we have $\mu = \ln 4$, so that $e^{n\mu} = N^{\mu/\ln 2} = N^2$. This can be used as an estimator for the median of $\rho(B)$:
\begin{corollary}\label{cor:median}
Let $B \sim \B_s(N,\Sigma_S)$, $\rho(B)$ be the GENP growth factor, and $M_n$ the median of $\rho(B)$. Then $M_n = N^{2+o(n^{-1/2})}$.
\end{corollary}
Hence, for $n$ sufficiently large, then $\rho(B)$ using GENP is approximately centered around $N^2$ even through $\E \rho(B) = \infty$ for any $n$.

\begin{remark}
Note $p(n)N^2$  matches the form of $M_n$ from \Cref{cor:median} for any polynomial $p$. Experiments suggest $\ln M_n - n\mu = \mathcal O(1)$ (if not $o(1)$). For instance, since \Cref{thm:gf} enables us to sample  growth factors for arbitrarily large Haar-butterfly matrices, we ran $10^6$ experiments for $n=2^{18}$ (corresponding to Haar-butterfly matrices of dimension $N \approx 1.611 \cdot 10^{78913}$). The sample median $\widehat M_n$ satisfied $\ln \widehat M_n - n\mu = -2.084063666115981$ while $\ln \widehat M_n$ and $n \mu$ are each approximately $3.63 \cdot 10^5$. This suggests  \Cref{cor:median} could potentially be improved to $M_n = N^{2 + \mathcal O(n^{-1})}$ or $N^{2+o(n^{-1})}$.
\end{remark}

Using pivoting, we have $\E \ln \rho(B)  = \mu = \ln 4 - \frac{4G}{\pi}$, so $\rho(B)$ is approximately centered around $N^{2-\frac{4G}{\pi \ln 2}} \approx  N^{0.31746612}$. This is not too far off from $\E \rho(B) \approx N^{0.34850387}$.\footnote{Jensen's inequality guarantees the median estimate $2^{\E \log_2 \rho(B)}$ is smaller than $\E \rho(B)$; however, since $\ln \rho(B)$ will have heavier right tails, then this median estimate should be an overestimate of the actual median of $\rho(B)$.} Similarly, \Cref{cor:clt} can be used to estimate quantiles for $\rho(B)$ (as well as for $\rho_o(B)$ and $\rho_\infty(B)$). 

Additionally, explicit results relating to $\kappa_\infty(B)$ can be computed for $B \sim \B_s(N,\Sigma_S)$:

\begin{theorem}
\label{prop: cond infty}

Let $B \sim \B_s(N,\Sigma_S)$ and $Y_j$ be iid $\operatorname{Arcsine}(0,1)$ for $j \ge 1$. Then
\begin{equation}
    \kappa_\infty(B) \sim \prod_{j=1}^n (1 + \sqrt{Y_j})
\end{equation}
with $1 \le \kappa_\infty(B) \le N$.
\end{theorem}

In the context of \Cref{thm:gf}, then $\sqrt{Y_j} = \frac{2|X_j|}{1 + X_j^2}$ for $X_j \sim \Cauchy(1)$ iid (see \Cref{lemma:gf_2x2,lemma: cond_2x2}).\footnote{Note if using GENP then $\kappa_\infty(B)\rho(B) \sim  \prod_{j=1}^n (1+|X_j|)^2$ for $X_j \sim \Cauchy(1)$.} 
Similarly, explicit average condition numbers as well as the average product of the condition number and growth factor when using GEPP or GERP (since $\kappa_\infty \ge 1$ this product is still not integrable in the GENP case), as found in \eqref{eq:ineq bound max} and \eqref{eq: gf bound inf}, can be computed:

\begin{corollary}
\label{cor: average cond}
Let $B \sim \B_s(N,\Sigma_S)$. Then $\E \kappa_\infty (B) = N^\xi$
for $\xi = \log_2(1 + \frac2\pi) \approx 0.71071919$. Using GEPP {or GERP}, then $\E \kappa_\infty(B) \rho(B) = N^{1+\zeta}$, $\E \kappa_\infty(B) \rho_o(B) = N^{1+\psi}$ and $\E \kappa_\infty(B) \rho_\infty(B) = N^{1+\phi}$
for $\zeta = \log_2(\frac2\pi(1+\log 2)) \approx 0.10821126$, $\psi = \log_2(1+\frac{\log 4-1}\pi) \approx 0.16730823$, and $\phi = \log_2(1 + \frac{\log 2}\pi) \approx 0.28763257$.
\end{corollary}

\noindent One consequence of \Cref{cor: average cond} is that using GEPP (or GERP) to solve the na\"ive model preconditioned by a Haar-butterfly matrix $B \sim \B_s(N,\Sigma_S)$ using double precision (i.e., $\epsilon = 2^{-52}$) generates a relative error that has an average upper bound from \eqref{eq:ineq bound max} of
\begin{equation}
    4N^2 \E \kappa_{\infty}(B)\rho(B) \epsilon = 4N^{3+\zeta} \epsilon = 2^{(3+\zeta)n-50}.
\end{equation}
Suppose one wanted to guarantee precision up to about $10^{-6}$. Since $(3+\zeta)n-50 < -20$ for $n<\frac{30}{3+\zeta} \approx 9.65185359$, the average $\ell_\infty$-relative error preconditioning the trivial linear system by a Haar-butterfly matrix using double precision can do no worse than $2^{-20}<10^{-6}$ for $N$ up to $2^{9} = 512$. Quad precision (i.e., $112-$bit precision) would maintain this bound for $N$ up $2^{28} =  2.68435456\cdot 10^8$. So computed solutions using butterfly preconditioners remain stable even for relatively large butterfly matrices.\footnote{Repeating the argument for $\kappa_\infty(B)\rho(B)$ using GENP as used in \Cref{cor:clt,cor:median} to find the approximate median yields about half the time GENP maintains $\ell_\infty$-relative error precision up to 1 for $N$ at most $2^{10} = 1028$.}

\begin{figure}[htp] 
\centering
\subfloat[$\rho^{\operatorname{GENP}}(B)$]{%
    \includegraphics[width=0.48\textwidth]{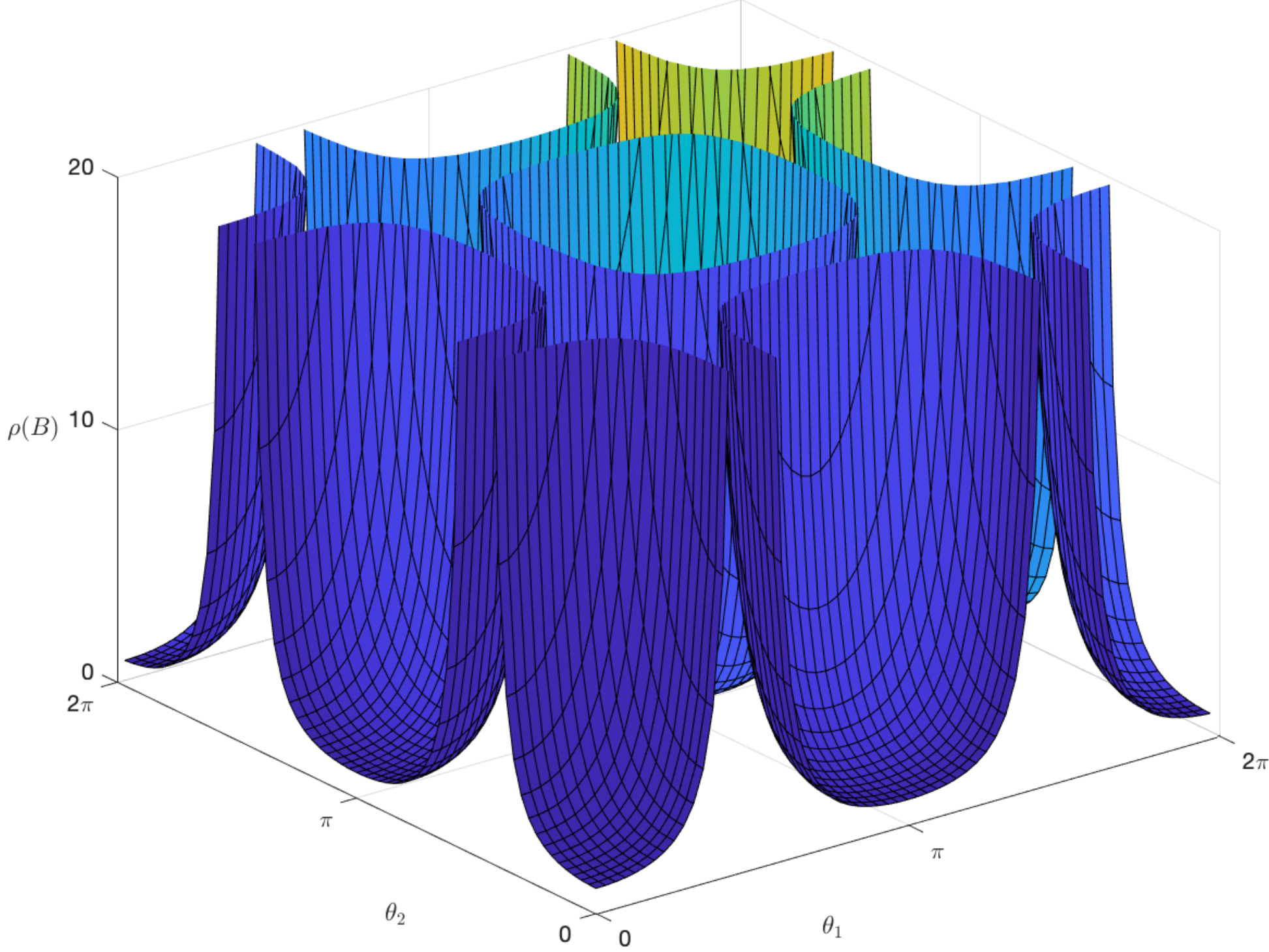}%
    \label{fig:4x4_gf_maxnorm_np}%
    }%
\quad %
\subfloat[$\rho^{\operatorname{GEPP}}(B)$]{%
    \includegraphics[width=0.48\textwidth]{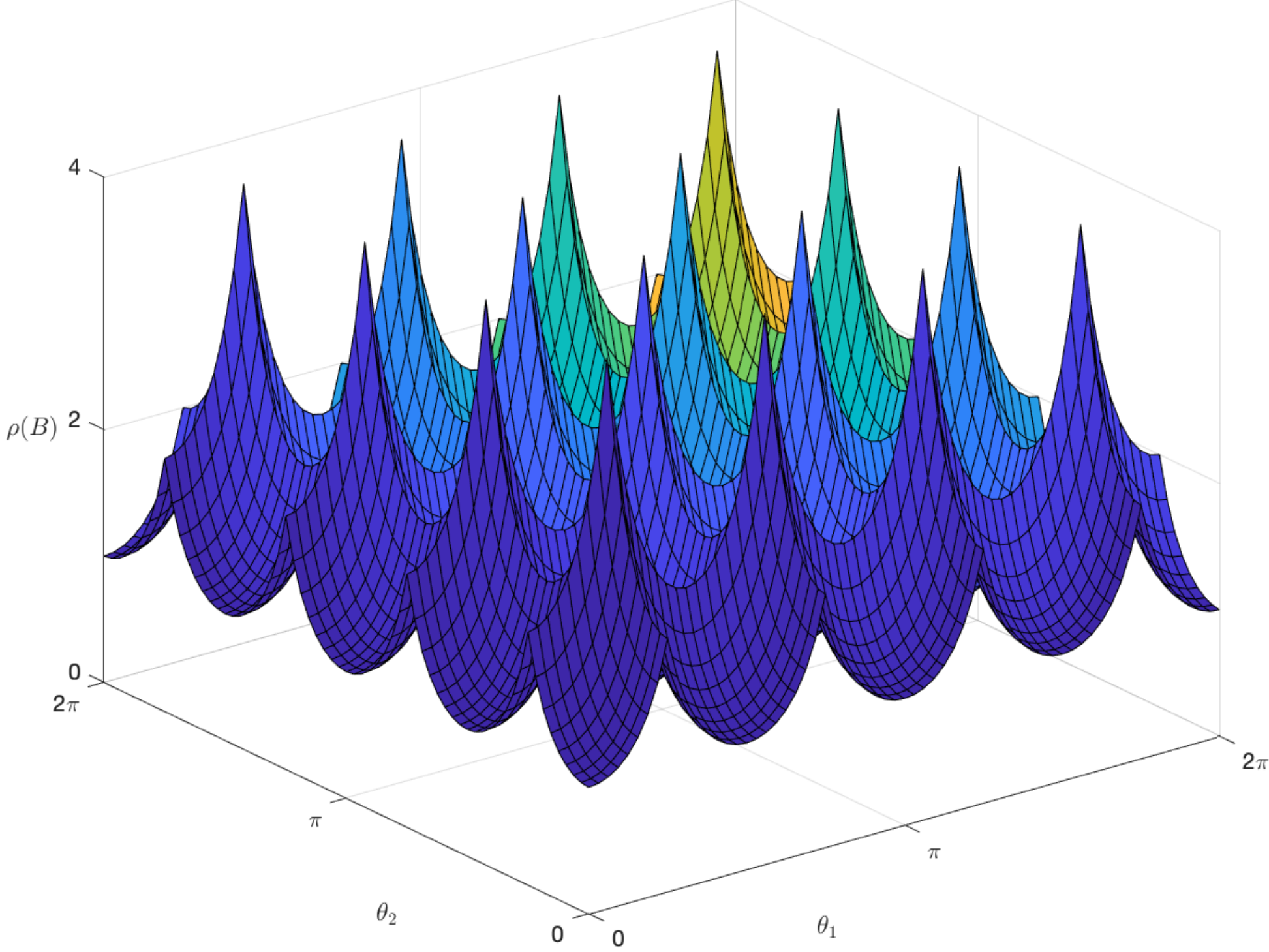}%
    \label{fig:4x4_gf_maxnorm_pp}%
    }%
\\
\subfloat[$\rho_o^{\operatorname{GENP}}(B)$]{%
    \includegraphics[width=0.48\textwidth]{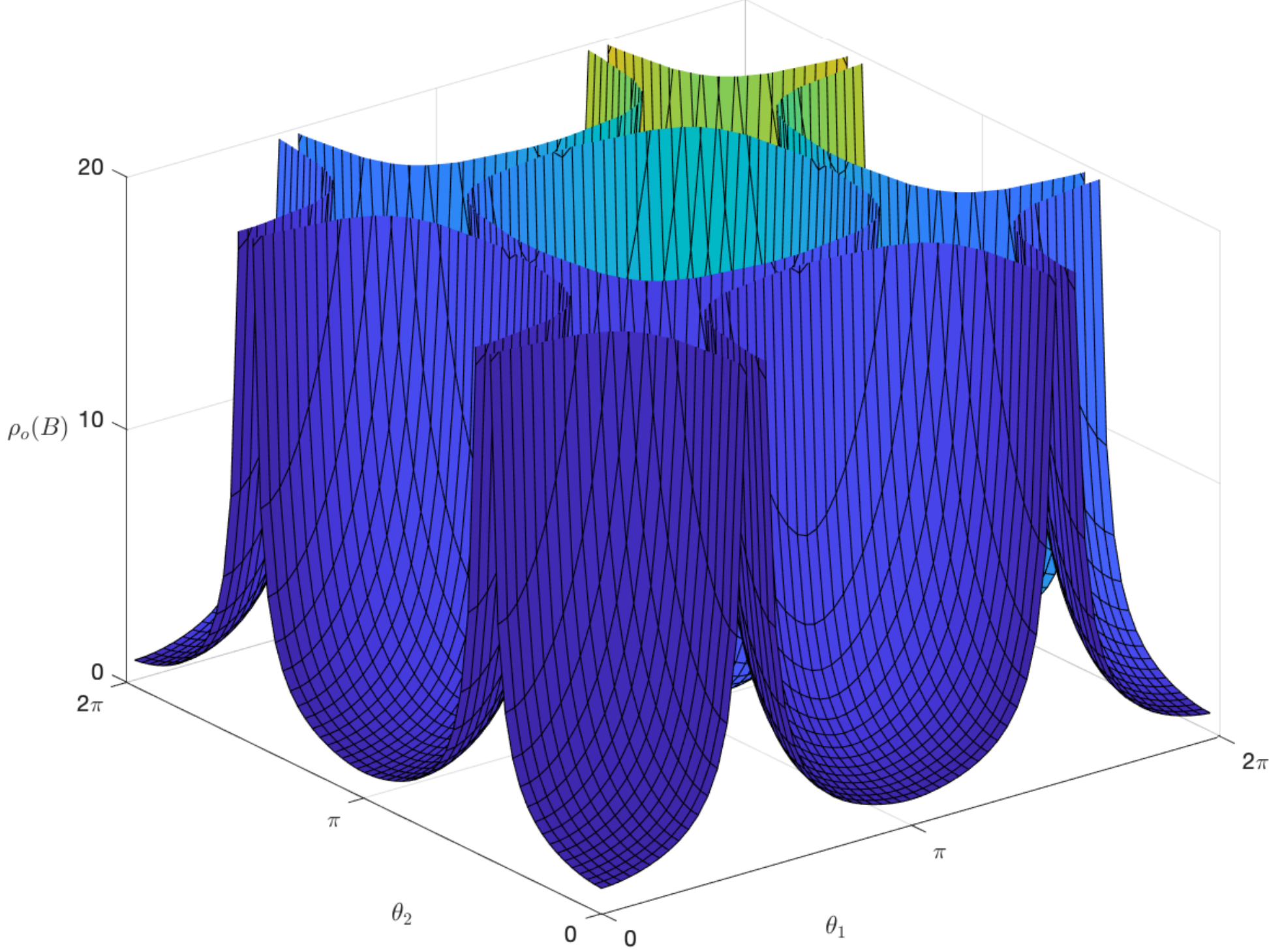}%
    \label{fig:4x4_gf_alt_np}%
    }%
\quad %
\subfloat[$\rho_o^{\operatorname{GEPP}}(B)$]{%
    \includegraphics[width=0.48\textwidth]{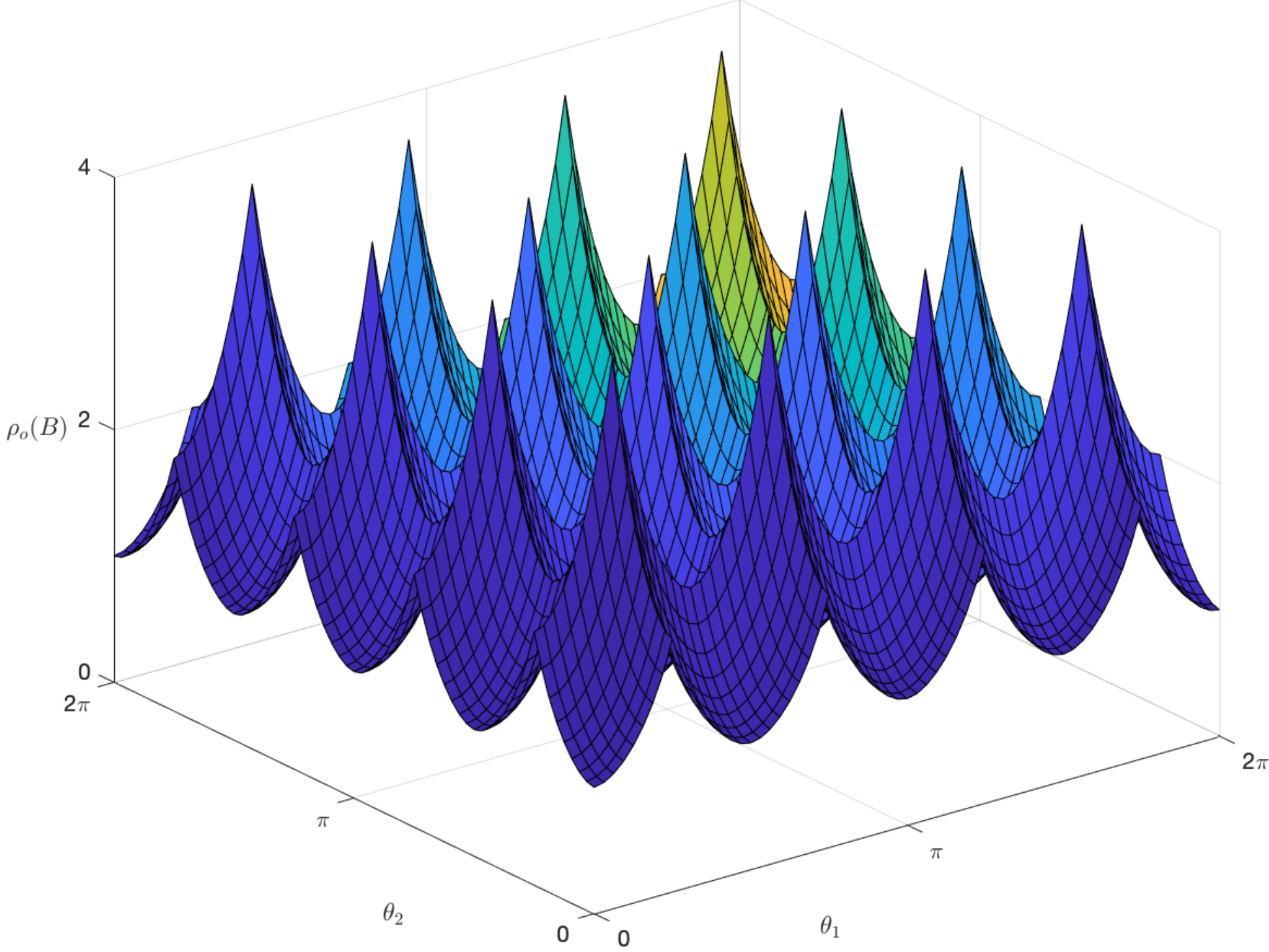}%
    \label{fig:4x4_gf_alt_pp}%
    }%
\\
\subfloat[$\rho_\infty^{\operatorname{GENP}}(B)$]{%
    \includegraphics[width=0.48\textwidth]{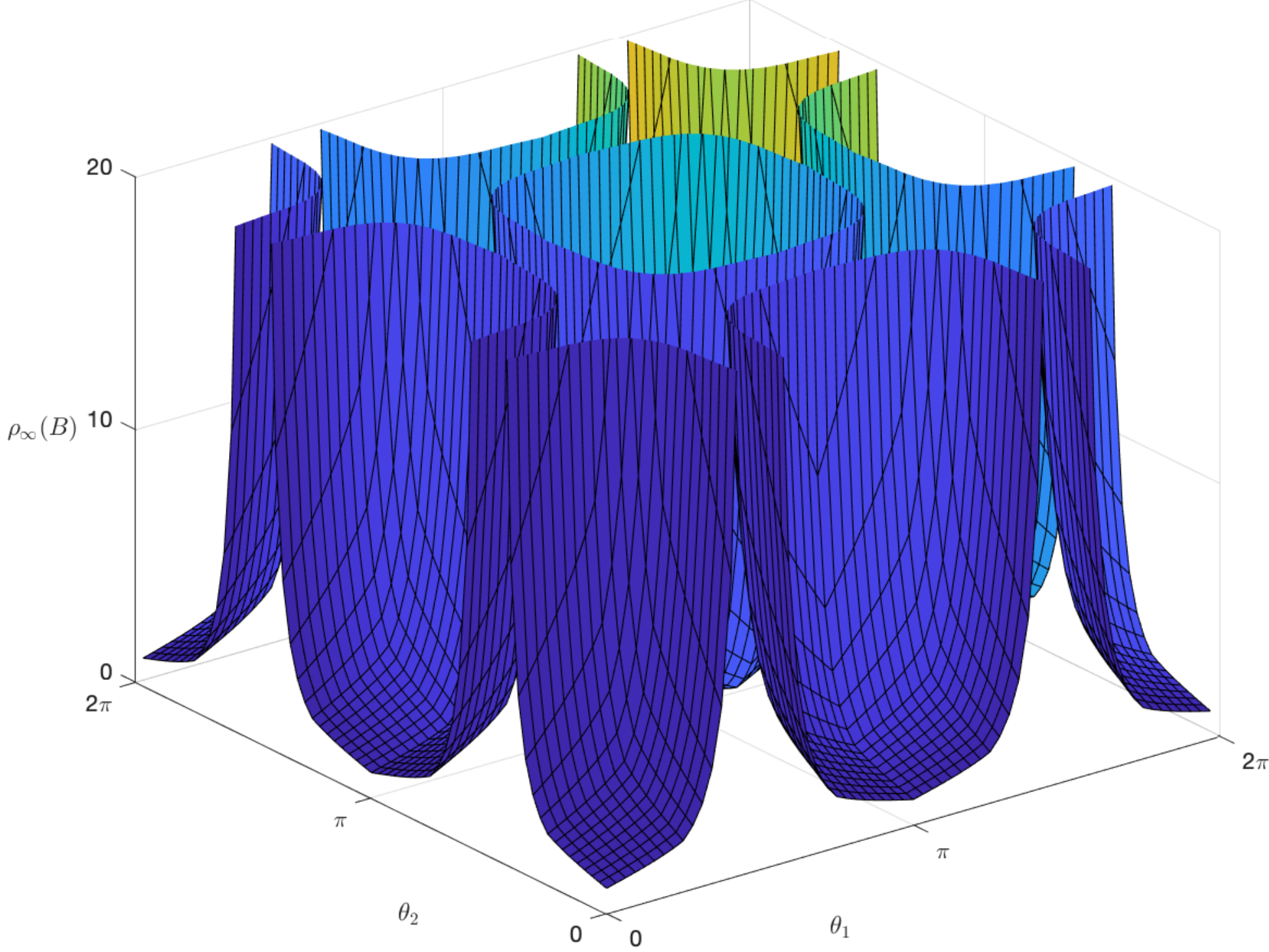}%
    \label{fig:4x4_gf_inf_np}%
    }%
\quad %
\subfloat[$\rho_\infty^{\operatorname{GEPP}}(B)$]{%
    \includegraphics[width=0.48\textwidth]{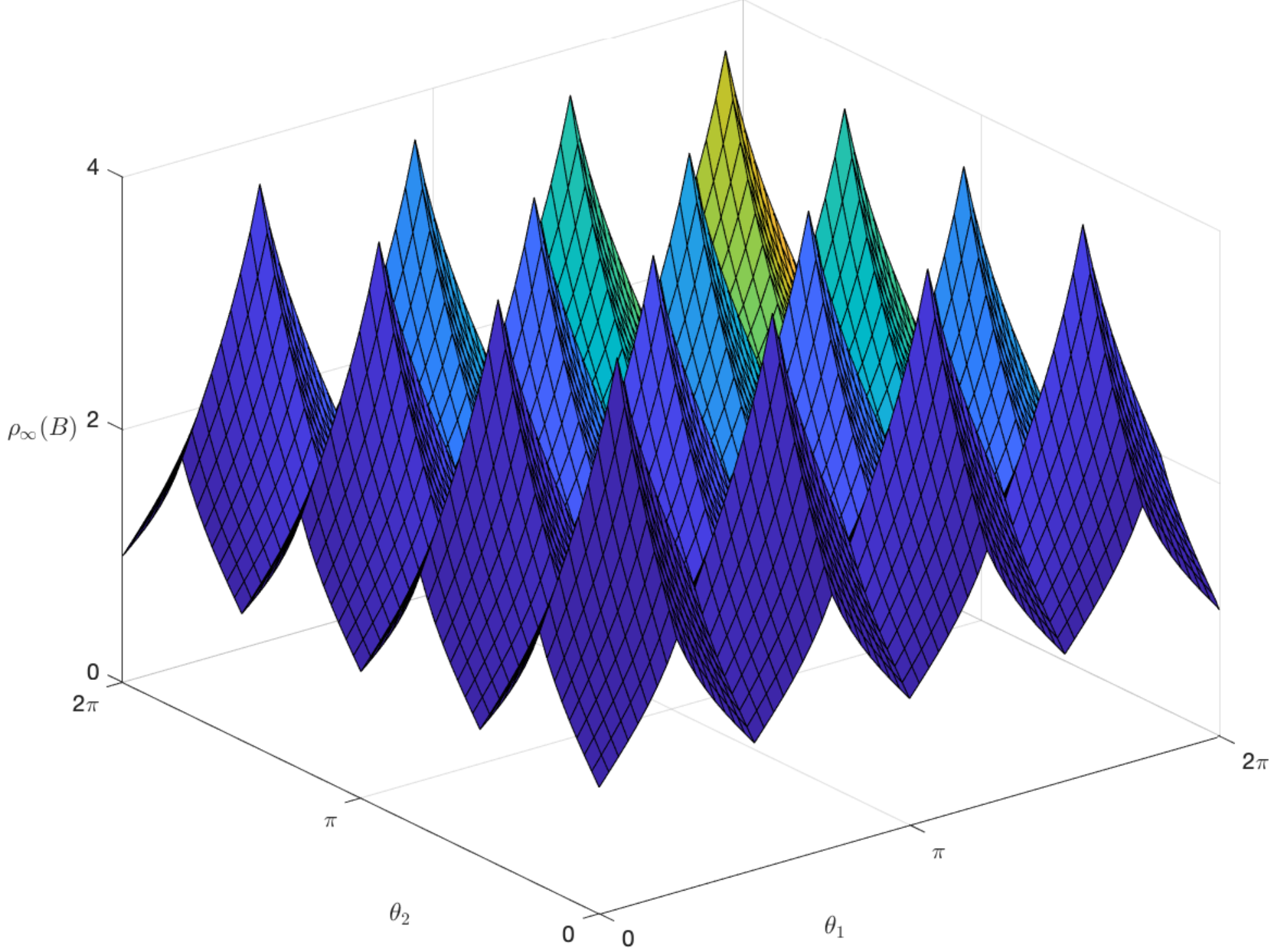}%
    \label{fig:4x4_gf_inf_pp}%
    }%
\caption{$\rho$, $\rho_o$, and $\rho_\infty$ using $B = B(\t_1,\t_2) \in \B_s(4)$}
\label{fig:4x4_gf}
\end{figure}



\Cref{fig:4x4_gf_maxnorm_np,fig:4x4_gf_alt_np,fig:4x4_gf_inf_np} show maps of $\rho(B)$, $\rho_o(B)$, and $\rho_\infty(B)$ using GENP for $B = B(\theta_1,\theta_2) \in \B_s(4)$, mapped against the input angles used to generate $B$. Note the singularity lines in \Cref{fig:4x4_gf_maxnorm_np,fig:4x4_gf_alt_np,fig:4x4_gf_inf_np} account for when $\theta_i = \pi \pm \frac\pi2$ for some $i$. These coincide with when $B_{11} = \cos\t_1\cos\t_2 = 0$ and GENP factorization fails on the first step (and in fact this shows this is the only step on which GENP can fail).
\Cref{fig:4x4_gf_maxnorm_pp,fig:4x4_gf_alt_pp,fig:4x4_gf_inf_pp} show maps of $\rho(B)$, $\rho_o(B)$, and $\rho_\infty(B)$ using GEPP (or GERP) for $B = B(\theta_1,\theta_2) \in \B_s(4)$.  Note each of the peaks in \Cref{fig:4x4_gf_maxnorm_pp,fig:4x4_gf_alt_pp,fig:4x4_gf_inf_pp} occur precisely at the scaled \emph{Hadamard matrices}, $B(\boldsymbol\t) \in \B_s(N)$ for $\boldsymbol \t \in (\frac\pi4+\frac\pi2\mathbb Z)^n$ so that $\sqrt{N} B(\boldsymbol \t) \in \{\pm 1\}^{N\times N}$ with orthogonal rows and columns. As such, butterfly models can be used as a continuous approximation of Hadamard matrices to derive other desirable properties. Hadamard matrices then are extreme points for these butterfly models with respect to growth factors. {Walsh}(-Hadamard) transformations are a popular choice to efficiently randomize the elements of a matrix \cite{St99}. With respect to the growth factor, then Walsh transforms appear to be the least desirable preconditioners in the butterfly ensembles. This will be  explored further in future work.

The relationship $1 \le \rho(B) \le \rho_\infty(B)$ can easily be viewed for this case, as \Cref{fig:4x4_gf_maxnorm_np} fits inside \Cref{fig:4x4_gf_inf_np} and \Cref{fig:4x4_gf_maxnorm_pp} fits inside \Cref{fig:4x4_gf_inf_pp}. Similarly, the relationships to $\rho_o(B)$ using GEPP can be viewed by iteratively stacking \Cref{fig:4x4_gf_maxnorm_pp,fig:4x4_gf_alt_pp,fig:4x4_gf_inf_pp}

\subsubsection{Numerical experiments}
\label{subsec: naive experiments}

\begin{figure}[htbp]
  \centering
\includegraphics[width=\textwidth]{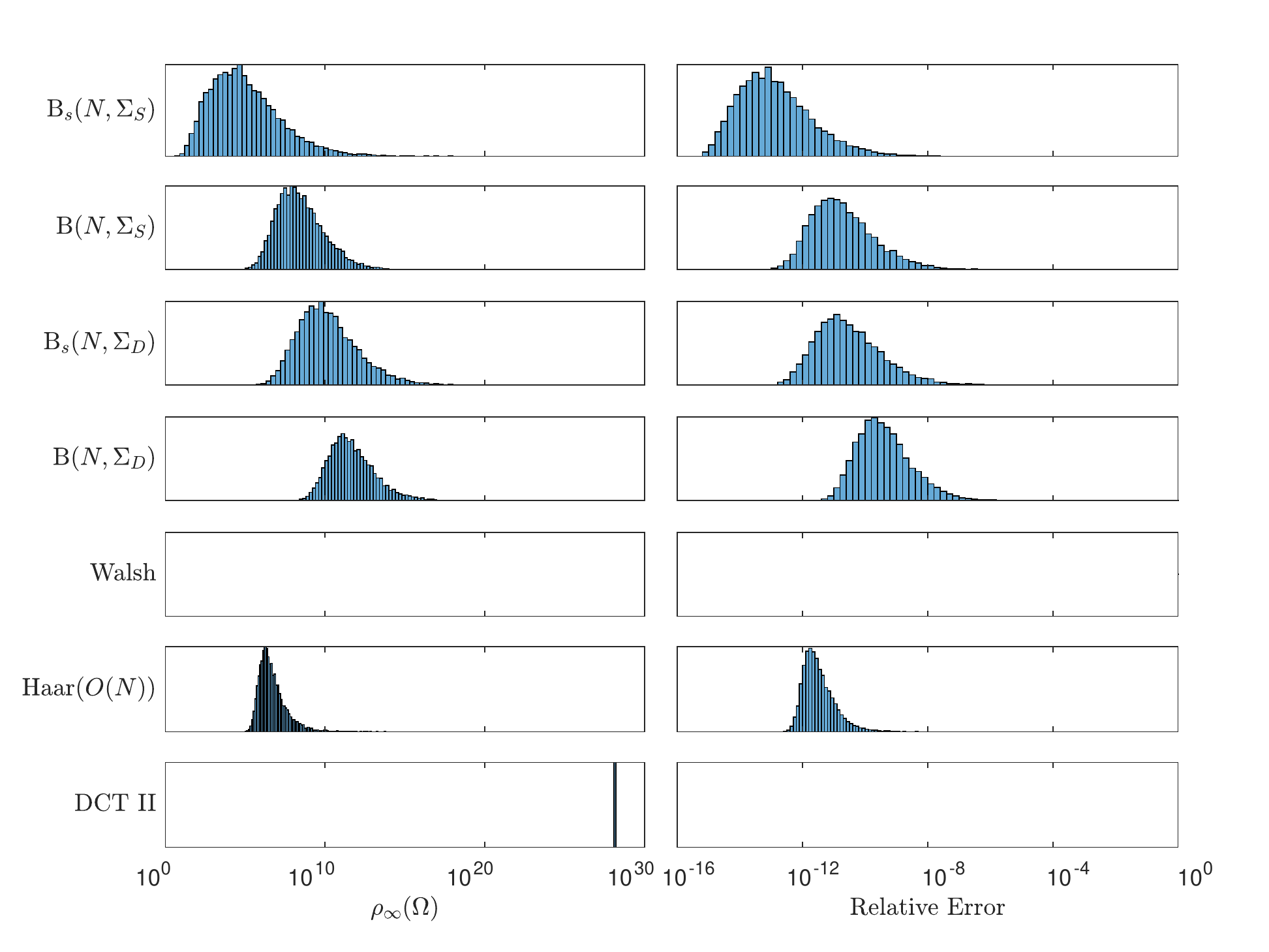}
  \caption{Na\"ive model: $\rho_\infty$ and relative errors using GENP, $N = 2^8$, 10,000 trials}
  \label{fig:hist_np}
\end{figure}

\begin{figure}[htbp]
  \centering
  \includegraphics[width=\textwidth]{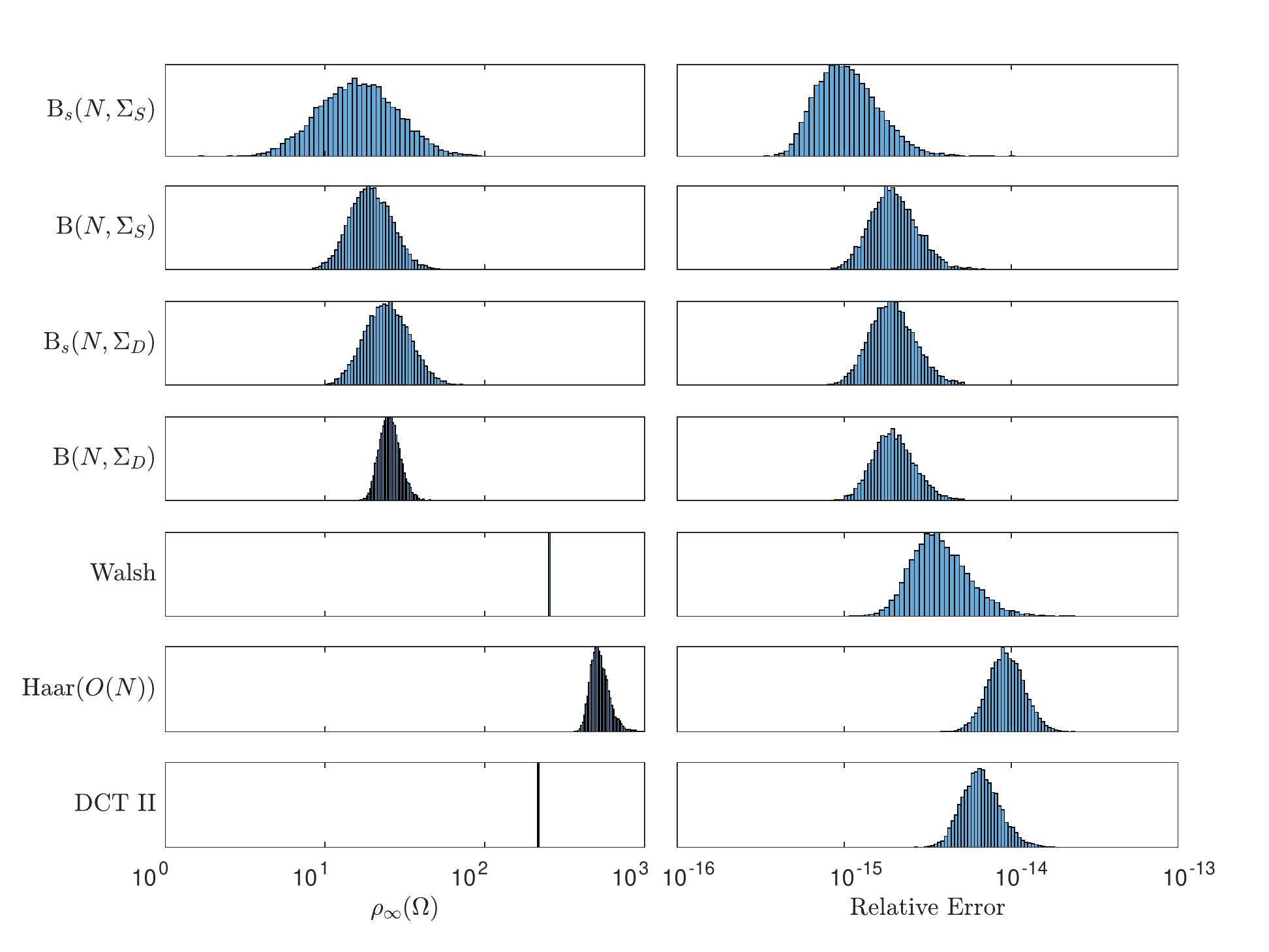}
  \caption{Na\"ive model: $\rho_\infty$ and relative errors using GEPP, $N = 2^8$, 10,000 trials}
  \label{fig:hist_pp}
\end{figure}

\begin{figure}[htbp]
  \centering
\includegraphics[width=\textwidth]{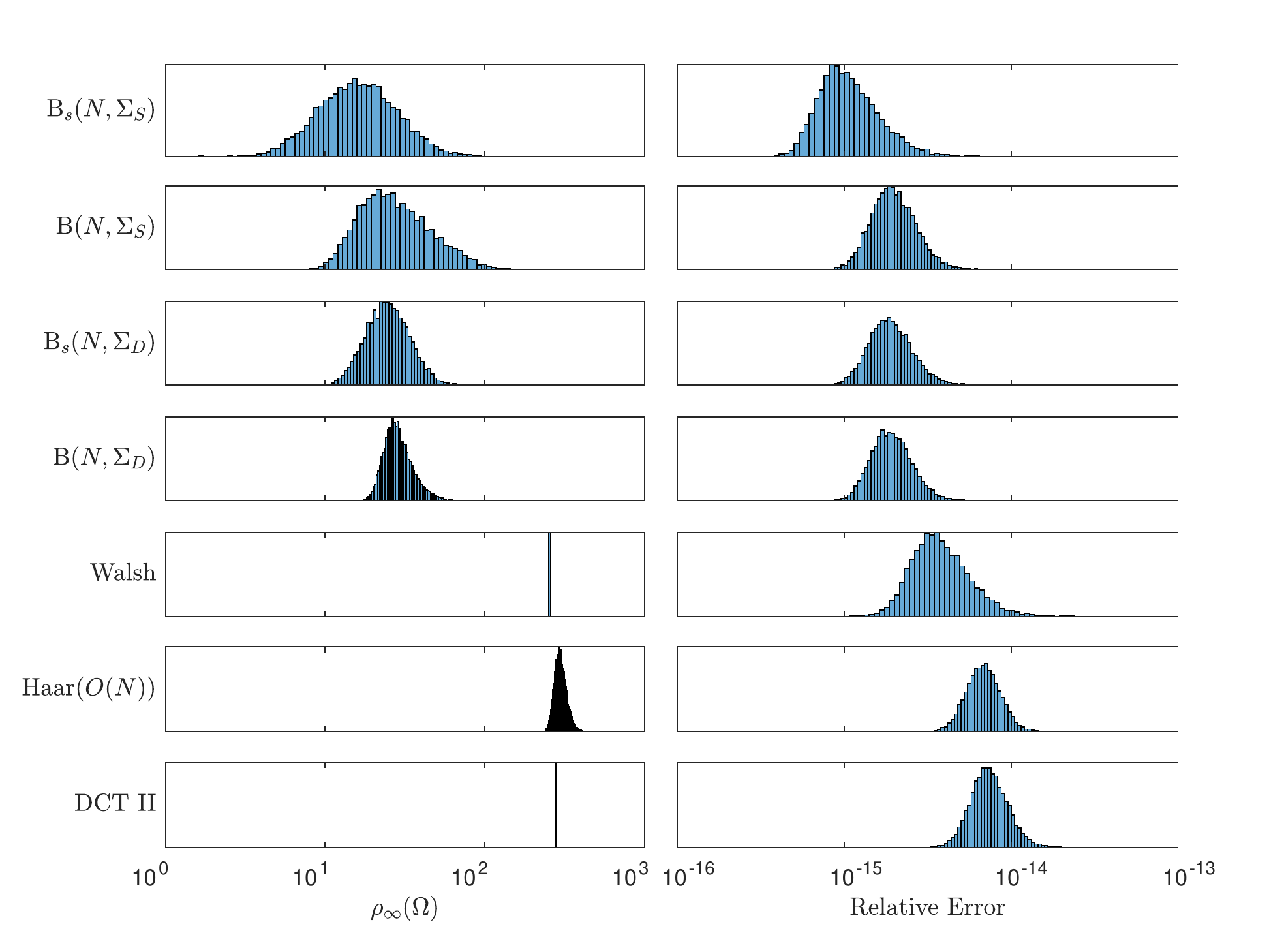}
  \caption{{Na\"ive model: $\rho_\infty$ and relative errors using GERP, $N = 2^8$, 10,000 trials}}
  \label{fig:hist_rp}
\end{figure}

\begin{figure}[htbp]
  \centering
  \includegraphics[width=\textwidth]{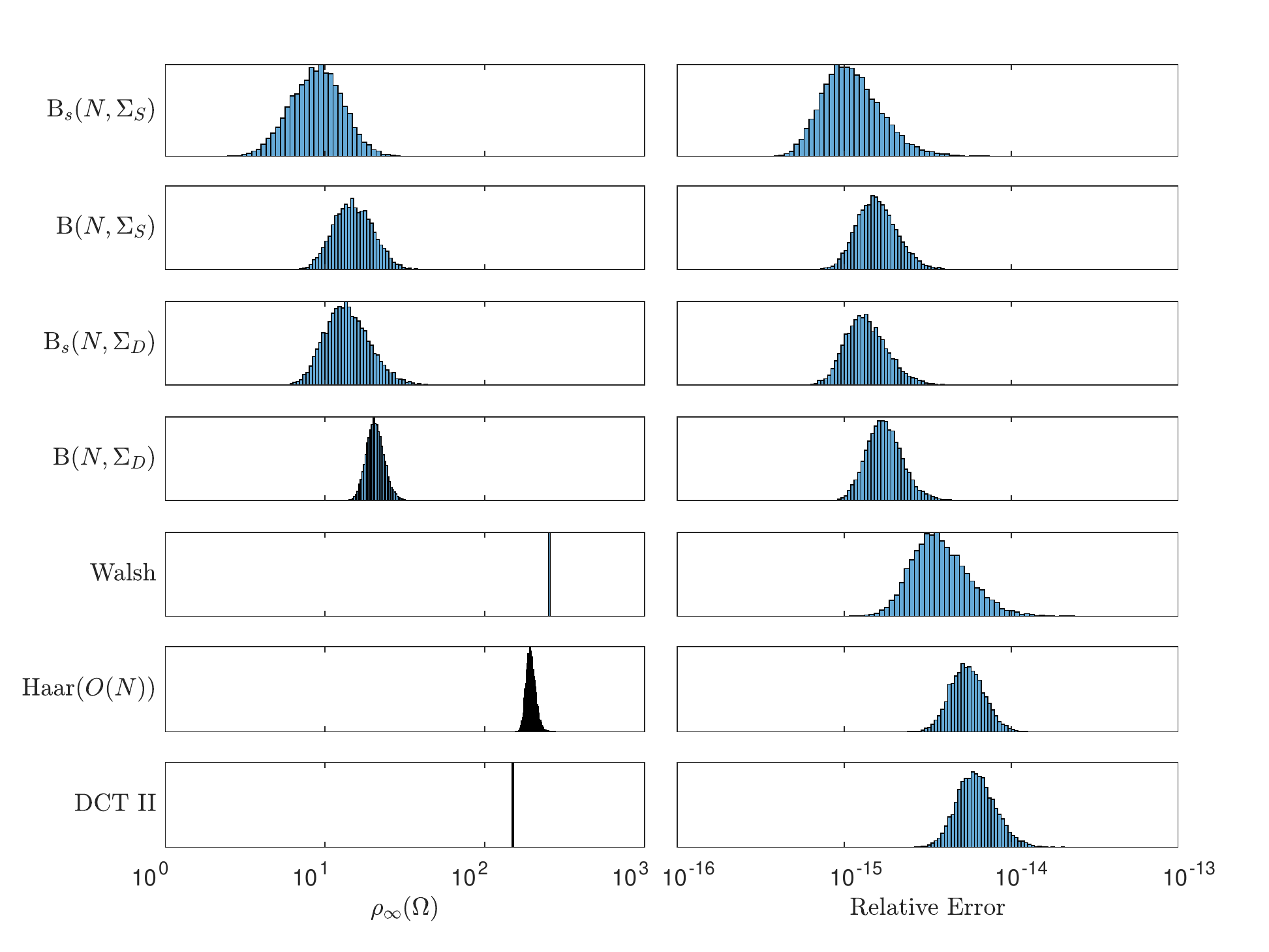}
  \caption{{Na\"ive model: $\rho_\infty$ and relative errors using GECP, $N = 2^8$, 10,000 trials}}
  \label{fig:hist_cp}
\end{figure}

In this section, we will run numerical experiments to compare the impact on the growth factor and relative errors in the na\"ive model when using a variety of random preconditioners $\Omega$ with GENP, GEPP, {GERP,} and GECP. These will only include the 1-sided transformations, so that the growth factors of the preconditioning matrices can be explored directly. Experiments are run using MATLAB in double precision, where $\epsilon = 2^{-52}$ ($\epsilon \approx 2.220446 \cdot  10^{-16})$. We will run $10,000$ trials for each preconditioner for the na\"ive model $\V I \V x = \V b$, where $\V x \sim \Uniform(\mathbb S^{N-1})$ and $N = 2^n$ for $n=2$ to $n=8$. To ease the following discussion, we choose $N=2^8$ as we feel it is representative of the performance we observed for other choices of $N$. For computational simplicity, we will only consider $\rho_\infty$. Additionally, one step of iterative refinement will be run.\footnote{Note we actually ran two steps of iterative refinement in all of our experiments, but only marginal benefits were gained.} \Cref{t:naive genp,t:naive gepp,t:naive gerp,t:naive gecp} will show the sample medians, means ($\bar x$) and standard deviations ($s$) for the trials for $n=8$. {Moreover, \Cref{fig:hist_np,fig:hist_pp,fig:hist_rp,fig:hist_cp} will summarize the growth factors and relative errors, resp. $\rho_\infty(\Omega)$ and $\frac{\|\V x - \hat{\V x}\|_\infty}{\|\V x\|_\infty}$.} 

The particular random preconditioners we will use in these experiments are
\begin{itemize}
\item $\B_s(N,\Sigma_S)$
\item $\B(N,\Sigma_S)$
\item $\B_s(N,\Sigma_D)$
\item $\B(N,\Sigma_D)$
\item Walsh transform
\item $\Haar(\O(N))$
\item Discrete Cosine Transform (DCT II).
\end{itemize}
{See \Cref{sec:prelim} for more information on the Walsh, Haar orthogonal and DCT transformations.\footnote{For the na\"ive model case, the random signs used for the Walsh and DCT experiments have no impact since the $LU$ factorization differs from the deterministic Walsh and DCT matrix factorizations by having the $U$ factor (and each intermediate $A^{(k)}$ factor) multiplied on the right by a diagonal random sign matrix, while $\|\cdot\|_{\max}$ and $\|\cdot\|_\infty$ are invariant under sign permutations; moreover, $\operatorname{Gin}(N,M)$ is invariant under orthogonal transformations. These do come into play for the worst-case model (see \Cref{subsec: WC experiments}).}}


\begin{remark}
The parameters needed to generate each butterfly matrix (i.e., the number of uniform angles) increases from $n$ for $\B_s(N,\Sigma_S)$, $N-1$ for both $\B(N,\Sigma_S)$ and $\B_s(N,\Sigma_D)$, to $\frac12Nn$ for $\B_s(N,\Sigma_D)$. The ordering of the butterfly models included in the summary tables and figures reflects this increasing  of number of these parameters. These compare to $\binom{N}2=\frac12N(N-1)$ uniform angles that could be used to sample $\Haar(\O(N))$, which can be realized by using Givens rotations to find the $QR$ factorization of $\operatorname{Gin}(N,N)$.
\end{remark}

\begin{table}[ht!]
\centering
{\tiny
\begin{tabular}{r|ccc|ccc|ccc}
     &\multicolumn{3}{c}{Growth factor: $\rho_\infty$} &\multicolumn{3}{|c|}{Relative error} &\multicolumn{3}{c}{Relative error + Iterative refinement}  \\
     & Median & $\bar x$ & $s$& Median & $\bar x$ & $s$& Median & $\bar x$ & $s$\\ \hline 
$\B_s(N,\Sigma_S)$	&	5.18e+04	&	5.09e+14	&	5.08e+16	&	1.00e-13	&	6.78e-10	&	6.43e-08	&	4.07e-16	&	4.15e-16	&	8.86e-17	\\
$\B(N,\Sigma_S)$	&	1.95e+08	&	1.16e+13	&	5.95e+14	&	1.51e-11	&	6.48e-09	&	3.85e-07	&	4.08e-16	&	4.18e-16	&	1.48e-16	\\
$\B_s(N,\Sigma_D)$	&	1.03e+10	&	1.17e+19	&	1.13e+21	&	2.39e-11	&	1.38e-01	&	1.38e+01	&	4.10e-16	&	2.40e-03	&	2.37e-01	\\
$\B(N,\Sigma_D)$	&	3.00e+11	&	1.43e+16	&	5.32e+17	&	3.16e-10	&	1.65e-01	&	7.58e+00	&	4.06e-16	&	5.29e-10	&	4.17e-08	\\
Walsh	&	   NaN	&	   NaN	&	   NaN	&	   NaN	&	   NaN	&	   NaN	&	   NaN	&	   NaN	&	   NaN	\\
$\Haar(\O(N))$	&	3.27e+06	&	2.68e+09	&	1.05e+11	&	2.61e-12	&	1.12e-11	&	8.54e-11	&	1.04e-15	&	1.07e-15	&	2.12e-16	\\
DCT II	&	1.33e+28	&	1.33e+28	&	0.00e+00	&	4.33e+12	&	5.48e+12	&	4.61e+12	&	3.89e+14	&	4.94e+14	&	4.16e+14	
\end{tabular}
}
\caption{Na\"ive model numerical experiments for GENP with 10,000 trials for $N=2^8$}
\label{t:naive genp}
\end{table}

\begin{table}[ht!]
\centering
{\tiny
\begin{tabular}{r|ccc|ccc|ccc}
     &\multicolumn{3}{c}{Growth factor: $\rho_\infty$} &\multicolumn{3}{|c|}{Relative error} &\multicolumn{3}{c}{Relative error + Iterative refinement}  \\
     & Median & $\bar x$ & $s$& Median & $\bar x$ & $s$& Median & $\bar x$ & $s$\\ \hline 
$\B_s(N,\Sigma_S)$	&	1.60e+01	&	1.87e+01	&	1.11e+01	&	1.00e-15	&	1.15e-15	&	5.72e-16	&	4.07e-16	&	4.16e-16	&	9.11e-17	\\
$\B(N,\Sigma_S)$	&	1.93e+01	&	2.02e+01	&	6.19e+00	&	1.94e-15	&	2.07e-15	&	6.91e-16	&	4.07e-16	&	4.16e-16	&	8.82e-17	\\
$\B_s(N,\Sigma_D)$	&	2.42e+01	&	2.58e+01	&	8.64e+00	&	1.91e-15	&	2.04e-15	&	6.85e-16	&	4.09e-16	&	4.18e-16	&	8.73e-17	\\
$\B(N,\Sigma_D)$	&	2.54e+01	&	2.59e+01	&	3.87e+00	&	1.96e-15	&	2.08e-15	&	6.58e-16	&	4.07e-16	&	4.15e-16	&	8.94e-17	\\
Walsh	&	2.56e+02	&	2.56e+02	&	0.00e+00	&	3.66e-15	&	4.09e-15	&	1.86e-15	&	3.28e-16	&	3.35e-16	&	6.02e-17	\\
$\Haar(\O(N))$	&	5.20e+02	&	5.32e+02	&	7.52e+01	&	9.27e-15	&	9.59e-15	&	2.42e-15	&	1.04e-15	&	1.07e-15	&	2.10e-16	\\
DCT II	&	2.14e+02	&	2.14e+02	&	1.85e-12	&	6.56e-15	&	6.86e-15	&	1.92e-15	&	3.48e-16	&	3.54e-16	&	6.02e-17	
\end{tabular}
}
\caption{Na\"ive model numerical experiments for GEPP with 10,000 trials for $N=2^8$}
\label{t:naive gepp}
\end{table}

\begin{table}[ht!]
\centering
{\tiny
\begin{tabular}{r|ccc|ccc|ccc}
     &\multicolumn{3}{c}{Growth factor: $\rho_\infty$} &\multicolumn{3}{|c|}{Relative error} &\multicolumn{3}{c}{Relative error + Iterative refinement}  \\
     & Median & $\bar x$ & $s$& Median & $\bar x$ & $s$& Median & $\bar x$ & $s$\\ \hline 
$\B_s(N,\Sigma_S)$	&	1.60e+01	&	1.87e+01	&	1.11e+01	&	1.02e-15	&	1.17e-15	&	5.69e-16	&	4.06e-16	&	4.17e-16	&	9.28e-17	\\
$\B(N,\Sigma_S)$	&	2.59e+01	&	3.12e+01	&	1.79e+01	&	1.95e-15	&	2.08e-15	&	6.78e-16	&	4.06e-16	&	4.16e-16	&	8.96e-17	\\
$\B_s(N,\Sigma_D)$	&	2.45e+01	&	2.59e+01	&	8.36e+00	&	1.87e-15	&	1.99e-15	&	6.31e-16	&	4.08e-16	&	4.18e-16	&	9.08e-17	\\
$\B(N,\Sigma_D)$	&	2.81e+01	&	2.95e+01	&	6.78e+00	&	1.90e-15	&	2.00e-15	&	5.94e-16	&	4.07e-16	&	4.17e-16	&	8.99e-17	\\
Walsh	&	2.56e+02	&	2.56e+02	&	0.00e+00	&	3.64e-15	&	4.11e-15	&	1.89e-15	&	3.28e-16	&	3.34e-16	&	6.00e-17	\\
$\Haar(\O(N))$	&	2.96e+02	&	3.00e+02	&	2.77e+01	&	6.80e-15	&	7.01e-15	&	1.68e-15	&	1.04e-15	&	1.07e-15	&	2.11e-16	\\
DCT II	&	2.79e+02	&	2.79e+02	&	5.63e-12	&	7.32e-15	&	7.65e-15	&	2.04e-15	&	3.48e-16	&	3.54e-16	&	6.13e-17	
\end{tabular}
}
\caption{{Na\"ive model numerical experiments for GERP with 10,000 trials for $N=2^8$}}
\label{t:naive gerp}
\end{table}

\begin{table}[ht!]
\centering
{\tiny
\begin{tabular}{r|ccc|ccc|ccc}
     &\multicolumn{3}{c}{Growth factor: $\rho_\infty$} &\multicolumn{3}{|c|}{Relative error} &\multicolumn{3}{c}{Relative error + Iterative refinement}  \\
     & Median & $\bar x$ & $s$& Median & $\bar x$ & $s$& Median & $\bar x$ & $s$\\ \hline 
$\B_s(N,\Sigma_S)$	&	8.90e+00	&	9.44e+00	&	3.54e+00	&	1.08e-15	&	1.22e-15	&	5.77e-16	&	4.08e-16	&	4.18e-16	&	9.10e-17	\\
$\B(N,\Sigma_S)$	&	1.50e+01	&	1.57e+01	&	4.63e+00	&	1.56e-15	&	1.64e-15	&	4.75e-16	&	4.06e-16	&	4.15e-16	&	8.75e-17	\\
$\B_s(N,\Sigma_D)$	&	1.35e+01	&	1.46e+01	&	5.20e+00	&	1.34e-15	&	1.43e-15	&	4.49e-16	&	4.08e-16	&	4.18e-16	&	8.95e-17	\\
$\B(N,\Sigma_D)$	&	2.06e+01	&	2.09e+01	&	2.70e+00	&	1.73e-15	&	1.81e-15	&	4.76e-16	&	4.07e-16	&	4.15e-16	&	8.72e-17	\\
Walsh	&	2.56e+02	&	2.56e+02	&	0.00e+00	&	3.66e-15	&	4.09e-15	&	1.86e-15	&	3.28e-16	&	3.35e-16	&	6.02e-17	\\
$\Haar(\O(N))$	&	1.93e+02	&	1.95e+02	&	1.35e+01	&	5.43e-15	&	5.61e-15	&	1.35e-15	&	1.04e-15	&	1.07e-15	&	2.12e-16	\\
DCT II	&	7.65e+01	&	7.31e+01	&	1.09e+01	&	6.15e-15	&	6.46e-15	&	1.84e-15	&	3.49e-16	&	3.55e-16	&	6.16e-17	
\end{tabular}
}
\caption{Na\"ive model numerical experiments for GECP with 10,000 trials for $N=2^8$}
\label{t:naive gecp}
\end{table}


\subsubsection{Discussion}
For the na\"ive model, note the Walsh matrix used by MATLAB has leading $2\times 2$ subblock of the form $\begin{bmatrix} 1 & 1 \\ 1 & 1\end{bmatrix}$. Hence, GENP will stop after one step of GE since the second pivot is 0. So for this model, the Walsh transform fails for the GENP experiments, and so does not produce any data used in \Cref{t:naive genp,fig:hist_np}. Although the DCT matrix used by MATLAB has nonvanishing principal minors, small pivots are introduced very early in GENP (e.g., $U_{22} =  -1.3311\cdot 10^{-5}$) so the resulting factorization has a very large growth factor, where the computed $\widehat \rho_\infty=1.3254\cdot 10^{28}$ (as seen in \Cref{t:naive genp}). Using GENP, DCT results in very large relative errors not included in the scaling of \Cref{fig:hist_np}.

In line with \Cref{cor:clt}, $\rho_\infty(B)$ then will look progressively {lognormal} for $B \sim \B(N,\Sigma_S)$. It happens this behavior  seems universal for other butterfly models as well as $\Haar(\O(N))$. This aligns partially with an observation of Trefethen and Schreiber in the iid model cases, which showed near universal behavior that could be modeled by  $\operatorname{Gin}(N,N)$ since only a few steps of GE resulted in near Gaussian entries \cite{TrSc90}. We note also how the supports for the Haar-butterfly and Haar orthogonal $\ell_\infty$-growth factors when using GEPP are essentially disjoint in \Cref{fig:hist_pp}. This aligns with \cref{eq: markov}, so that $n=8$ appears to be large enough to exhibit this behavior.

Overall for GENP, the butterfly models have moderate growth factors and relative errors, with the one step of iterative refinement essentially correcting the propagated errors introduced by (unnecessarily) using GE to solve the na\"ive model. Surprisingly, the simplest butterfly model $\B_s(N,\Sigma_S)$ did the least damage among the preconditioners in the GENP experiments. $\Haar(\O(N))$ outperformed the remaining butterfly models. {For GEPP, GERP, and GECP, the butterfly models outperformed the remaining models (i.e., they did the least damage) while the Walsh transformation outperformed both the DCT  $\Haar(\O(N))$ transformations, with DCT having smaller errors than $\Haar(\O(N))$ only in the GEPP case.} 

In comparing different pivoting strategies, one sees that GENP when combined with one step of iterative refinement outperforms (in terms of relative errors) the GECP relative error experiments  and were on par with the GECP and one step of iterative refinement combination. In line with the purpose of the na\"ive model, these suggest the Haar-butterfly transformations do the least damage when solving the trivial linear system. Additionally, the added computational costs needed to run GECP seemed to go to waste, since GEPP had equivalent performance in conjunction with $\B_s(N,\Sigma_S)$.

\begin{remark}
{Also, we can connect the prior experiments to the theoretical results from \Cref{subsec:gf_haar_b}. Notably, our experiments align with the result in \Cref{thm:gf} that Haar-butterfly matrices have identical GEPP and GERP factorizations, as \Cref{fig:hist_pp,fig:hist_rp} as well as \Cref{t:naive gepp,t:naive gerp} both support. Moreover, using \Cref{thm:gf}, we can compute the mean and the standard deviation for $\rho_\infty(B)$ exactly for GEPP and GERP for $B \sim \B_s(N,\Sigma_S)$: in both cases, we have $\E \rho_\infty(B) = N^\gamma = \left(1+\frac{\log 4}\pi\right)^8\approx 18.619399429$ for $N = 2^8$  (cf. \Cref{cor: exp gf}) and a similar computation yields 
\begin{equation}
    \operatorname{Var} \rho_\infty(B) = \left(4\cdot \frac{1+ \log 2}\pi\right)^8 - \left(1+\frac{\log 4}\pi\right)^{16} \approx 119.803995751
\end{equation}
so that the population standard deviation is $\sigma_{\rho_\infty(B)} = \sqrt{\operatorname{Var} \rho_\infty(B)} \approx 10.9455011649$. These essentially match the respective sample mean and standard deviation computations of 18.7152 and 11.1003 from \Cref{t:naive gepp} and 18.6695 and 11.058 from \Cref{t:naive gerp}.
}

{
Next, we can connect the theoretical and numerical median estimates for Haar-butterfly matrices. Since we are scaling the axis to show logarithmic growth factors, the middle of the then approximate Gaussian curves aligns with the median rather than the average of $\rho_\infty$ for both the GENP (when the average is actually infinite) and GEPP (and GERP) experiments. Additionally, in line with \Cref{cor:median} where we showed the median of the max-norm growth factor, $\rho(B)$, for $B \sim \B_s(N,\Sigma_S)$ is approximately $N^2$, we have $\mu_\infty = \E \ln(\rho_\infty(B)) = \frac{2G}\pi + \frac54 \ln2$ so that $\frac{\mu_\infty}{\ln 2} \approx 2.0912669407$ provides an approximate median, $M_{n,\infty}$, for the $\ell_\infty$-growth factor, $\rho_\infty(B)$, where $M_{n,\infty} = N^{\mu_\infty /\ln 2+o(n^{-1/2})}$. \Cref{fig:hist_np} shows the marker representing $\E \log_{10}(\rho_\infty(B)) = \frac{n\mu_\infty}{\ln 10} \approx 5.036$ (so about halfway between $10^0$ and $10^{10}$ in \Cref{fig:hist_np}) is not yet a good approximation for the GENP median $M_{n,\infty}$ when $n=8$ since the plot of $\rho_\infty(B)$ is not sufficiently lognormal: the estimate $N^{\mu_\infty/\ln 2} \approx 108,710$ is about twice as large as the sample median $\widehat M_n = 51,776$ from \Cref{t:naive genp}. However, the GEPP and GERP plots for $\rho_\infty$ in \Cref{fig:hist_pp,fig:hist_rp} appears much closer to lognormal than the GENP case, so $\mu_\infty = \frac{\ln 2}2$ yields the marker at $\E \log_{10}(\rho_\infty(B)) = \frac{n\mu_\infty}{\ln 10} = \frac{4 \ln 2}{\ln 10} \approx 1.20411998$ (between $10^1$ and $10^2$ in \Cref{fig:hist_pp,fig:hist_rp}) where $n=8$ provides a much better estimator for the median of \Cref{fig:hist_pp,fig:hist_rp} using GEPP or GERP: here $10^{\frac{4 \ln 2}{\ln 10}}=2^4 = 16$ provides a very good approximation for the computed sample medians 16.0059 from \Cref{t:naive gepp} and 15.987 from \Cref{t:naive gerp}. 
}
\end{remark}

\subsection{Worst-case model}
\label{sec:wc}

Wilkinson established a weak form of backward stability of GEPP by showing $\rho(A) \le 2^{m-1}$ for any nonsingular $A \in \mathbb R^{m \times m}$ \cite{Wi61,Wi65}. In \cite[pg. 202]{Wi61}, Wilkinson further shows this bound on worst-case growth factor is sharp, using the following construction:
\begin{equation}
    A_m = \V I_m - \sum_{i>j} \V E_{ij} + \sum_{i=1}^{m-1} \V E_{im}.
\end{equation}
By design, GEPP would carry out without using any pivoting on $A_m$ at any intermediate GE steps, so that the $A_m=LU$ factorizations coincide for GENP and GEPP, where 
    $L = \V I_m - \sum_{i>j} \V E_{ij}$ {and} $U = \V I_m -\V E_{mm}+ \sum_{i=1}^m 2^{i-1} \V E_{im}.
$ 
It follows 
    $\rho(A_m) = |U_{mm}| = 2^{m-1}.    
$ 
For example, 
\begin{equation}
    A_4 = \begin{bmatrix} 1 && &1\\-1&1&&1\\-1&-1&1&1\\-1&-1&-1&1\end{bmatrix} = \begin{bmatrix} 1\\-1&1\\-1&-1&1\\-1&-1&-1&1\end{bmatrix}\begin{bmatrix} 1 & &&1\\&1&&2\\&&1&4\\&&&8\end{bmatrix}
\end{equation}
has $\rho(A_4) = 8=2^3$. 

It happens that also $\rho_\infty(A_m) = 2^{m-1}$ (while $\rho_o(A_m) = 1 + \frac2m(2^{m-1}-1)$), although this is not the upper bound of $\rho_\infty$.\footnote{For example, multiplying the last column of $A_m$ by a constant $c\ge 1$ does not change $\rho$ but $\rho_\infty$ would become $\frac{mc}{m+c-1}2^{m-1}$; letting $c \to \infty$ has $\rho_\infty$ monotonically increase to $m2^{m-1}$.} See \cite{CoPe07} for additional discussions relating to explicit relationships and bounds between different growth factors.

\begin{remark}
Note GECP, {GERP,} or any {\textit{column}} pivoting scheme would result in $PAQ=LU$ for $P = \V I$ and $Q = P_{(2 \ m)(3 \ m) \cdots (m-1 \ m)} = P_{(2 \ m \ m-1 \ \cdots  \ 3)}$ for
    $L = \V I + \sum_{i>j \ge 2} \V E_{ij} - \sum_{i=2}^m \V E_{ij}$ {and}  $U = \V E_{11} + 2\V E_{22} - 2\sum_{i=3}^m \V E_{ii} + \sum_{i=1}^{m-1} \V E_{i,i+1}$, 
where $\rho(A) = 2$ (while $\rho_o(A_m) = 3-\frac2m$ and $\rho_\infty(A_m) = 3-\delta_{m2}$). For example,
\begin{equation}
    A_4 P_{(2 \ 4 \ 3)} = \begin{bmatrix}
    1 & 1\\-1 & 1 & 1\\-1&1&-1&1\\-1&1&-1&-1
    \end{bmatrix} = \begin{bmatrix}
    1\\-1 & 1\\-1 & 1&1\\-1&1&1&1
    \end{bmatrix}
    \begin{bmatrix}
    1 &1\\&2&1\\&&-2&1\\&&&-2
    \end{bmatrix}.
\end{equation}
{In particular, we see GE with \textit{column} partial pivoting leads to very small growth factors for the Wilkinson model. This aligns with the observation by Cort\'ez and Pe\~na that matrices appearing in practical observations for which the (row) GEPP growth factors are very large then have small growth factors when using partial pivoting by columns   \cite{CoPe07}.}
\end{remark}

\begin{remark}
\begin{figure}[t]
  \centering
\includegraphics[width=0.8\textwidth]{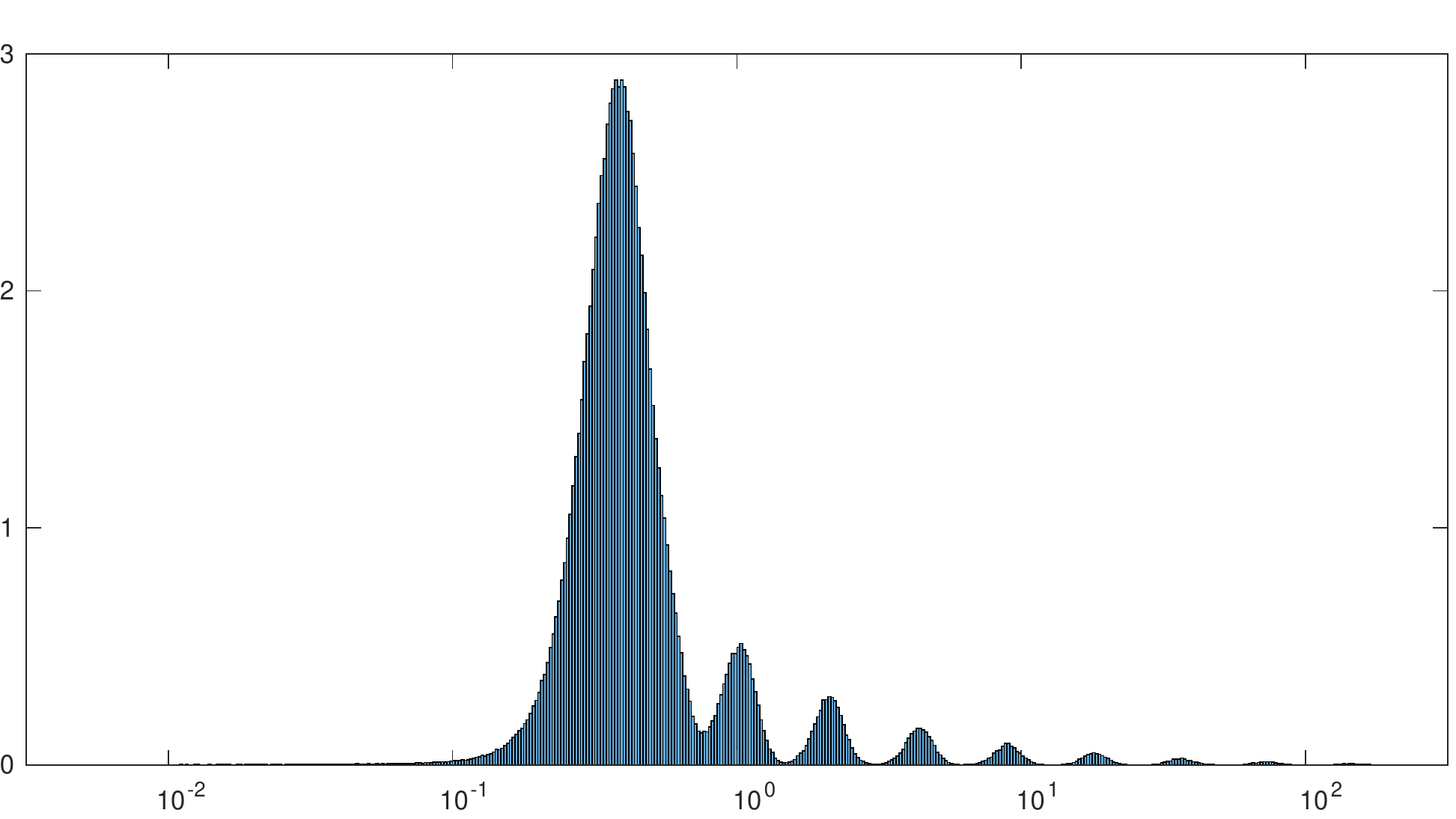}
  \caption{$\ell_2$-relative error using GENP or GEPP to solve $A_N \V x = \V b$ with $\V x \sim \Uniform(\mathbb S^{N-1})$, $10^6$ trials for $N = 2^6$}
  \label{fig:RE_wc_raw}
\end{figure}

Using GENP or GEPP with $A_N$ (which are equivalent for this matrix), accuracy is lost very quickly when using double precision. A direct computation shows $\kappa_{\infty}(A_N) = N$, so the relative error upper bound in \cref{eq:ineq bound max} is 
\begin{equation}
    4N^2\kappa_{\infty}(A_N)\rho(A_N)\epsilon = 4N^32^{N-1} \epsilon = 2^{2+3n+N-1-52} = 2^{2^n+3n-51}.
\end{equation}
Since $2^n+3n-51 < 0$ for $n < 17-\frac{W(131072 \ln(2)/3)}{\ln 2} \approx 5.151584709$ (using the Lambert $W$ function), then the relative errors can exceed 1 for $n\ge 6$.\footnote{Using instead \cref{eq: gf bound inf} similarly yields an upper bound on $n$ of approximately 5.489464984.}  The instability of computed solutions in these scenarios can lead to interesting distributions. For example, \Cref{fig:RE_wc_raw} shows the resulting summary image for $n=6$ of the $\ell_2$-relative error, $\frac{\|\V x - \hat{\V x}\|_2}{\|\V x\|_2}$, running $10^6$ trials using GENP (or GEPP) to solve $A_N\V x = \V b$ for $\V x \sim \Uniform(\mathbb S^{N-1})$.\footnote{Note the $\ell_2$-relative error satisfies \cref{eq:ineq bound max,eq: gf bound inf} with an additional $\sqrt N$ factor. The $\ell_\infty$-relative errors generates a similar  but more rugged image than \Cref{fig:RE_wc_raw}. Moreover, note the $\ell_2$-relative error when using $\V x \in \mathbb S^{N-1}$ is equivalent to the absolute error since $\|\V x\|_2 = 1$.}
\end{remark}

\begin{remark}
{Unlike in the na\"ive model, establishing the exact distributional properties of $UA_NV^*$ for iid $U,V \sim \B_s(N,\Sigma_S)$ is not as tractable since the worst case model loses the Kronecker product forms. One would need to apply the 2-sided Haar-butterfly transform to a Kronecker product matrix to be able to yield full distributional descriptions like in \Cref{thm:gf} (which through the invariance properties of the Haar measure equivalently uses $UAV^*$ for $A = \V I_N$, which can be realized as the $n$-fold Kronecker product of $\V I_2$).}
\end{remark}

\subsubsection{Numerical experiments}
\label{subsec: WC experiments}

For the worst-case model, our goal is to better understand how much of a dampening effect preconditioning by butterfly matrices can have. We will study the two-sided preconditioning by looking at models as outlined in \eqref{eq: precond linear} when studying the linear system $A_N \V x = \V b$ for $\V x \sim \Uniform(\mathbb S^{N-1})$ for $N = 2^n$ and $n = 2$ to 8. We are again running $10,000$ trials computing $\rho_\infty$, relative errors and relative errors after one step of iterative refinement for GENP, GEPP, {GERP,} and GECP. Hence, we are starting with a linear system with a maximal growth factor $\rho$, and we will study how the growth factor  and the corresponding relative error of computed solutions are impacted through using GE to solve the equivalent linear systems
\begin{equation}
    U A_N V^* \V y = U \V b \quad \mbox{and} \quad \V x = V^*\V y
\end{equation}
for independent random transformations $U,V$. For computational simplicity, we will compare $\rho_\infty$, for which we note $\rho_\infty(A_N) = \rho(A_N) = 2^{N-1}$, using the same random transformations used in \Cref{subsec: naive experiments}: $\B_s(N,\Sigma)$ and $\B(N,\Sigma)$ for $\Sigma = \Sigma_S$ and $\Sigma_D$, along with the Walsh transform, Haar orthogonal transform and Discrete Cosine Transform. (See \Cref{sec:prelim,subsec: naive experiments} for more explicit descriptions of the transformations used in these experiments.) \Cref{t:wc genp,t:wc gepp,t:wc gerp,t:wc gecp} will summarize the sample medians, means and standard deviations for the $n=8$ trials.


\begin{table}[ht!]
\centering
{\tiny
\begin{tabular}{r|ccc|ccc|ccc}
     &\multicolumn{3}{c}{Growth factor: $\rho_\infty$} &\multicolumn{3}{|c|}{Relative error} &\multicolumn{3}{c}{Relative error + Iterative refinement}  \\
     & Median & $\bar x$ & $s$& Median & $\bar x$ & $s$& Median & $\bar x$ & $s$\\ \hline 
$\B_s(N,\Sigma_S)$	&	3.30e+05	&	1.93e+13	&	1.93e+15	&	2.20e-12	&	9.10e-10	&	5.72e-08	&	2.60e-15	&	2.87e-15	&	1.24e-15	\\
$\B(N,\Sigma_S)$	&	3.16e+05	&	9.33e+08	&	5.51e+10	&	2.87e-12	&	3.71e-11	&	1.78e-09	&	2.74e-15	&	2.92e-15	&	9.17e-16	\\
$\B_s(N,\Sigma_D)$	&	4.04e+05	&	3.45e+09	&	1.94e+11	&	3.32e-12	&	3.16e-11	&	9.37e-10	&	2.64e-15	&	2.86e-15	&	1.11e-15	\\
$\B(N,\Sigma_D)$	&	2.27e+05	&	2.10e+10	&	1.93e+12	&	2.33e-12	&	6.63e-06	&	6.56e-04	&	2.78e-15	&	7.69e-14	&	7.39e-12	\\
Walsh
&	   1.86e+28	&	   7.38e+42	&	   7.19e+44	&	   8.63e+00	&	   1.60e+12	&	   4.37e+13	&	   8.16e-01	&	   6.33e+16	&	   6.15e+18	\\ 
$\Haar(\O(N))$	&	2.17e+05	&	2.30e+08	&	1.10e+10	&	1.99e-12	&	8.88e-12	&	6.98e-11	&	7.73e-15	&	8.05e-15	&	1.89e-15	\\
DCT II	&	3.91e+05	&	7.16e+08	&	3.34e+10	&	6.38e-12	&	2.06e-09	&	1.90e-07	&	2.57e-15	&	5.40e-15	&	2.74e-13	
\end{tabular}
}
\caption{Worst-case model numerical experiments for GENP with 10,000 trials for $N=2^8$ (excluding 505 Walsh trials where GENP failed)}
\label{t:wc genp}
\end{table}

\begin{table}[ht!]
\centering
{\tiny
\begin{tabular}{r|ccc|ccc|ccc}
     &\multicolumn{3}{c}{Growth factor: $\rho_\infty$} &\multicolumn{3}{|c|}{Relative error} &\multicolumn{3}{c}{Relative error + Iterative refinement}  \\
     & Median & $\bar x$ & $s$& Median & $\bar x$ & $s$& Median & $\bar x$ & $s$\\ \hline 
$\B_s(N,\Sigma_S)$	&	2.67e+01	&	2.97e+01	&	1.51e+01	&	6.85e-15	&	8.94e-15	&	8.81e-15	&	2.59e-15	&	2.89e-15	&	1.27e-15	\\
$\B(N,\Sigma_S)$	&	5.66e+01	&	5.79e+01	&	1.53e+01	&	8.75e-15	&	9.38e-15	&	3.45e-15	&	2.75e-15	&	2.92e-15	&	9.17e-16	\\
$\B_s(N,\Sigma_D)$	&	5.20e+01	&	5.17e+01	&	1.07e+01	&	7.11e-15	&	7.93e-15	&	3.90e-15	&	2.64e-15	&	2.87e-15	&	1.10e-15	\\
$\B(N,\Sigma_D)$	&	7.09e+01	&	6.90e+01	&	1.23e+01	&	8.16e-15	&	8.47e-15	&	2.09e-15	&	2.77e-15	&	2.93e-15	&	8.87e-16	\\
Walsh	&	7.46e+01	&	7.17e+01	&	1.16e+01	&	8.36e-15	&	8.60e-15	&	1.97e-15	&	1.99e-15	&	2.11e-15	&	6.83e-16	\\
$\Haar(\O(N))$	&	7.47e+01	&	7.19e+01	&	1.13e+01	&	1.09e-14	&	1.11e-14	&	2.29e-15	&	7.73e-15	&	8.05e-15	&	1.91e-15	\\
DCT II	&	7.65e+01	&	7.31e+01	&	1.09e+01	&	8.45e-15	&	8.68e-15	&	1.90e-15	&	2.57e-15	&	2.67e-15	&	6.82e-16	
\end{tabular}
}
\caption{Worst-case model numerical experiments for GEPP with 10,000 trials for $N=2^8$}
\label{t:wc gepp}
\end{table}

\begin{table}[ht!]
\centering
{\tiny
\begin{tabular}{r|ccc|ccc|ccc}
     &\multicolumn{3}{c}{Growth factor: $\rho_\infty$} &\multicolumn{3}{|c|}{Relative error} &\multicolumn{3}{c}{Relative error + Iterative refinement}  \\
     & Median & $\bar x$ & $s$& Median & $\bar x$ & $s$& Median & $\bar x$ & $s$\\ \hline 
$\B_s(N,\Sigma_S)$	&	2.59e+01	&	2.66e+01	&	6.96e+00	&	3.05e-15	&	3.22e-15	&	9.63e-16	&	2.59e-15	&	2.88e-15	&	1.26e-15	\\
$\B(N,\Sigma_S)$	&	5.33e+01	&	5.27e+01	&	1.14e+01	&	4.72e-15	&	4.88e-15	&	1.31e-15	&	2.75e-15	&	2.92e-15	&	9.09e-16	\\
$\B_s(N,\Sigma_D)$	&	3.85e+01	&	3.85e+01	&	6.42e+00	&	3.16e-15	&	3.29e-15	&	8.73e-16	&	2.63e-15	&	2.84e-15	&	1.09e-15	\\
$\B(N,\Sigma_D)$	&	6.12e+01	&	5.99e+01	&	6.67e+00	&	5.13e-15	&	5.27e-15	&	1.16e-15	&	2.76e-15	&	2.93e-15	&	8.90e-16	\\
Walsh	&	7.34e+01	&	7.32e+01	&	4.19e+00	&	5.52e-15	&	5.64e-15	&	1.16e-15	&	1.98e-15	&	2.10e-15	&	6.78e-16	\\
$\Haar(\O(N))$	&	7.32e+01	&	7.29e+01	&	4.36e+00	&	9.28e-15	&	9.55e-15	&	2.02e-15	&	7.73e-15	&	8.05e-15	&	1.93e-15	\\
DCT II	&	7.33e+01	&	7.30e+01	&	4.26e+00	&	5.75e-15	&	5.86e-15	&	1.16e-15	&	2.56e-15	&	2.67e-15	&	6.81e-16	
\end{tabular}
}
\caption{{Worst-case model numerical experiments for GERP with 10,000 trials for $N=2^8$}}
\label{t:wc gerp}
\end{table}

\begin{table}[ht!]
\centering
{\tiny
\begin{tabular}{r|ccc|ccc|ccc}
     &\multicolumn{3}{c}{Growth factor: $\rho_\infty$} &\multicolumn{3}{|c|}{Relative error} &\multicolumn{3}{c}{Relative error + Iterative refinement}  \\
     & Median & $\bar x$ & $s$& Median & $\bar x$ & $s$& Median & $\bar x$ & $s$\\ \hline 
$\B_s(N,\Sigma_S)$	&	1.36e+01	&	1.44e+01	&	4.69e+00	&	2.71e-15	&	2.87e-15	&	8.61e-16	&	2.52e-15	&	2.81e-15	&	1.31e-15	\\
$\B(N,\Sigma_S)$	&	4.17e+01	&	4.18e+01	&	9.74e+00	&	3.80e-15	&	3.94e-15	&	1.05e-15	&	2.62e-15	&	2.86e-15	&	1.18e-15	\\
$\B_s(N,\Sigma_D)$	&	1.88e+01	&	1.92e+01	&	4.10e+00	&	2.68e-15	&	2.84e-15	&	8.14e-16	&	2.45e-15	&	2.70e-15	&	1.23e-15	\\
$\B(N,\Sigma_D)$	&	4.66e+01	&	4.58e+01	&	4.92e+00	&	3.93e-15	&	4.06e-15	&	9.14e-16	&	2.63e-15	&	2.87e-15	&	1.20e-15	\\
Walsh	&	6.46e+01	&	6.46e+01	&	2.99e+00	&	4.57e-15	&	4.66e-15	&	9.37e-16	&	2.10e-15	&	2.23e-15	&	8.00e-16	\\
$\Haar(\O(N))$	&	6.43e+01	&	6.43e+01	&	3.10e+00	&	8.76e-15	&	9.03e-15	&	1.93e-15	&	5.87e-15	&	6.30e-15	&	1.92e-15	\\
DCT II	&	6.44e+01	&	6.44e+01	&	3.03e+00	&	4.85e-15	&	4.95e-15	&	9.77e-16	&	2.33e-15	&	2.48e-15	&	8.11e-16	
\end{tabular}
}
\caption{Worst-case model numerical experiments for GECP with 10,000 trials for $N=2^8$}
\label{t:wc gecp}
\end{table}

\begin{figure}[htbp]
  \centering
\includegraphics[width=0.8\textwidth]{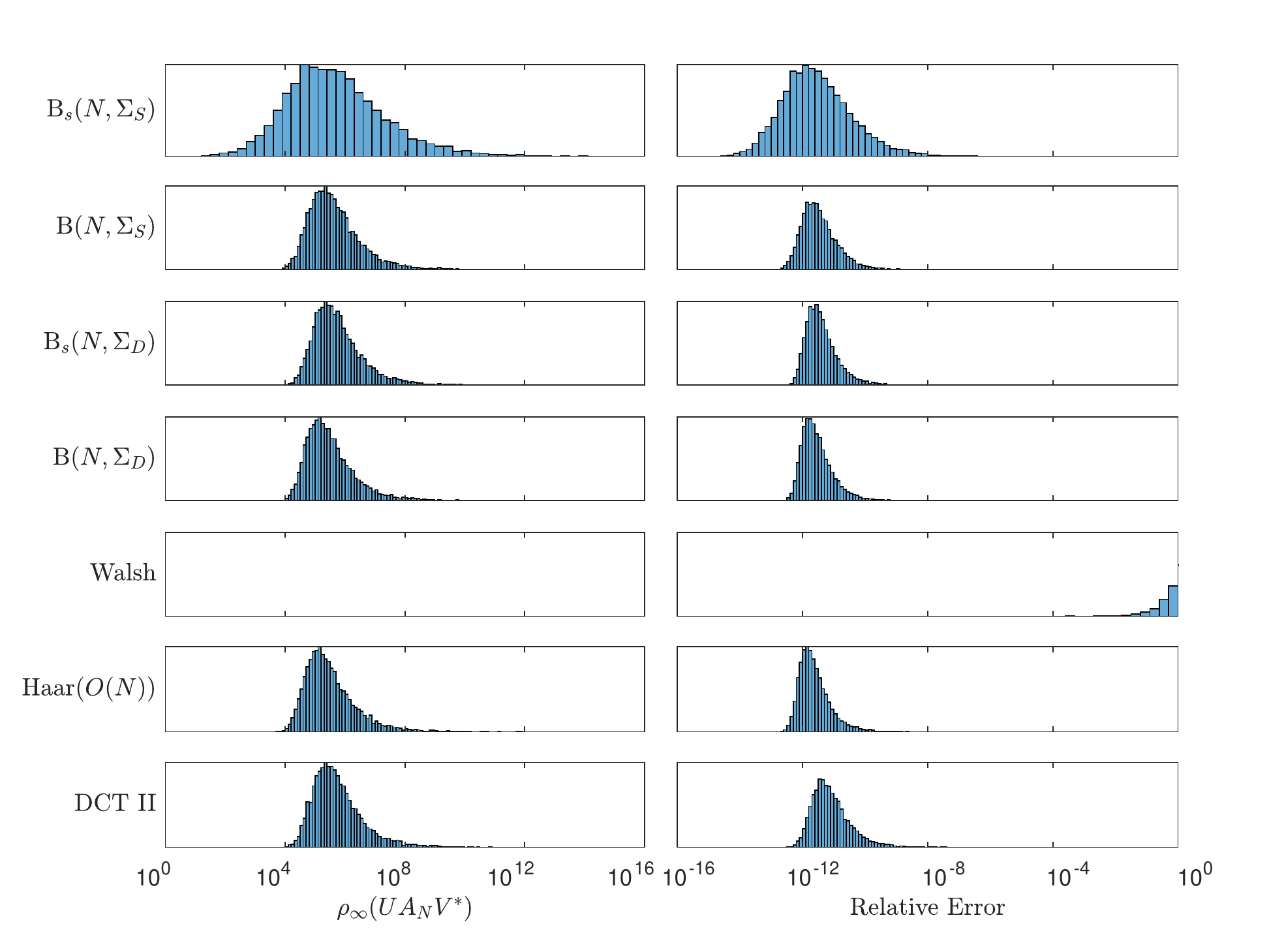}
  \caption{Worst-case model: $\rho_\infty$ and relative errors using GENP, $N = 2^8$, 10,000 trials}
  \label{fig:hist_np_wc}
\end{figure}
\begin{figure}[htbp]
  \centering
  \includegraphics[width=0.8\textwidth]{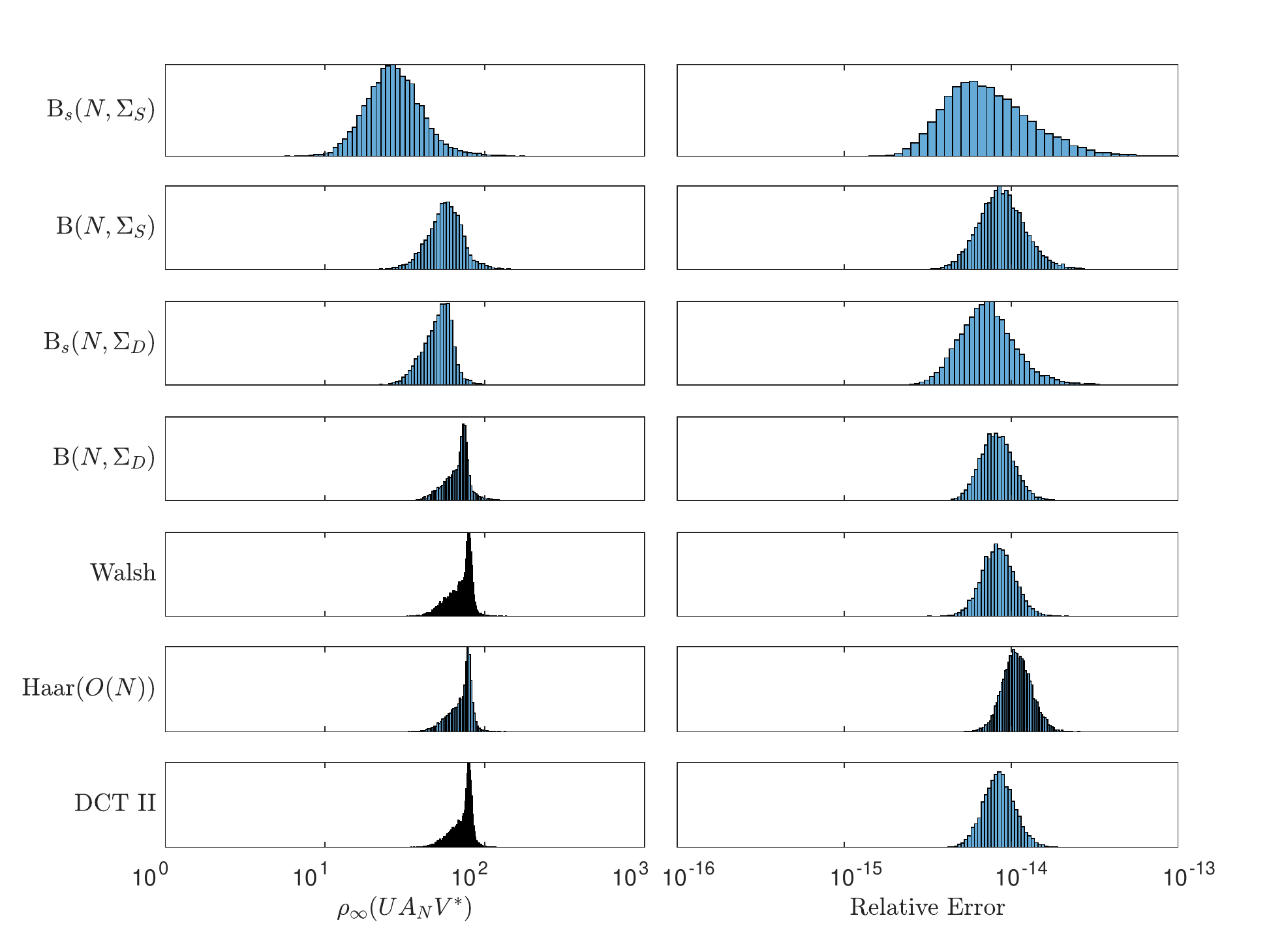}
  \caption{Worst-case model: $\rho_\infty$ and relative errors using GEPP, $N = 2^8$, 10,000 trials}
  \label{fig:hist_pp_wc}
\end{figure}
\begin{figure}[htbp]
  \centering
\includegraphics[width=0.8\textwidth]{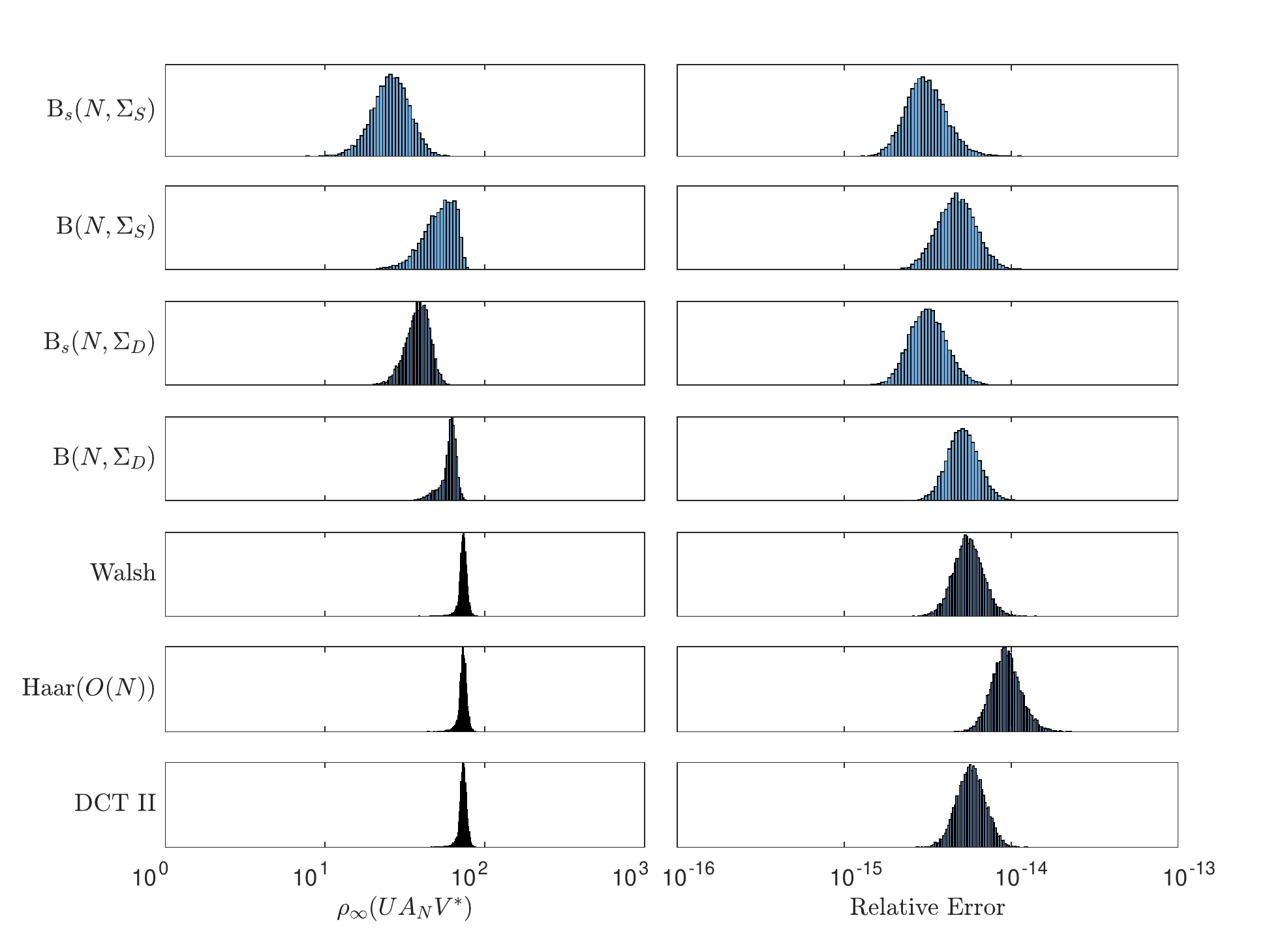}
  \caption{{Worst-case model: $\rho_\infty$ and relative errors using GERP, $N = 2^8$, 10,000 trials}}
  \label{fig:hist_rp_wc}
\end{figure}
\begin{figure}[htbp]
  \centering
  \includegraphics[width=0.8\textwidth]{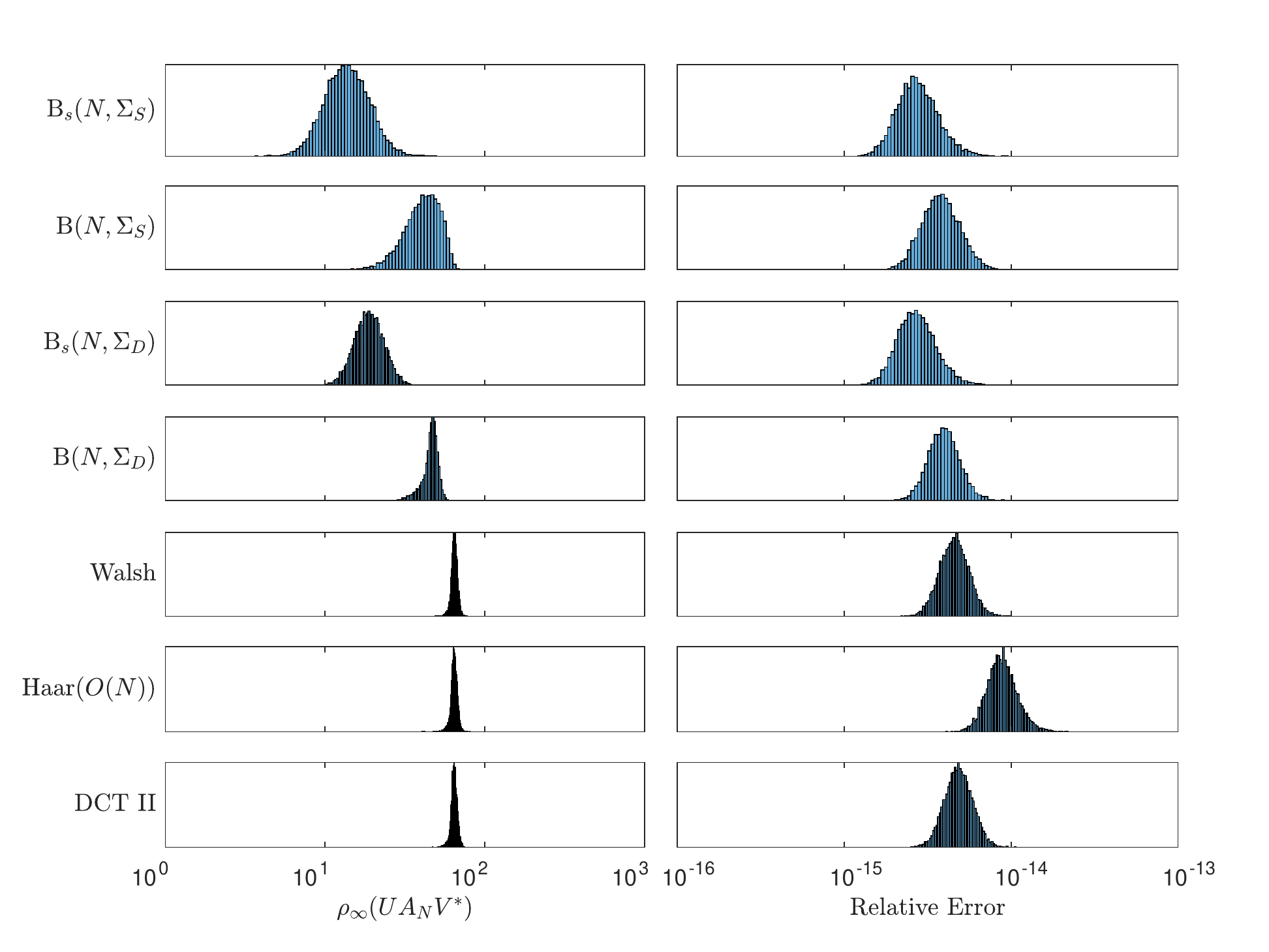}
  \caption{{Worst-case model: $\rho_\infty$ and relative errors using GECP, $N = 2^8$, 10,000 trials}}
  \label{fig:hist_cp_wc}
\end{figure}

\subsubsection{Discussion}

Recall using GEPP (or GENP), we have $\rho_\infty(A_N) = 2^{N-1}$. In particular, for $n = 8$ then $\rho_\infty(A_N) = 2^{255} \approx 5.7896 \cdot 10^{76}$. \Cref{t:wc gepp} shows each random preconditioner has beneficial dampening impacts in terms of reducing the growth factors and the associated relative errors. \Cref{fig:hist_pp_wc} further illustrates this dampening impact with respect to the {growth factors. Each preconditioner has  growth factor approximately less than $10^{2}$, far away from the worst-case scenario.} 

Using GENP, all \textit{except} the Walsh transform have strong dampening impacts on the growth factor, as illustrated in \Cref{t:wc genp} and \Cref{fig:hist_np_wc}. Of interest for the Walsh transform is that the resulting matrix $UA_NV^*$ can now be block degenerate. For our 10,000 trials, 505 using the Walsh transform resulted in block degenerate matrices. (\Cref{t:wc genp} shows the summary of the results for the block nondegenerate trials.) Overall, even when the transformed matrix was block nondegenerate, the results were nowhere near the performance seen in the other random preconditioners considered and led to highly unstable computed solutions. Using GENP, the Haar orthogonal transformations had the best performance on average, while the butterfly models still outperformed the DCT. Among the butterfly models, $\B(N,\Sigma_S)$ and $\B_s(N,\Sigma_D)$ had the smallest growth factors and best precision results. Adding one step of iterative refinement then led to also $\B_s(N,\Sigma_S)$ matching precision of the GEPP, GERP, and GECP experiments. Note also that although $\B_s(N,\Sigma_S)$ had a larger average relative error than most of the other models, its median was  smaller than all but the Haar orthogonal case. Moreover, comparing the histograms in \Cref{fig:hist_np_wc} indicates about a fifth of the time the relative error for GENP is smaller than the left edge of the bulk for each of the remaining models. So even though the Haar-butterfly models led to some samples with the largest relative errors (excluding the Walsh model), they also produced a portion of samples that had higher accuracy than any other model.

For the GERP and GECP experiments, the worst-case model actual is closer in spirit to the na\"ive model since $\rho_\infty(A_N) = 3$ is almost optimal. So these experiments measure how much this preconditioning can mess things up. In this light, all of the butterfly models outperformed the remaining random transformations.

Also of note with respect to the $\B_s(N,\Sigma_S)$ trials is that the logarithmic growth factors in \Cref{fig:hist_np_wc} both look approximately lognormal, while the $\Haar(\O(N))$ logarithmic growth factor generates a non-lognormal histogram, especially in the GEPP case. The remaining butterfly models lie in between and appear to resemble interpolations between these two models.




\subsection{Conclusions}
\label{sec:conclusions}

While Parker had shown one could remove the need of pivoting in GE through the use of butterfly matrices (see \Cref{thm:parker}), we wanted to explore whether this was a good idea in practice. In spite of the benefits one can gain in reducing computation time, one could lose stability when removing pivoting with GE. This idea had been explored previously in \cite{LiLuDo20}, where the authors ran experiments using a combination of butterfly transformations and one step of iterative refinement in a mixed-precision setup to establish how much stability was regained. Our results for both the na\"ive and worst-case models align with theirs, in that one step of iterative refinement seems necessary to stabilize computed solutions when using GENP. We additionally took the direction of comparing performance between different pivoting schemes. The added complexity needed to compute GECP over GEPP did not result in significant gains in precision. In the case of $\B_s(N,\Sigma_S)$, GEPP had nearly identical precision in computed solutions as GECP in the na\"ive model. Additionally, we show that using Haar-butterfly matrices with GENP does the least damage in the trivial linear system. Although one step of iterative refinement more than corrects computed errors using only GENP, the errors are still not too extreme. Future work can also further explore mixed-precision or high precision scenarios.

Our tests also compared butterfly transformations against other transforms taken from randomized numerical linear algebra. For both the na\"ive and worst-case experiments using each pivoting scheme, the butterfly models achieved the highest precision results in the GEPP setup. GENP (with iterative refinement) does remain unstable, which was particularly exemplified in experiments using the Walsh transform. This does not imply no one should use Walsh transforms, since our models are purposely extreme with respect to the max-norm growth factors. In fact, the Walsh transforms do have other desirable properties not shared by the random butterfly models (such as through Hoeffding's inequality to efficiently mix entries in a vector). For the particular setup explored in our experiments, the Walsh transform matches the worst-case growth factor behavior of the Haar-butterfly models as established in \Cref{thm:gf}.  Relationships between the Walsh transform and butterfly matrices, which can be viewed as a continuous generalization of a particular subclass of Hadamard matrices, will be further explored in future work.

Of note is \Cref{thm:gf}, which gives the full distribution of the random growth growth factors for Haar-butterfly matrices using GENP, GEPP and GERP. This is a substantial advancement in the understanding of random growth factors for dense matrices, which had previously focused on empirical first moment estimates and asymptotic bounds for certain cases (see \cite{TrSc90,HiHi20}). Future work will explore whether this result can be extended to GECP.

\begin{figure}[htbp]
  \centering
  \includegraphics[width=0.8\textwidth]{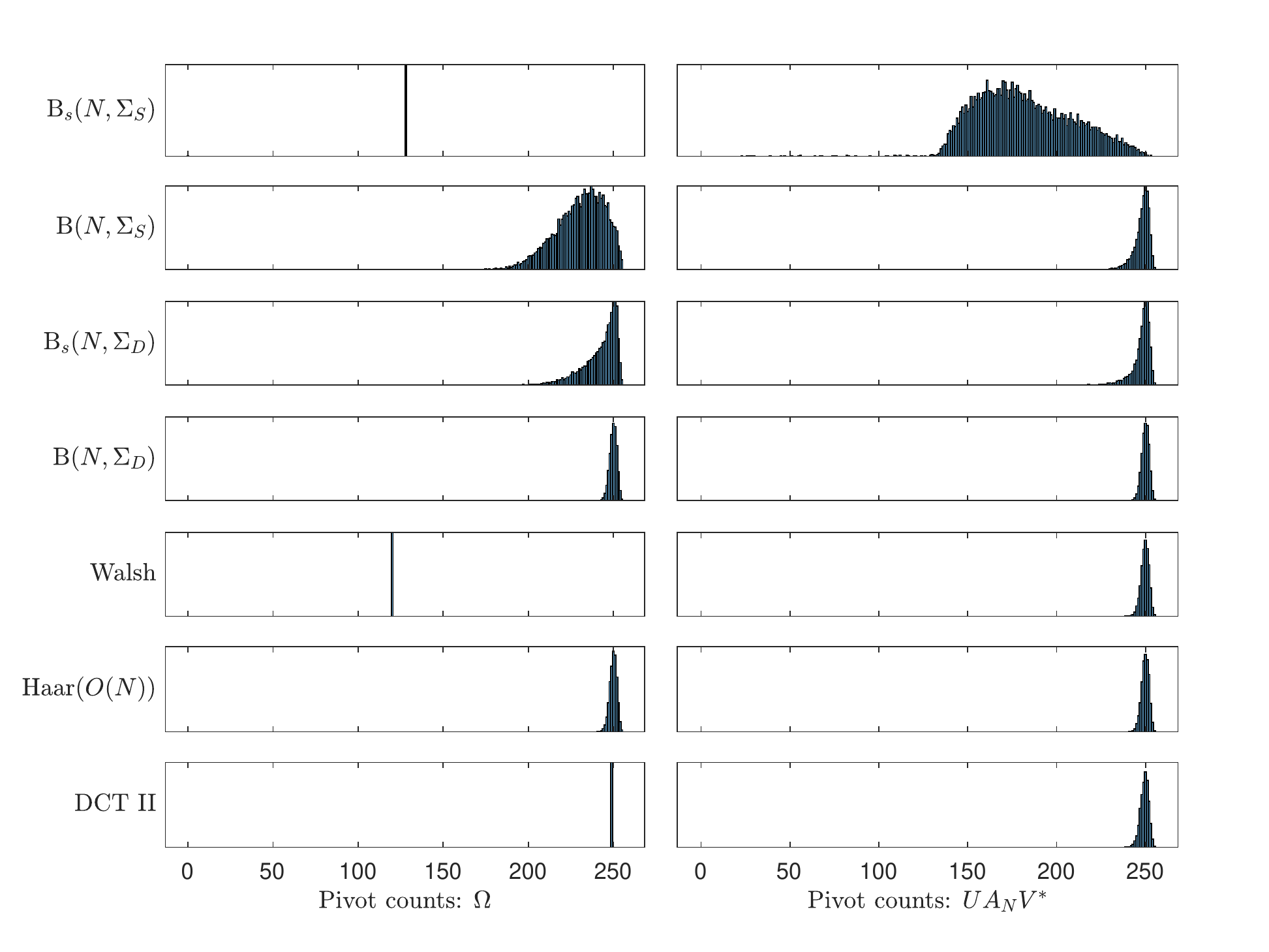}
  \caption{{GEPP pivot movements using the na\"ive and worst case models, using $N = 2^8$ and 10,000 trials}}
  \label{fig:pivots}
\end{figure}

{A related question that will be explored further in a future work regards the number of pivot movements these transformations yield using different pivoting strategies. Parker's original motivation was to provide a means to remove the need for pivoting all together when using GE. Reducing the number of pivot movements would be beneficial for high performance machines to limit the communication overhead for coordinating data movements \cite{baboulin,Pa95}. Our current article explored further the implication of maintaining different pivoting strategies in addition to using these random transformations. So a natural future question is to ask how much actual data movement is there under these transformations. \Cref{fig:pivots,t:pivots} show summaries of the actual pivot movements when using GEPP for our experiments from \Cref{sec:naive,sec:wc}. Both models used matrices that would use no pivot movements with GEPP, so this study is akin to the na\"ive setting, by looking at how much additional movement is being introduced through these transformations. In the worst case model, only the Haar-butterfly transformation had pivot movements that were not near the upper bound for $N-1=255$. In the na\"ive model setting, the Haar-butterfly and Walsh transforms maintained movements no larger than $N/2=128$. The theory and numerical results for GEPP data movements will be further extended in a future study.}

\begin{table}[ht!]
\centering
{
\begin{tabular}{r|ccc|ccc}
     &\multicolumn{3}{c}{Pivot counts: $\Omega$} &\multicolumn{3}{|c}{Pivot counts: $UA_NV^*$}\\
     & Median & $\bar x$ & $s$& Median & $\bar x$ & $s$\\ \hline 
$\B_s(N,\Sigma_S)$	&	128 & 127.56 &   7.45 &  179 &  181.96 &   27.13	\\ 
$\B(N,\Sigma_S)$	& 232 &  230.31 &   14.25 &  249 &  247.97 &   4.29		\\ 
$\B_s(N,\Sigma_D)$	& 245 &  242.01 &   10.38 &   249 &  247.50 &    5.57		\\ 
$\B(N,\Sigma_D)$	&	250 & 249.86 & 2.13 & 250 & 249.89 & 2.11	\\
Walsh  & 120 & 120 & 0 & 250 & 249.66 & 2.27		\\ 
$\Haar(\O(N))$	& 250 & 249.88 & 2.11 & 250 & 249.82 & 2.13		\\ 
DCT II	& 249 & 249 & 0 & 250 & 249.42 & 2.35 
\end{tabular}
}
\caption{{Pivot counts for numerical experiments for GEPP with 10,000 trials for $N=2^8$}}
\label{t:pivots}
\end{table}

Another direction for future work in regards to the impact of preconditioning linear systems by random transformations can explore non-extreme models. Our two models, the na\"ive and worst-case, initiate with a linear systems whose growth factors, respectively, minimize and maximize the GEPP max-norm growth factor $\rho$. These extreme models can further be explored by using ensembles of random matrices with small and large growth factors. For instance, a general structure of matrices with maximal GEPP growth fact is given in \cite{HiHi89}, and can be used to design random ensembles with maximal growth factor $\rho = 2^{N-1}$.

\section*{Acknowledgments}
The authors would like to thank Mike Cranston and Roman Vershynin for many helpful thoughts and insights during this project.

\appendix

\section{Linear algebra background}
\label{sec: background}

Define $A\oplus B \in \mathbb R^{(n_1+n_2)\times (m_1+m_2)}$ to be the block diagonal matrix with blocks $A \in \mathbb R^{n_1\times m_1}$ and $B\in \mathbb R^{n_2\times m_2}$, i.e.,
\begin{equation}
    \label{eq:direct_sum_def}
    A \oplus B = \begin{bmatrix}
    A \\ &B
    \end{bmatrix}.
\end{equation}
%
Define $A \otimes B \in \mathbb R^{n_1n_2 \times m_1m_2}$ to be the \emph{Kronecker product} of $A \in \mathbb R^{n_1 \times m_1}$ and $B \in \mathbb R^{n_2\times m_2}$, given by
\begin{equation}
    \label{eq:kronecker_def}
    A \otimes B = \begin{bmatrix}
    A_{11} B & \cdots & A_{1,m_1} B\\
    \vdots & \vdots & \vdots\\
    A_{n_1,1} B & \cdots & A_{n_1,m_1} B
    \end{bmatrix}.
\end{equation}
Note $(A\otimes B)^* = A^*\otimes B^*$ and $(A\otimes B)^T = A^T \otimes B^T$. Recall $A\otimes B = P(B \otimes A)Q$ for  perfect shuffle permutation matrices $P,Q$ (cf. \cite{perfectShuffle}). If $A$ and $B$ are both square then $Q = P^T$, so that $A \otimes B$ and $B \otimes A$ are conjugate. Also, recall the \emph{mixed-product property} of the Kronecker product: if the products $AC$ and $BD$ can be computed, then the product of the Kronecker products is the Kronecker product of the products, i.e., 
\begin{equation}
    \label{eq:mixed_prod}
    (A\otimes B)(C\otimes D) = AC \otimes BD.
\end{equation}
Certain classes of matrices, $\mathcal S$, are closed under the Kronecker product operator: if $A,B \in \mathcal S$ then $A \otimes B \in \mathcal S$. Straightforward arguments show that this holds for $\mathcal S$ being orthogonal, unitary, permutation, diagonal, lower-triangular or upper-triangular matrices.

The Kronecker product, in particular, allows simplified matrix norm calculations:

 \begin{lemma}
\label{lemma: norm_mult}
If $\|\cdot\|$ is an induced $\ell^p$-matrix norm or $\|\cdot\|_{\max}$, then $\|A \otimes B\| = \|A\|\|B\|$ for $A \in \mathbb R^{n_1\times m_1}$ and $B \in \mathbb R^{n_2\times m_2}$.
\end{lemma}

\begin{proof}

This is a trivial calculation for $\|\cdot\| =\|\cdot\|_{\max}$. If $\|\cdot\|$ is an induced matrix norm, then the multiplicative property is established in \cite[Theorem 8]{LaFa72}.
\end{proof}
These are not the only norms on which Kronecker products satisfy this multiplicative property. Other straightforward computations show $p$-Schatten norms $|\trace(A^p)|^{1/p}$ work (since $\trace(A \otimes B) = \trace(A)\trace(B)$ for square $A,B$), including the Frobenius norm, as well as any standard $p$-norm applied to the vectorization of the matrix $A \in \mathbb R^{n\times m}$ as an element of $\mathbb R^{nm}$ (e.g., $\|A\|_{\max} = \| \operatorname{vec}(A)\|_\infty$). We will restrict our attention in this paper specifically to $\|\cdot\|_{\max}$ and (the induced) $\|\cdot\|_\infty$.

Additionally, the Kronecker product allows for straightforward matrix factorizations determined directly from the corresponding factorizations of each factor. For example:
\begin{lemma}
\label{lemma: kron_factor}
If $P_jA_jQ_j = L_jU_j$ for permutation matrices $P_j,Q_j$ and unit lower triangular $L_j$ and upper triangular $U_j$ and $\displaystyle A = \bigotimes_{j=1}^n A_j$, then $PAQ = LU$ for permutation matrices $\displaystyle P = \bigotimes_{j=1}^n P_j$, $\displaystyle Q = \bigotimes_{j=1}^n Q_j$ and unit lower triangular $\displaystyle L = \bigotimes_{j=1}^n L_j$ and upper triangular $\displaystyle U = \bigotimes_{j=1}^n U_j$.
\end{lemma}

\begin{proof}
This follows directly from the mixed-product property \eqref{eq:mixed_prod}.
\end{proof}

\subsection{Gaussian elimination}
\label{subsec:ge}

Gaussian elimination ({GE}) remains the most prominent approach to solving linear systems $A\V x = \V b$ for $A \in \mathbb R^{n\times n}$. GE without pivoting ({GENP}), when it succeeds, results in the triangular factorization $A = LU$ for unit lower triangular $L$ with
\begin{equation}
    L_{ij} = \frac{A^{(j)}_{ij}}{A^{(j)}_{jj}}
\end{equation}
for $i>j$, where $A^{(k)}$ denotes the matrix before the $k^{th}$ step in GENP with zeros below the first $k-1$ diagonal entries, and $U = A^{(n)}$ is upper triangular with $U_{jj} = A^{(j)}_{jj}$. GENP can be carried out in $\frac23n^3+\mathcal O(n^2)$ FLOPs. 

If $A$ is block nondegenerate, then $A$ has an $LU$ factorization using GENP. If $A$ is block degenerate, then pivoting is necessary. Even when not essential, different pivoting schemes are often employed  to control  errors from using floating-point arithmetic. The most popular modified version of GE is GE with partial pivoting ({GEPP}). This involves an additional scan at GE step $k$ to find the entry of max norm  below the diagonal of $A^{(k)}$ within the column and then a possible row swap to ensure the magnitude of the entry on the $k^{th}$ diagonal is at least as large as those below it. This results in a $PA=LU$ factorization where $P$ is a permutation matrix. By construction, the $L$ factor from GEPP satisfies 
\begin{equation}\label{eq: L cond gepp}
|L_{ij}| \le 1 \quad \text{for any} \quad i>j. 
\end{equation}
GE with complete pivoting ({GECP}) involves a scan through the entire lower untriangularized remaining block of $A^{(k)}$ followed by possible row and column swaps to ensure the magnitude of the entry on the $k^{th}$ diagonal is at least as large as that of all the entries in the lower-left $(n-k+1)\times(n-k+1)$ subblock of $A^{(k)}$. This results in a $PAQ=LU$ factorization for permutation matrices $P,Q$. For this paper, we will assume a pivot search chooses the pivot with minimal taxi cab distance with respect to the row and column indices to the pivot position (i.e., the main leading diagonal entry in the remaining untriangularized subblock), prioritizing minimal row index distance over column index in the case of a tie between candidates.

GE with rook pivoting ({GERP}) involves iteratively scanning first through the column below the diagonal to find a max norm entry, and then within that specific row to find the max norm entry, and repeating between column and row scans, until this process completes with the candidate pivot that maximizes both its row and column. This is followed then by the associated row and column swaps to move the resulting candidate to the pivot position. The name of the pivoting scheme is derived from the limitation on the pivot scans to paths a rook piece could make in a game of chess. See \cite{PoNe00} for further discussion regarding GERP. We will not explore additional numerical experiments for GERP beyond highlighting  connections to our chosen models. Note for this paper we will assume GERP always sequences column scans before row scans at each intermediate step.


The total operational costs of these pivoting schemes differs only in the added scans to find each pivot and associated row and column swaps. $\mathcal O(n^2)$ and $\mathcal O(n^3)$ additional scans are needed, respectively, for GEPP and GECP, while GERP ranges from twice the GEPP complexity (with respect to added scans) to the full GECP complexity. In \cite{PoNe00}, computations with iid models show that on average the number of scans needed for GERP is $\mathcal O(n^2)$ as well, typically accounting for only a factor of 3 more scans than for GEPP.

Using GENP, a standard result establishes the existence and uniqueness of an GENP $LU$ factorization. The particular form presented in \cite[Theorem 9.1]{Hi02} is given below:
\begin{theorem}
\label{thm: unique genp}
  There exists a unique $LU$ factorization of $A \in \mathbb R^{n\times n}$ using GENP if and only if $A_{:k,:k}$ is nonsingular for $k=1,\ldots,n-1$. If $A_{:k,:k}$ is singular for some $k<n$, an $LU$ factorization may exist but it is not unique.
\end{theorem}
With pivoting, the $LU$ factorization is sensitive to row and column permutations: If $PA=LU$ is the $LU$ factorization for $A$ using GEPP, then for $B = QA$ for another permutation matrix $Q$, then we do not necessarily have $(PQ^T)B=LU$ using GEPP. This non-uniqueness of $LU$ factors follows in the event there are any ``ties'' at an intermediate GE pivot search step; in such a case, we have $|A^{(j)}_{ij}| = |A^{(j)}_{jj}|$ so that $|L_{ij}| = 1$ for some $i>j$. 
For example, if $A = \begin{bmatrix}1 & x\\1 & y\end{bmatrix}$ and $B = P_{(1 \ 2)}A = \begin{bmatrix}1 & y\\1 & x\end{bmatrix}$, then using GEPP, we have
\begin{equation}
    A = \begin{bmatrix} 1 & 0 \\ 1 & 1\end{bmatrix}\begin{bmatrix} 1 & x \\ 0&y-x\end{bmatrix} \quad \mbox{and} \quad B = \begin{bmatrix} 1 & 0\\ 1 & 1\end{bmatrix}\begin{bmatrix} 1 & y\\0&x-y \end{bmatrix}.
\end{equation}
%
When ties are avoided, GEPP results in unique $L$ and $U$ factors.
\begin{theorem}
    \label{thm:GEPP_unique}
    Let $A$ be a nonsingular square matrix. Then the $L$ and $U$ factors in the GEPP factorization $PA=LU$ are invariant under row permutations on $A$ iff $|L_{ij}|<1$ for all $i>j$.
\end{theorem}

\begin{proof}
By considering $B = PA$, we can assume $P = \V I$. Suppose first $|L_{ij}| = 1$ for some $i>j$. Let $j$ and then $i$ be minimal such that this occurs. Since  $|B^{(j)}_{jj}| = |B^{(j)}_{ij}|$, then using the row transposition permutation ${(i\ j)}$ would yield a different $P_{\sigma}$ factor with $\sigma(j)  \ne j$, so that $P_{\sigma} \ne \V I$. The $L$ and $U$ factors must differ for $P_\sigma B \ne B$ by \Cref{thm: unique genp}. Now assume $|L_{ij}| < 1$ for all $i > j$. Let $P_\sigma B = L'U'$ be another GEPP factorization of $B$ for some $\sigma \in S_n$. Suppose $\sigma$ is a nontrivial permutation. Let $i$ be the first non-fixed point of $\sigma$, and note then $i < \min(\sigma(i), \sigma^{-1}(i))$. Then
\begin{equation}
    |(L')_{\sigma(i),i}| = \left|\frac{(P_\sigma B)^{(i)}_{\sigma(i),i}}{(P_\sigma B)^{(i)}_{ii}} \right| = \left|\frac{B^{(i)}_{ii}}{B^{(i)}_{\sigma^{-1}(i),i}} \right| = \frac1{|L_{\sigma^{-1}(i),i}|} > 1.
\end{equation}
This contradicts \eqref{eq: L cond gepp}. It follows necessarily $\sigma$ must be the trivial permutation so that $P_\sigma = \V I$. The uniqueness of the $L$ and $U$ factors  follows from \Cref{thm: unique genp}.
\end{proof}
\Cref{thm:GEPP_unique} can be generalized to other pivoting strategies (e.g., GECP), such that the specific $L$ and $U$ factors are unique whenever no ties are encountered in each intermediate pivot search step.

\section{Proofs for \Cref{subsec:gf_haar_b}}
\label{sec:thm proofs}

Note first using \eqref{eq:sbm_def} then we can find the GENP and GEPP factorizations of simple butterfly matrices directly:

\begin{proposition}
    \label{prop: ge factors}
    Let $B = B(\boldsymbol \t) \in \B_s(N)$. If $\cos\t_i \ne 0$ for all $i$, then $B$ has GENP factorization $B = L_{\boldsymbol \t}U_{\boldsymbol \t}$ where $L_{\boldsymbol\t} = \bigotimes_{j=1}^n L_{\t_{n-j+1}}$, $U_{\boldsymbol \t} = \bigotimes_{j=1}^n U_{\t_{n-j+1}}$, and
    \begin{equation}
    \label{eq: genp factors}
        L_{\t} = \begin{bmatrix}
            1 &0 \\-\tan\t & 1
        \end{bmatrix} \quad \mbox{and} \quad U_{ \t} =  \begin{bmatrix}
            \cos\t & \sin\t\\
            0&\sec\t
        \end{bmatrix}
    \end{equation}
    
    Let $\boldsymbol \t ' \in [0,2\pi)^n$ be such that
    \begin{equation}
        \t_i' = \left\{\begin{array}{cl}
        \t_i & \mbox{if $|\tan\t_i|\le 1$,}\\
        \frac\pi2 - \t_i & \mbox{if $|\tan\t_i| > 1$.}
        \end{array}\right.
    \end{equation}
    If $|\tan \t_i| \ne 1$ for any $i$, then the GEPP  factorization of $B$ is $PB = LU$ where  $P = P_{\boldsymbol \t} = \bigotimes_{j=1}^n P_{\t_{n-j+1}}$ for 
    \begin{equation}
    P_\t = \left\{
    \begin{array}{cl}
    \V I_2 & \mbox{if $|\tan \t_i| \le 1$,}\\
    P_{(1 \ 2)} 
    & \mbox{if $|\tan \t_i| > 1$,}  
    \end{array}\right.
    \end{equation}
    $L = L_{\boldsymbol\t'}$, $U = U_{\boldsymbol \t'} D_{\boldsymbol \t}$, $D_{\boldsymbol \t} = \bigotimes_{j=1}^n D_{\t_{n-j+1}}$, and 
    \begin{equation}
        D_\t = \left\{
        \begin{array}{cl}
        \V I_2 & \mbox{if $|\tan \t_i| \le 1$,}\\
        (-1)\oplus 1 
        & \mbox{if $|\tan \t_i| > 1$.}  
        \end{array}\right.
    \end{equation}
    Moreover, for all $k$ we have $(PB)^{(k)} = B(\boldsymbol \t')^{(k)}D_{\boldsymbol \t}$ where $B(\boldsymbol\t') \in \B_s(N)$.
\end{proposition}



\begin{proof}
First consider the GENP case. Note if $\cos\t\ne 0$, then $B(\t) \in \SO(2)$ has an $LU$ factorization with
    $B(\t) 
    = L_\t U_\t.$
The conclusion then follows by \Cref{lemma: kron_factor} and \eqref{eq:sbm_def}. 

For the GEPP case: Let $PB = LU$ be the GEPP factorization of $B$.  Note first using GEPP for $B(\t) \in \SO(2)$, then a row pivot is needed only if $|\tan\t| > 1$. Using \Cref{thm:GEPP_unique}, if we can establish $|L_{ij}| < 1$ for all $i>j$, then the format for $P$  follows immediately from \Cref{lemma: kron_factor}. 
Note 
\begin{equation}
    P_{(1 \ 2)} = \begin{bmatrix}
       0 &1\\1&0
    \end{bmatrix} = \begin{bmatrix}
        0&1\\-1&0
    \end{bmatrix}\begin{bmatrix}
        -1&0\\0&1
    \end{bmatrix}=B\left(\frac\pi2\right)(-1 \oplus 1).
\end{equation}
For 
\begin{equation}
    e_j = \left\{
    \begin{array}{cl}
    1 & \mbox{if $|\tan \t_j| \le 1$,}\\
    0 & \mbox{if $|\tan \t_j| > 1$,} 
    \end{array}\right.
\end{equation}
we have
\begin{equation}
    P_{(1 \ 2)}^{e_j} = B\left(\frac\pi2 e_j\right)((-1)^{e_j}\oplus 1), \quad \mbox{while also} \quad 
    (-1 \oplus 1)B(\t)(-1 \oplus 1) 
    = B(-\t).
\end{equation}
Additionally using the mixed-product property and \eqref{eq:sbm_def}, we have $PB = B'D$ 
for $D = D_{\boldsymbol\t}$  and $B' = B(\boldsymbol \t') \in \B_s(N)$, with $\boldsymbol \t'$ such that 
\begin{equation}
\t_j' = \frac\pi2 {e_{j}}+ (-1)^{e_{j}}\t_{j} = %
\left\{
\begin{array}{cl} 
\t_j & \mbox{if $|\tan\t_j| \le 1$}\\ 
\frac\pi2 - \t_j & \mbox{if $|\tan\t_j|>1$}.
\end{array} 
\right.
\end{equation}
If $B' = L'U'$ is the GENP factorization of $B'$ then $B'D = L'(U'D)$ is the GENP factorization of $B'D$. It follows that
\begin{equation}
    (PB)^{(k)} = (B'D)^{(k)} = (L')^{(k)}B'D = {B'}^{(k)}D.
\end{equation}
The final factorization follows from the GENP case applied to $B'$. The unique form of the factors in $PB=LU$ with $L = L'$ and $U = U'D$ follows from \Cref{thm:GEPP_unique} since $|\tan \t_k'| < 1$ for all $k$ so that $|L_{ij}|<1$ for all $i>j$ by \Cref{lemma: kron_factor}.
\end{proof}

Having this explicit $LU$ factorization of $B \in \B_s(N)$ allows us to also construct each of the intermediate matrices $B^{(k)}$.

\begin{lemma}
\label{lemma:Bk}
    Suppose $A \in \mathbb R^{N/2\times N/2}$ has an $LU$ factorization using GENP. Let
    \begin{equation}
    B = \begin{bmatrix}
    \cos\theta A & \sin\theta A\\-\sin\theta A & \cos\theta A
    \end{bmatrix} = B(\t) \otimes A
    \end{equation}
    for $\cos\theta \ne 0$. Then $B$ has an $LU$ factorization using GENP. Moreover, if $k \le N/2$, then
    \begin{equation}
    \label{eq:Bk first half}
        B^{(k)} = \begin{bmatrix}
        \cos\theta A^{(k)} & \sin\theta A^{(k)}\\
        -\sin\theta \begin{bmatrix}
        \V0&\V0\\\V0&\V I_{N/2-k+1}
        \end{bmatrix} A^{(k)}
        & \sec\theta \left(A - \sin^2\theta \begin{bmatrix}
        \V0&\V0\\\V0&\V I_{N/2-k+1}
        \end{bmatrix}A^{(k)} \right)
        \end{bmatrix}.
    \end{equation}
    If $k=N/2+j$ for $j \ge 1$, then
    \begin{equation}
    \label{eq:Bk second half}
        B^{(k)} = \begin{bmatrix}
        \cos\theta A^{(N/2)} & \sin\theta A^{(N/2)}\\\V0&\sec\theta A^{(j)}
        \end{bmatrix}.
    \end{equation}
\end{lemma}

\begin{proof}
Let $A = L'U'$ be the GENP factorization of $A$. 
Then $B$ has an $LU$ factorization using GENP, $B=LU$, where $B(\t)=L_{\t}U_{\t}$, and
\begin{equation}
\label{eq: LU kron}
    L = L_\t \otimes L' = \begin{bmatrix}
        L' & \V 0\\ -\tan\t L' & L'
    \end{bmatrix} \quad \mbox{and} \quad U = U_\t \otimes U' = \begin{bmatrix}
        \cos\t U' & \sin\t U'\\\V 0 &\sec\t U'
    \end{bmatrix}
\end{equation}
by \Cref{lemma: kron_factor}. Recall 
\begin{equation}
\label{eq:def Lk}
B^{(k)} = L_{k}B^{(k-1)} = L_{k}L_{k-1}\cdots L_{1} B =: L^{(k)} B,    
\end{equation}
where $L^{(1)} := \V I$, $L^{(N)} = L^{-1}$, and for $1 \le k < N$
\begin{equation}
    L_k = \V I - \sum_{i \ge k} L_{i,k-1} \V E_{i,k-1}.
\end{equation}
It follows
\begin{align}
\label{eq:Bk_Lk}
    {L^{(k)}}^{-1}
        &= L_1^{-1}\cdots L_k^{-1}= (\V I + \sum_{i>1} L_{i1}\V E_{i1})\cdots (\V I + \sum_{i>k-1} L_{i,k-1} \V E_{i,k-1}) 
        =
        \V I + \sum_{i > j, k>j} L_{ij} \V E_{ij} \nonumber \\
        &= \left[\begin{array}{c|c}L_{:,:k-1} & \begin{array}{c} \V0\\\hline \V I_{N-k+1}\end{array} \end{array}\right].
\end{align}

If $k \le N/2$, then by \eqref{eq: LU kron} we have
\begin{align*}
    {L^{(k)}}^{-1}
        &=\begin{bmatrix}
        {L'^{(k)}}^{-1} &\V0\\ \V0& {L'^{(k)}}^{-1}
        \end{bmatrix}
        \begin{bmatrix}
        \V I_{N/2} &\V0\\
        -\tan\theta \begin{bmatrix}
        \V I_{k-1}&\V0\\\V0&\V0
        \end{bmatrix} & \V I_{N/2}
        \end{bmatrix}
        \begin{bmatrix}
        \V I_{N/2} &\V0\\ \V0& L'^{(k)}
        \end{bmatrix},
\end{align*}
and so
\begin{align*}
    L^{(k)}
        &=\begin{bmatrix}
        \V I_{N/2}&\V0 \\ \V0& {L'^{(k)}}^{-1}
        \end{bmatrix}
        \begin{bmatrix}
        \V I_{N/2}&\V0\\
        \tan\theta \begin{bmatrix}
        \V I_{k-1}&\V0\\\V0&\V0
        \end{bmatrix} & \V I_{N/2}
        \end{bmatrix}
        \begin{bmatrix}
        {L'^{(k)}} &\V0\\\V0 & {L'^{(k)}}
        \end{bmatrix}.
\end{align*}
It follows
\begin{align}
\label{eq:Bk first half_1}
    B^{(k)}
    &=L^{(k)}B =  \begin{bmatrix}
        \V I&\V0 \\\V0 & {L'^{(k)}}^{-1}
        \end{bmatrix}
        \begin{bmatrix}
        \V I&\V0\\
        \tan\theta \begin{bmatrix}
        \V I_{k-1}&\V0\\\V0&\V0
        \end{bmatrix} & \V I
        \end{bmatrix}
        \begin{bmatrix}
        {L'^{(k)}} &\V0\\\V0 & {L'^{(k)}}
        \end{bmatrix}B \nonumber \\
    &=
        \begin{bmatrix}
        \cos\theta A^{(k)}&\sin\theta A^{(k)} \\ -\sin\theta{L'^{(k)}}^{-1}\begin{bmatrix}
        \V 0&\V0\\\V0&\V I_{N/2-k+1}
        \end{bmatrix}A^{(k)} & \sec\theta  {L'^{(k)}}^{-1}\begin{bmatrix}
        \V I_{k-1}&\V0 \\\V0& \cos^2\theta \V I_{N/2-k+1}
        \end{bmatrix} A^{(k)}
        \end{bmatrix}.
\end{align}
We can write
$$
    L'^{(k)} = \begin{bmatrix}
     L_1 &\V0\\  L_2 & \V I_{N/2-k+1}
    \end{bmatrix} = \begin{bmatrix}
    L_1&\V0\\\V0&\V I_{N/2-k+1}
    \end{bmatrix}\begin{bmatrix}
    \V I_{k-1}&\V0\\L_2 & \V I_{N/2-k+1}
    \end{bmatrix}
$$
and so
$$
    {L'^{(k)}}^{-1}
        =\begin{bmatrix}
         L_1^{-1} &\V0\\ - L_2  L_1^{-1} & \V I_{N/2-k+1}
        \end{bmatrix} = \begin{bmatrix}
        \V I_{k-1} &\V0\\-L_2 & \V I_{N/2-k+1}
        \end{bmatrix}\begin{bmatrix}
        L_1^{-1}&\V0\\\V0&\V I_{N/2-k+1}
        \end{bmatrix}.
$$
First, we see
\begin{align}
\label{eq:Bk first half_0}
    {L'^{(k)}}^{-1} \begin{bmatrix}
        \V 0&\V0\\\V0&\V I_{N/2-k+1}
        \end{bmatrix} = \begin{bmatrix}\V 0&\V0\\\V0&\V I_{N/2-k+1}\end{bmatrix}.
\end{align}
Next, note
\[
    A^{(k)} = L'^{(k)}A = \begin{bmatrix}
    L_1&\V0\\L_2&\V I
    \end{bmatrix} \begin{bmatrix} A_{:k-1,:k-1} & A_{:k-1,k:}\\ A_{k:,:k-1} & A_{k:,k:}
    \end{bmatrix} = \begin{bmatrix}
    L_1A_{:k-1,:k-1} & L_2 A_{:k-1,k:} \\ \V0& L_2 A_{:k-1,k:} + A_{k:,k:}
    \end{bmatrix},
\]
where we further note 
\begin{equation}\nonumber 
    A^{(k)}_{k:,:k-1} = L_2A_{:k-1,:k-1}+A_{k:,:k-1} = \V0
\end{equation} 
since $A^{(k)}$ has zeros below the first $k-1$ diagonals. It follows
\begin{align}
\label{eq:Bk first half_2}
    {L'^{(k)}}^{-1} \begin{bmatrix}
    \V I &\V0\\\V0&\cos^2\theta \V I_{N/2-k+1}
    \end{bmatrix} A^{(k)}
    &=\begin{bmatrix}
    \V I&\V0 \\-L_2&\V I
    \end{bmatrix}\begin{bmatrix}
    L_1^{-1}&\V0 \\\V0&\V I
    \end{bmatrix}\begin{bmatrix}
    L_1A_{:k-1,:k-1} & L_2 A_{:k-1,k:} \\ \V0& \cos^2\theta A^{(k)}_{k:,k:}
    \end{bmatrix} \nonumber \\
    &= A - \sin^2\theta \begin{bmatrix}
    \V 0 &\V0\\\V0&\V I_{N/2-k+1}
    \end{bmatrix} A^{(k)},
\end{align}
using also $-L_2A_{:k-1,k:} = A_{k:,k:} - A^{(k)}_{k:,k:}$. Combining  \eqref{eq:Bk first half_1}, \eqref{eq:Bk first half_0} and \eqref{eq:Bk first half_2} then yields \eqref{eq:Bk first half}.

If $k = N/2+j$ for $j \ge 1$, then again using \eqref{eq: LU kron} we have
\begin{align*}
    {L^{(k)}}^{-1}
        &=\left[
        \begin{array}{c|c}
        L'&\V0\\ \hline 
        -\tan\theta L' & \begin{array}{c|c}
        L'_{:,1:j-1} & \begin{array}{c} \V 0 \\ \hline \V I_{N/2-j+1}
        \end{array}
        \end{array}
        \end{array}
        \right] 
        =
        \begin{bmatrix} \V I&\V0 \\ -\tan\theta \V I & \V I
        \end{bmatrix}
        \begin{bmatrix}
        L'&\V0\\\V0 & {L'^{(j)}}^{-1}
        \end{bmatrix}
\end{align*}
so that
\begin{equation}\nonumber
    L^{(k)} = 
        \begin{bmatrix}
        {L'}^{-1} &\V0\\\V0 & {L'^{(j)}}
        \end{bmatrix}
        \begin{bmatrix} \V I &\V0\\ \tan\theta \V I & \V I
        \end{bmatrix} = \begin{bmatrix}
        {L'}^{-1}&\V0\\ \tan\theta L'^{(j)}&L'^{(j)}
        \end{bmatrix}
\end{equation}
and hence
\begin{align}
    B^{(k)}
        &= L^{(k)} B = \begin{bmatrix}
        {L'}^{-1} &\V0\\ \V0& {L'^{(j)}}
        \end{bmatrix}
        \begin{bmatrix} \cos \theta A & \sin\theta A \\ \V0& \sec\theta A
        \end{bmatrix} = 
        \begin{bmatrix} 
        \cos \theta U' & \sin\theta U' \\ \V0& \sec\theta A^{(j)}
        \end{bmatrix}.
\end{align}
\eqref{eq:Bk second half} follows then by noting $U' = A^{(N/2)}$.
\end{proof}

For $B \in \B_s(N)$, one direct consequence of \Cref{lemma:Bk} and induction is
\begin{equation}
    \begin{bmatrix} \V 0 & \V I_{N-k+1}\end{bmatrix} |B^{(k)}| \begin{bmatrix} \V 0\\\V I_{N-k+1}\end{bmatrix} \V e_1  = \begin{bmatrix} \V 0 & \V I_{N-k+1}\end{bmatrix} |B^{(k)}|^T \begin{bmatrix} \V 0\\\V I_{N-k+1}\end{bmatrix} \V e_1
\end{equation}
for all $k$. In particular, there is symmetry in the leading column and row of the remaining untriangularized block of $|B^{(k)}|$. (A more involved argument can establish the entire untriangularized block of $|B^{(k)}|$ is symmetric.) This results in:

\begin{corollary}\label{cor: gerp}
Let $B \in \B_s(N)$ with $PB = LU$ the GEPP factorization of $B$. Then this is also the GERP factorization of $B$.
\end{corollary}

The remaining technical piece involves establishing the maximal growth encountered using GENP, GEPP and GERP is found in the final GE step. 

\begin{proposition}
\label{prop:max Uk}
Let $B \in \B_s(N)$ for $n\ge 1$ and  $\eta,\varepsilon \in \mathbb R$ satisfying $|\eta|,|\varepsilon|,|\eta-\varepsilon| \le 1$. Let $PBQ = LU$ be the $LU$ factorization of $B$ using GENP (with $P =Q= \V I$), GEPP (with $Q = \V I)$ or GERP. Then for all $k=1,2,\ldots,N$, 
\begin{equation}
\label{eq:max Uk ineq}
    \left\|\begin{bmatrix}
    \V 0  &\V I_{N-k+1}
    \end{bmatrix}(\eta PBQ - \varepsilon (PBQ)^{(k)})\begin{bmatrix}
    \V 0 \\ \V I_{N-k+1}
    \end{bmatrix} \right\|_{\max} \le \|U\|_{\max}.
\end{equation}
In particular,  
\begin{equation}
\label{eq:max Uk}
    \max_k \|(PBQ)^{(k)}\|_{\max} = \|U\|_{\max}.
\end{equation}
\end{proposition}
\begin{proof}
First, note it suffices to establish the GENP case since we can reduce to the case when $P = Q = \V I$: We have $PBQ = B'D$ for diagonal $D$ with diagonal entries in $\{\pm 1\}$ and $B' \in \B_s(N)$ with  $(PBQ)^{(k)} =  {B'}^{(k)}D$ using \Cref{prop: ge factors} if using GEPP along with \Cref{cor: gerp} if using GERP. 
It follows then
\begin{align*}
    &\left\|\begin{bmatrix}
    \V 0  &\V I_{N-k+1}
    \end{bmatrix}(\eta PBQ - \varepsilon (PBQ)^{(k)})\begin{bmatrix}
    \V 0 \\ \V I_{N-k+1}
    \end{bmatrix} \right\|_{\max}\\
    &\hspace{1pc}=\left\|\begin{bmatrix}
    \V 0  &\V I_{N-k+1}
    \end{bmatrix}(\eta {B'} - \varepsilon {B'}^{(k)})D\begin{bmatrix}
    \V 0 \\ \V I_{N-k+1}
    \end{bmatrix} \right\|_{\max} \\
    &\hspace{1pc}=\left\|\begin{bmatrix}
    \V 0  &\V I_{N-k+1}
    \end{bmatrix}(\eta {B'} - \varepsilon {B'}^{(k)})\begin{bmatrix}
    \V 0 \\ \V I_{N-k+1}
    \end{bmatrix} \right\|_{\max}
\end{align*}
using the fact $\|\cdot\|_{\max}$ is invariant under row or column $\pm 1$ multiplications. This establishes the GENP case subsumes the GEPP and GERP cases. We will now assume $P = Q = \V I$.

Next, note how \eqref{eq:max Uk} follows: since 
$$
B^{(k)} = \begin{bmatrix}
    U_{:k-1,:k-1} & U_{:k-1,k:}\\
    \V 0 & \begin{bmatrix}
        \V 0 & \V I_{N-k+1}
    \end{bmatrix} B^{(k)} \begin{bmatrix}
        \V 0\\\V I_{N-k+1}
    \end{bmatrix}
\end{bmatrix}
$$
then $\|B^{(k)}\|_{\max} \le \|U\|_{\max}$ using \eqref{eq:max Uk ineq} with $\eta = 0$ and $\varepsilon = -1$.

To prove \eqref{eq:max Uk ineq}, we will once again use induction on $n$. Note first the result always holds for $k = 1$: Since
\begin{equation}
\label{eq:BleU}
    \|B\|_{\max} = \prod_{j=1}^n \max(|\cos\t_j|,|\sin\t_j|) \le 1 \le \prod_{j=1}^n |\sec\t_j| = |U_{NN}| \le \|U\|_{\max}
\end{equation} 
using \Cref{prop: ge factors} and $|\eta - \varepsilon| \le 1$, then
\begin{equation}
    \|\V I_N(\eta B - \varepsilon B^{(1)})\V I_N\|_{\max} = |\eta - \varepsilon| \|B\|_{\max} \le \|U\|_{\max}.
\end{equation}
So we can consider only $k \ge 2$. For $n=1$ and $k = 2$,
$$\left\|\begin{bmatrix}
0&1
\end{bmatrix}(\eta B - \varepsilon B^{(2)})\begin{bmatrix}
0\\1
\end{bmatrix}\right\|_{\max} = |\eta \cos^2\t - \varepsilon||\sec\t| \le  |\sec\t|= \|U\|_{\max},$$
where we note 
\begin{equation}
\label{eq: eta prime ineq}
    |\eta \cos^2\t - \varepsilon| \le \max(|\varepsilon|,|\eta - \varepsilon|) \le 1.
\end{equation}
(This can be established by considering the cases for $\eta \ge 0$ and $\eta < 0$ separately.) This completes the base case.

Now assume the result holds for $\B_s(N/2)$ for $n \ge 2$. Let $B = B(\t,A) \in \B_s(N)$ with $A = L'U' \in \B_s(N/2)$. Note 
\begin{equation}
\label{eq:AleU}
\|A\|_{\max} \le \|U'\|_{\max} \le |\sec\t|\|U'\|_{\max} = \|U\|_{\max},
\end{equation}
using \eqref{eq:BleU} for the first inequality and \Cref{prop: ge factors} for the last equality.

For $k \le N/2$, then for $\V I = \V I_{N/2-k+1}$ when not indicated otherwise, we have
\begin{align*}
    &\begin{bmatrix}
    \V 0  &\V I_{N-k+1}
    \end{bmatrix}(\eta B - \varepsilon B^{(k)})\begin{bmatrix}
    \V 0 \\ \V I_{N-k+1}
    \end{bmatrix} \\
    &\hspace{1pc}= \begin{bmatrix}
    \cos\t\begin{bmatrix}
    \V 0 & \V I
    \end{bmatrix}(\eta A - \varepsilon A^{(k)})\begin{bmatrix}
    \V 0 \\ \V I
    \end{bmatrix} & \sin\t\begin{bmatrix}
    \V 0 & \V I
    \end{bmatrix}(\eta A - \varepsilon A^{(k)})\\
    -\sin\t\left(\eta A - \varepsilon \begin{bmatrix}
    \V0&\V0\\\V0&\V I
    \end{bmatrix}A^{(k)}\right)\begin{bmatrix}
    \V 0 \\ \V I
    \end{bmatrix}
    &- \sec\t\left((\varepsilon-\eta\cos^2\theta)A - \varepsilon\sin^2\t  \begin{bmatrix}
    \V0&\V0\\\V0&\V I
    \end{bmatrix}A^{(k)}\right)
    \end{bmatrix}
\end{align*}
using \Cref{lemma:Bk}. Let $\eta' = \varepsilon - \eta \cos^2\t$ and $\varepsilon' = \varepsilon \sin^2 \t$. First, note $|\eta'| \le 1$ by \cref{eq: eta prime ineq}  while also $|\varepsilon'| = \sin^2 \t |\varepsilon| \le 1$ and $|\eta ' - \varepsilon'| = \cos^2\t |\eta - \varepsilon| \le 1$. By the inductive hypothesis, we have 
\begin{align*}
    \left\|\begin{bmatrix} \V0 &\V I
\end{bmatrix}(\eta A - \varepsilon A^{(k)}) \begin{bmatrix}
\V0 \\\V I
\end{bmatrix} \right\|_{\max} &\le \|U'\|_{\max} \le \|U\|_{\max} \quad \mbox{and}\\
|\sec\t|\left\|\begin{bmatrix}
    \V0&\V I
    \end{bmatrix}\left(\eta ' A - \varepsilon'  A^{(k)}\right)\begin{bmatrix}
    \V0\\\V I
    \end{bmatrix} \right\|_{\max} &\le |\sec\t|\|U'\|_{\max} = \|U\|_{\max}.
\end{align*}
Moreover, $|\eta\sin\t|\|A\|_{\max} \le \|U'\|_{\max} \le \|U\|_{\max}$ and 
\begin{equation}
|\sec\t||\varepsilon - \eta \cos^2\t|\|A\|_{\max} \le |\sec\t|\|U'\|_{\max} = \|U\|_{\max}    
\end{equation}
by \eqref{eq:AleU}. It follows
\begin{align}
&\left\|\begin{bmatrix}
    \V 0  &\V I_{N-k+1}
    \end{bmatrix}(\eta B - \varepsilon B^{(k)})\begin{bmatrix}
    \V 0 \\ \V I_{N-k+1}
    \end{bmatrix}\right\|_{\max}\\
    &\hspace{1pc} =\max\left(
    \begin{array}{c}|\eta \sin\theta|\|A_{N/2-k+1:,:k-1}\|_{\max},\\|\eta \sin\theta|\| A_{:k-1,N/2-k+1:}\|_{\max},\\
    \max(|\cos\t|,|\sin\t|)\left\|\begin{bmatrix} \V0 &\V I
\end{bmatrix}(\eta A - \varepsilon A^{(k)}) \begin{bmatrix}
\V0 \\\V I
\end{bmatrix} \right\|_{\max},\\
|\sec\t||\varepsilon - \eta \cos^2\theta |\| A\|_{\max},\\
|\sec\t|\left\|\begin{bmatrix}
    \V0&\V I
    \end{bmatrix}\left((\varepsilon-\cos^2\theta \eta)A - \varepsilon\sin^2\t  A^{(k)}\right)\begin{bmatrix}
    \V0\\\V I
    \end{bmatrix} \right\|_{\max}
\end{array}\right)\\
&\hspace{1pc}\le \|U\|_{\max}.
\end{align}

For $k = N/2+j$ and $j \ge 1$ so that $N-k+1=N/2-j+1$, let $\eta ' =\eta \cos^2\t$ and $\varepsilon' = \varepsilon$ so then $|\eta'|,|\varepsilon'|,|\eta'-\varepsilon'| \le 1$. Hence,
\begin{align}
&\left\|\begin{bmatrix}
    \V 0  &\V I_{N-k+1}
    \end{bmatrix}(\eta B - \varepsilon B^{(k)})\begin{bmatrix}
    \V 0 \\ \V I_{N-k+1}
    \end{bmatrix}\right\|_{\max}\\
    &\hspace{1pc}= |\sec\t|\left\|\begin{bmatrix}
    \V 0  &\V I_{N/2-j+1}
    \end{bmatrix}(\eta \cos^2\t A - \varepsilon  A^{(j)})\begin{bmatrix}
    \V 0 \\ \V I_{N/2-j+1}
    \end{bmatrix} \right\|_{\max}\\
    &\hspace{1pc}\le |\sec\t|\|U'\|_{\max} = \|U\|_{\max}
\end{align}
follows by the inductive hypothesis.
\end{proof}

We can now establish the multiplicativity of each of the Haar-butterfly matrices growth factors with respect to the Kronecker factors.
\begin{proposition}
\label{lemma: gf mult}
If $B=B(\boldsymbol \t) \in \B_s(N)$, then
\begin{equation}
    \kappa_\infty(B) = \prod_{j=1}^n \kappa_\infty(B(\t_j)).
\end{equation}
Using GENP, GEPP, or GERP, then
\begin{align}
    \rho(B) &= \prod_{j=1}^n \rho(B(\t_j))\\
    \rho_o(B) &= \prod_{j=1}^n \rho_o(B(\t_j))\\
    \rho_\infty(B) &= \prod_{j=1}^n \rho_\infty(B(\t_j)).
\end{align}
\end{proposition}
\begin{proof}
Let $PBQ = LU$ be the $LU$ factorization of $B$ using GENP, GEPP, or GERP. \Cref{cor: gerp} yields the GEPP and GERP factorizations align for $B$. By \Cref{prop:max Uk}, then 
\begin{equation}
    \rho(B) = \frac{\|L\|_{\max} \cdot \displaystyle\max_k \|B^{(k)}\|_{\max}}{\|B\|_{\max}}=\frac{\|L\|_{\max} \|U\|_{\max}}{\|B\|_{\max}}.
\end{equation}
Since $B^{-1} = \bigotimes_{j=1}^n B(\t_{n-j+1})^{-1}$ and $|\bigotimes_{j=1}^n A_j| = \bigotimes_{j=1}^n |A_j|$, then the result follows from the mixed-product property and \Cref{lemma: kron_factor,lemma: norm_mult}.
\end{proof}

It thus remains only to establish the case for $n = 1$. Note for $B=B(\t) \in \B_s(2) = \B(2) = \SO(2)$ where $\cos\t \ne 0$ and $\sin \t \ne 0$, we have
\begin{align}
    B =L_\t U_\t 
    \quad \mbox{and} \quad
    PB =   B\left(\frac\pi2-\t\right)(-1 \oplus 1)
\end{align}
for $P = P_{(1 \ 2)}$ so that $|PB| = |L_{\frac\pi2-\t} U_{ \frac\pi2-\t}|$, while also
\begin{equation}
|L_\t||U_\t| = \begin{bmatrix}
|\cos\t| & |\sin\t|\\
|\sin\t| & |\cos\t|(1 + 2\tan^2\t)
\end{bmatrix}
\end{equation}
It follows directly
\begin{lemma}
\label{lemma:gf_2x2}
    Let $B = B(\theta) \in \B_s(2)$. Let $f(\t) = |\tan \t|$ if using GENP when $\cos\t \ne 0$ and $f(\t) = \min(|\tan\t|,|\cot\t|)$ if using GEPP. Then
    \begin{align}
    \rho(B) &= 1+f(\t)^2\\
        \rho_o(B) &= 1+\frac{2f(\t)^2}{1+f(\t)}\\
    \rho_\infty(B) &= 1+\max(f(\t),f(\t)^2)
\end{align}
for all $\t$. Moreover, using GENP then $1 \le \rho(B) \le \rho_\infty(B)$ for all $\t$, with $1=\rho(B)=\rho_o(B) = \rho_\infty(B)$ for $\cos\t = 0$, $\rho(B) = \rho_\infty(B)$ if $|\tan\t| \ge 1$, and $1<\rho(B) < \rho_\infty(B)$ otherwise, while also $\rho(B) < \rho_o(B)$ for $|\tan\t| \in (0,1)$. Using GEPP, then $1 \le \rho(B) \le \rho_o(B) \le \rho_\infty(B) \le 2$, with $1 = \rho(B) = \rho_o(B) = \rho_\infty(B)$ for $\cos\t = 0$ or $\sin \t = 0$, $\rho(B) = \rho_o(B) = \rho_\infty(B) = 2$ for $|\tan \t| = 1$, and strict inequalities otherwise.
\end{lemma}

\begin{proof}
    The GENP case is straightforward using the Pythagorean identity $\sec^2\theta = 1+\tan^2\theta$, where we note $\|B\|_{\max} \le  \|U_\t\|_{\max}$ while also $\|B\|_{\max} = |\cos\t| \|L_\t\|_{\max}$ and $\|B\|_\infty = |\cos\t|\|L_\t\|_\infty$. For GEPP, note a pivot occurs only if $|\tan\t|>1$, in which case $|\tan \t|$ would be replaced with $|\tan(\frac\pi2-\t)| = |\cot\t|$, and hence in general by $\min(|\tan\t|,|\cot\t|)$, for the computations of the growth factors from the GENP case. The remaining inequalities follow from  $\frac{2x^2}{1+x} \le x\max(1,x)$ for $x \ge 0$, where equality only holds for $x = 0$ or $x = 1$ (note this is equivalent to the trivial bound $\rho_o \le \rho_\infty$ in the GENP case).
\end{proof}

The structure of simple scalar butterfly matrices enables us to explicitly derive more properties of these matrices other than the growth factor. For instance, the following lemma shows a straightforward computation of the $\infty$-condition number:
\begin{lemma}
\label{lemma: cond_2x2}
Let $B = B(\t) \in \B_s(2)$. Then
\begin{equation}
    \kappa_\infty(B) = 1 + |\sin(2\t)|.
\end{equation}
\end{lemma}
\begin{proof}
Since $\|B(\t)^{-1}\|_\infty = \|B(-\t)\|_\infty = \|B(\t)\|_\infty$, we can compute directly $$
\kappa_\infty(B) = \|B\|_\infty \|B^{-1}\|_\infty = \|B\|^2_\infty = (|\cos\t|+|\sin\t|)^2 =  1+|\sin(2\t)|.
$$
\end{proof}

The remaining pieces relate to straightforward properties of $\t \sim \Uniform([0,2\pi))$. First, note the following simple lemma, which is used sporadically throughout the following arguments:
\begin{lemma}
\label{lemma:simple}
If $X \sim \Uniform(0,1)$ and $Y$ is independent of $X$, then $(X+Y)\pmod{1} \sim X$.
\end{lemma}
\begin{proof}
Since $X+Y\mid Y \sim \Uniform(Y,Y+1)\mid Y$ and so $(X+Y)\pmod{1} \mid Y \sim \Uniform(0,1)$, then for $t \in (0,1)$,
\begin{equation}\nonumber
    \P((X+Y)\hspace{-8pt}\pmod{1} \le t) = \E\P((X+Y)\hspace{-8pt}\pmod{1} \le t \mid Y) = \E \P(X \le t) = \P(X \le t).
\end{equation}
\end{proof}

Next, recall for $X \sim \Cauchy(1)$, then
\begin{equation}
    \P(X \le t) = \frac1\pi \arctan t + \frac12.
\end{equation}
Note $X \sim -X$ (say, since $\arctan x$ is an odd function), and so
\begin{equation}
\label{eq: cdf abs cauchy}
    \P(|X| \le t) = 2\P(X \le t) - 1 = \frac2\pi \arctan t.
\end{equation}
In particular, $\P(|X| \le 1) = \frac12$. We see:

\begin{lemma}
\label{lemma: tan cauchy}
If $\t \sim \Uniform([0,2\pi))$, then $\tan\t,\cot\t \sim \Cauchy(1)$.
\end{lemma}
\begin{proof}
Recall if $Y \sim \Uniform(0,1)$, then $\tan(\pi(Y-\frac12)) \sim \Cauchy(1)$. Note $\pi(Y-\frac12)\pmod \pi \sim \pi Y \sim \Uniform(0,\pi)$ by \Cref{lemma:simple}. Since $\t \pmod \pi \sim \pi Y$ and by the periodicity of $\tan x$ also $\tan x = \tan (x \pmod \pi)$, then
\begin{equation}
    \tan \t  \sim \tan (\pi Y) \sim  \tan\left(\pi\left(Y-\frac12\right)\right) \sim \Cauchy(1).
\end{equation}
Hence, we have also $\cot\t \sim \cot(\t-\frac\pi2) = -\tan\t  \sim \Cauchy(1) $ by \Cref{lemma:simple}.
\end{proof}


\begin{lemma}
\label{lemma: tan cutoff cauchy}
If $\t \sim \Uniform([0,2\pi))$ and $X \sim \Cauchy(1)$, then 
\begin{equation}
\min(|\tan\t|,|\cot\t|) \sim |X|\mid {|X| \le 1}.
\end{equation}
\end{lemma}
\begin{proof}
Using \Cref{lemma: tan cauchy}, then for $t \in (0,1]$, we have
\begin{align*}
    &\P(\min(|\tan\t|,|\cot \t|) \le t) 
    =   1-\P\left(\min\left(|X|,\frac1{|X|}\right) \ge t\right)= 1 - \P\left(|X| \ge t, \frac1{|X|} \ge t\right)\\
    &\hspace{0.9pc} 
    =1+\P(|X| \le t)-\P\left(|X| \le \frac1t\right)
    =1+\frac2\pi\left(\arctan t - \arctan \frac1t\right) 
    =1+\frac2\pi\left(2\arctan t - \frac\pi2\right)\\
    &\hspace{0.9pc}=\frac4\pi \arctan t = \frac{\P(|X| \le t)}{\P(|X| \le 1)}
    =\P(|X| \le t \mid |X| \le 1)
\end{align*}
using also the fact $\arctan t$ and $\arctan \frac1t$ comprise complementary angles when $t >0$.
\end{proof}

The $n=1$ growth factor case now follows:
\begin{lemma}
\label{prop:gf_rbm_2x2}
    Let $B \sim \B_s(2,\Sigma_S)$ and $X \sim \operatorname{Cauchy}(1)$. If using GENP, let $Y = |X|$ and if using GEPP, let $Y = |X| \mid |X| \le 1$. Then
    \begin{align}
        \rho(B) &\sim 1+Y^2\\
        \rho_o(B) & \sim 1 + \frac{2Y^2}{1 + Y}\\
        \rho_\infty(B) &\sim 1+\max(Y,Y^2).
    \end{align}
\end{lemma}
\begin{proof}
This follows directly from \Cref{lemma:gf_2x2,
lemma: tan cauchy,lemma: tan cutoff cauchy}
\end{proof}

Additionally, we can explicitly derive the distribution of the $\infty$-condition numbers of random butterfly matrices. First recall if $Y \sim \operatorname{Arcsine}(0,1)$ then for $t \in [0,1]$ we have
\begin{equation}
    \P(Y \le t) = \frac2\pi \arcsin \sqrt t.
\end{equation}

\begin{lemma}
\label{lemma: arcsine}
If $\t \sim \Uniform([0,2\pi))$ then $\sin^2(2\t) \sim \operatorname{Arcsine}(0,1).$
\end{lemma}
\begin{proof}
Let $\varphi \sim \Uniform([0,\pi))$. Note first $2\t\pmod \pi \sim \t \pmod \pi \sim \varphi$, so since $|\sin(x)|$ has period $\pi$ then $|\sin(2\t)| \sim \sin \varphi$. Now note for $t \in [0,1]$
\begin{align*}
    \P(|\sin2\t| \le t) 
    &= \P(\sin\varphi \le t) = \P(\varphi \in [0,\arcsin t] \cup [\pi-\arcsin t,\pi]) 
    = \frac2\pi \arcsin t.
\end{align*}
Hence, $\displaystyle \P(\sin^2(2\t) \le t) = \P(|\sin (2\t)| \le \sqrt t) = \frac2\pi \arcsin \sqrt t$.
\end{proof}
We have then
\begin{lemma}
    \label{prop:cond_2x2}
    If $B \sim \B_s(2,\Sigma_S)$ and $Y \sim \operatorname{Arcsine}(0,1)$, then $$\kappa_\infty(B) \sim 1+\sqrt Y.$$
\end{lemma}
\begin{proof}
Use \Cref{lemma: cond_2x2,lemma: arcsine}.
\end{proof}

Lastly, we will state a theorem needed to establish the minimality of the GEPP max-norm growth factor as given at the end of \Cref{thm:gf}:
\begin{theorem}[\cite{HiHi89}]
\label{thm: lower bd gf}
If $A \in \mathbb R^{n\times n}$ is nonsingular, then $ \|A\|_{\max} \|A^{-1}\|_{\max} \ge \frac1n$ and for any permutation matrices $P,Q$ such that $PAQ$ has a GENP factorization of the form $PAQ=LU$, then $\rho(PAQ) \ge \frac{\|L\|_{\max}}{\|A\|_{\max} \|A^{-1}\|_{\max}}$ for the GENP growth factor.
\end{theorem}
Note this theorem was also used to establish \cref{eq: haar lb}.

Now we can sum up these results to establish the main statements from \Cref{subsec:gf_haar_b}.
\begin{proof}[Proof of \Cref{thm:gf}]
Use \Cref{lemma: gf mult,lemma:gf_2x2,prop:gf_rbm_2x2} with the uniqueness results in the random models following from \Cref{thm: unique genp,thm:GEPP_unique}. The last statement that the GEPP growth factor minimizes the growth factor among all pivoting strategies follows by using \Cref{lemma:gf_2x2} to see 
\begin{equation}
\|U_\t\|_{\max}^2=\rho(B(\t))=\frac{\|L_\t\|_{\max} \|U_\t\|_{\max}}{\|B(\t)\|_{\max}} = \frac{\|U_\t\|_{\max}}{\|B(\t)\|_{\max}}    
\end{equation}
so that $\|B(\t)\|_{\max} = \|B(\t)^{-1}\|_{\max} = \|U_\t\|_{\max}^{-1}$; \Cref{lemma: kron_factor,prop:max Uk} yield then $\|B\|_{\max} = \|B^{-1}\|_{\max} = \|U\|_{\max}^{-1} = \rho^{\operatorname{GEPP}}(B)^{-1/2}$ for $PB=LU$ the GEPP factorization. Hence, if $PBQ=L_oU_o$ is a GENP factorization, then 
\begin{equation}
\rho^{\operatorname{GENP}}(PBQ) \ge \frac{\|L_o\|_{\max}}{\|B\|_{\max} \|B^{-1}\|_{\max}} = \|L_o\|_{\max} \cdot \rho^{\operatorname{GEPP}}(B) \ge \rho^{\operatorname{GEPP}}(B)    
\end{equation}
by \Cref{thm: lower bd gf}.
\end{proof}


\begin{proof}[Proof of \Cref{prop: cond infty}]
Use \Cref{lemma: gf mult,prop:cond_2x2}
\end{proof}

\bibliographystyle{siamplain}
\bibliography{references}

\begin{thebibliography}{10}

\bibitem{cnn}
{\sc K.~Alizadeh, A.~Pradhu, A.~Farhadi, and M.~Rastegari}, {\em Butterfly
  transform: An efficient {FFT} based neural architecture design}, Proc. of the
  Conf. on CV and Pat. Recog.,  (2020),
  \url{https://doi.org/10.1109/CVPR42600.2020.01204}.

\bibitem{baboulin}
{\sc M.~Baboulin, X.~S. Li, and F.-H. Rouet}, {\em Using random butterfly
  transformations to avoid pivoting in sparse direct methods}, In: Proc. of
  Int. Con. on Vector and Par. Proc.,  (2014),
  \url{https://doi.org/10.1007/978-3-319-17353-5_12}.

\bibitem{CoPe07}
{\sc V.~Cort\'es and J.~Pe{\~n}a}, {\em Growth factor and expected growth
  factor of some pivoting strategies}, Journal of Comp. and App. Math., 202
  (2007), pp.~292--303, \url{https://doi.org/10.1016/j.cam.2006.02.040}.

\bibitem{Cryer}
{\sc C.~W. Cryer}, {\em Pivot size in {G}aussian elimination}, Numer. Math., 12
  (1968), pp.~335--345, \url{https://doi.org/10.1007/BF02162514},
  \url{https://doi.org/10.1007/BF02162514}.

\bibitem{fast_alg}
{\sc T.~Dao, A.~Gu, M.~Eichhorn, A.~Rudra, and C.~Re}, {\em Learning fast
  algorithms for linear transforms using butterfly factorizations}, Proc. Mach.
  Learn. Res.,  (2019), pp.~1517--1527.

\bibitem{perfectShuffle}
{\sc P.~Diaconis, R.~L. Graham, and W.~M. Kantor}, {\em The mathematics of
  perfect shuffles}, Adv. in App. Math., 4 (1983), pp.~175--196,
  \url{https://doi.org/10.1016/0196-8858(83)90009-X}.

\bibitem{EdelmanGECP}
{\sc A.~Edelman}, {\em The complete pivoting conjecture for gaussian
  elimination is false}, The Mathematica Journal, 2 (1992), pp.~58--61.

\bibitem{Hi02}
{\sc N.~J. Higham}, {\em Accuracy and Stability of Numerical Algorithms, Second
  Edition}, SIAM, Philadelphia, PA, 2002.

\bibitem{HiHi89}
{\sc N.~J. Higham and D.~J. Higham}, {\em Large growth factors in gaussian
  elimination with pivoting}, SIAM J. Matrix Anal. Appl., 10 (1989),
  pp.~155--164, \url{https://doi.org/10.1137/0610012}.

\bibitem{HiHi20}
{\sc N.~J. Higham, D.~J. Higham, and S.~Pranesh}, {\em Random matrices
  generating large growth in {LU} factorization with pivoting}, SIAM J. Matrix
  Anal. Appl., 42 (2021), pp.~185--201,
  \url{https://doi.org/10.1137/20M1338149}.

\bibitem{LaFa72}
{\sc P.~Lancaster and H.~K. Farahat}, {\em Norms on direct sums and tensor
  products}, Math. of Comp., 26 (1972), pp.~401--414,
  \url{https://doi.org/10.2307/2005167}.

\bibitem{butterflynet}
{\sc Y.~Li, X.~Cheng, and J.~Lu}, {\em Butterfly-net: Optimal function
  representation based on convolutional neural networks}, Commun. Comput.
  Phys., 28 (2020), pp.~1838--1885,
  \url{https://doi.org/10.4208/cicp.OA-2020-0214}.

\bibitem{LiLuDo20}
{\sc N.~Lindquist, P.~Luszczek, and J.~Dongarra}, {\em Replacing pivoting in
  distributed gaussian elimination with randomized techniques}, 2020 IEEE/ACM
  ScalA,  (2020), \url{https://doi.org/10.1109/ScalA51936.2020.00010}.

\bibitem{MaTr20}
{\sc P.-G. Martinsson and J.~A. Tropp}, {\em Randomized numerical linear
  algebra: Foundations and algorithms}, Acta Numerica, 29 (2020), pp.~403--572,
  \url{https://doi.org/10.1017/S0962492920000021}.

\bibitem{Mezz}
{\sc F.~Mezzadri}, {\em How to generate random matrices from the classical
  compact groups}, Notices of the American Mathematical Society, 54 (2007),
  pp.~592 -- 604.

\bibitem{Pa95}
{\sc D.~S. Parker}, {\em Random butterfly transformations with applications in
  computational linear algebra}, Tech. rep., UCLA,  (1995).

\bibitem{PoNe00}
{\sc G.~Poole and L.~Neal}, {\em The rook's pivoting strategy}, Journal of
  Comp. and App. Math., 123 (2000), pp.~353--369,
  \url{https://doi.org/10.1016/S0377-0427(00)00406-4}.

\bibitem{stewart}
{\sc G.~W. Stewart}, {\em The efficient generation of random orthogonal
  matrices with an application to condition estimators}, SIAM J. Numer. Anal.,
  17 (1980), pp.~403--409, \url{https://doi.org/10.1137/0717034}.

\bibitem{St99}
{\sc G.~Strang}, {\em The discrete cosine transform}, SIAM Review, 41 (1999),
  pp.~135--147, \url{https://doi.org/10.1137/S0036144598336745}.

\bibitem{TrSc90}
{\sc L.~Trefethen and R.~Schreiber}, {\em Average case stability of gaussian
  elimination}, SIAM J. Matrix Anal. Appl., 11 (1990), pp.~335--360,
  \url{https://doi.org/10.1137/0611023}.

\bibitem{Tr19}
{\sc T.~Trogdon}, {\em On spectral and numerical properties of random butterfly
  matrices}, Applied Math. Letters, 95 (2019), pp.~48--58,
  \url{https://doi.org/10.1016/j.aml.2019.03.024}.

\bibitem{Tr11}
{\sc J.~A. Tropp}, {\em Improved analysis of the subsampled randomized hadamard
  transform}, Adv. Adapt. Data Anal., 3 (2011), pp.~115--126,
  \url{https://doi.org/10.1142/S1793536911000787}.

\bibitem{We40}
{\sc A.~Weil}, {\em L'int\'egration dans les groupes topologiques et ses
  applications, \emph{Actualit\'es Scientifiques et Industrielles}}, vol.~869,
  Paris: Hermann, 1940.

\bibitem{Wi61}
{\sc J.~Wilkinson}, {\em Error analysis of direct methods of matrix inversion},
  J. Assoc. Comput. Mach., 8 (1961), pp.~281--330,
  \url{https://doi.org/10.1145/321075.321076}.

\bibitem{Wi65}
{\sc J.~Wilkinson}, {\em The Algebraic Eigenvalue Problem}, Oxford University
  Press, London, UK, 1965.

\end{thebibliography}

\end{document}